\RequirePackage{fix-cm}
\documentclass[smallextended]{svjour3}       

\smartqed  
\usepackage[utf8]{inputenc}
\usepackage[T1]{fontenc}

\usepackage{indentfirst}
\usepackage[english]{babel}

\usepackage{array}
\usepackage{color} 

\usepackage{graphicx} 
\usepackage{amsmath}
\usepackage{amssymb}
\usepackage{amsfonts}	
\usepackage{moreverb}

\usepackage{tipa}
\usepackage{upgreek}

\usepackage{textcomp}
\usepackage[dvipsnames]{xcolor}
\usepackage{mathtools}
\usepackage{setspace}

\usepackage{caption}
\usepackage{subcaption}
\usepackage{booktabs}
\usepackage{float}
\allowdisplaybreaks

\usepackage{hyperref} 
\hypersetup{
	colorlinks=true,                          
	linkcolor=blue, 
	citecolor=red, 
	urlcolor=blue, 
	pdftitle={A semi-implicit thermodynamically compatible finite volume scheme for the compressible Euler equations},
	pdfauthor={Michael Dumbser, Gabriella Puppo, Sara Rinaldi, Andrea Thomann},
	pdfsubject={Journal of Scientific Computing},
	pdfkeywords={hyperbolic and thermodynamically compatible (HTC) scheme; asymptotic-preserving (AP) semi-implicit (SI) scheme; predictor-corrector scheme; low-Mach-number limit; compressible Euler equations in entropy formulation} }

\journalname{Journal of Scientific Computing}
\begin{document} 
\title{A thermodynamically compatible semi-implicit finite volume scheme for the compressible Euler equations}  

\titlerunning{A thermodynamically compatible semi-implicit FV scheme for the Euler equations}        

\author{Michael Dumbser \and Gabriella Puppo \and \\ Sara Rinaldi $^*$ \and Andrea Thomann}


\institute{
		M. Dumbser \at
	Laboratory of Applied Mathematics, DICAM, University of Trento, Via Mesiano 77, 38123 Trento, Italy\\  
	Department of Mathematics, Southern University of Science and Technology, \\ Xueyuan Avenue 1088, 518055 Shenzhen, Guangdong, China \\ \email{michael.dumbser@unitn.it}  \\
	G. Puppo \at
	Department of Mathematics, La Sapienza Universit{\`a} di Roma, Piazza Aldo Moro 5, 00185 Roma, Italy\\
	\email{gabriella.puppo@uniroma1.it} \\
	S. Rinaldi $^*$ \at
	Laboratory of Applied Mathematics, DICAM, University of Trento, Via Mesiano 77, 38123 Trento, Italy\\ \email{sara.rinaldi-1@unitn.it}\\
	A. Thomann    \at
	Universit\'e de Strasbourg, CNRS, Inria, IRMA, Strasbourg F-67000, France \\ 
	Laboratory of Applied Mathematics, DICAM, University of Trento, Via Mesiano 77, 38123 Trento, Italy
	\email{andrea.thomann@inria.fr}              
}

\date{Received: date / Accepted: date}

\maketitle

		
		\begin{abstract}
			In this work we construct a new semi-implicit (SI) hyperbolic and thermodynamically compatible (HTC) scheme for the compressible Euler equations of gasdynamics. The challenge therein lies in solving the entropy inequality as the primary evolution equation instead of the total energy conservation law, which in our framework is an additional conservation law that needs to be fulfilled by the numerical scheme as the consequence of a compatible discretization. 
			The proposed method has a predictor-corrector character, where we first solve the isentropic problem without considering the entropy production. Therein, the entropy and the nonlinear convective terms in the momentum flux are discretized explicitly, while the mass flux and the pressure terms are taken implicitly.
			The discrete momentum equation is then inserted into the discrete mass conservation equation, leading to a mildly nonlinear system for the pressure, which can be efficiently solved with a nested Newton-type algorithm. Based on this prediction, a non-negative production term in the entropy is added, ensuring a cell entropy inequality for the numerical scheme.
			To guarantee global conservation of total energy at the discrete level, a new modified global Abgrall-type flux correction is introduced.
			The presented numerical method allows time steps that are restricted only by the material wave speed and not by the sound speed. This makes the scheme particularly useful for flows in the low-Mach-number regime. 
			Numerical results are shown in one and two space dimensions to assess the theoretical properties of the new scheme.
		\end{abstract}
		

\section{Introduction} 
In 1961, Godunov \cite{God1961} established, for the first time, a rigorous mathematical connection between symmetric hyperbolicity in the sense of Friedrichs \cite{Friedrichs1958} and thermodynamic compatibility.
In \cite{God1961} he showed that the total energy conservation law is the consequence of a hyperbolic  system of partial differential equations (PDE) which has an underlying variational formulation. By defining the so-called main field \cite{Ruggeri81} or thermodynamic dual variables, which are the partial derivatives of the total energy potential with respect to the conserved variables, the total energy conservation law can be obtained via a suitable linear combination of all other equations with the aid of the thermodynamic dual variables.

This structure has since been recognized as a unifying framework, now known as the family of hyperbolic and thermodynamically compatible (HTC) or even stronger symmetric hyperbolic and thermodynamically compatible (SHTC) systems \cite{God1961}. A wide class of continuum models falls into this family, including magnetohydrodynamics (MHD)  \cite{God1972MHD,GodunovRomenski72,Godunov:1995a,Godunov:2003a,Godunov2012}, nonlinear hyperelasticity \cite{GPRmodel,Rom1998}, compressible multiphase flows \cite{RomenskiTwoPhase2010}, the turbulent shallow water equations \cite{Gavrilyuk2018,SWETurbulence}, relativistic gasdynamics \cite{GRGPR} and a hyperbolic relaxation model for the compressible Navier-Stokes equations \cite{PeshkovRomenski2014}. 
In \cite{SHTC-GENERIC-CMAT}, a connection between HTC systems and Hamiltonian mechanics was established, while an extension to continuum mechanics with torsion was presented in \cite{ilyapeshkov_2019_continuum}.
In HTC systems, the entropy-density equation, endowed with a non-negative production term in accordance with the second law of thermodynamics, is typically taken as a primary evolution equation. 

In numerical analysis, in recent years, the development of numerical schemes that are also thermodynamically compatible at the discrete level has been of great interest. 
In this way, a discrete extra conservation law is a consequence of a compatible discretization of the governing PDE system. 
Usually, when deriving a numerical scheme for hyperbolic conservation laws, the direct discretization of the total energy is favored over the entropy. 
The entropy inequality is then enforced as a consequence by means of entropy stable fluxes as in the seminal work of Tadmor \cite{Tadmor1,Tadmor2003}.
 Many other contributions have been made to the development of high-order entropy-stable schemes, for instance, in the papers of Mishra \textit{et al.} \cite{CastroFjordholm,FjordholmMishraTadmor,Hiltebrand2014}, Gassner \textit{et al.} \cite{GassnerEntropyGLM,GassnerSWE,GassnerEntropyEuler,Rueda2021,GassnerEntropyALE} and Shu
\textit{et al.} \cite{ShuEntropyMHD1,ShuEntropyMHD2}. 
In \cite{Abgrall2018,AbgrallBT2018,AbgrallOeffnerRanocha} Abgrall and collaborators presented a general framework that allows the construction of thermodynamically compatible schemes for overdetermined hyperbolic PDE systems by discretizing the total energy. 
In recent years a dual approach to existing entropy-stable schemes has been developed in a series of papers \cite{HTCAbgrall,HTCAbgrall2,HTCMHD,SWETurbulence,HTCGPR}, in the following denoted by HTC schemes, referring to a discrete version of the underlying class of HTC systems. The key feature is that HTC schemes discretize the \textit{entropy inequality directly}, while the \textit{total energy conservation} law is obtained as a \textit{consequence}. In the above references on HTC schemes the construction of space semi-discrete and fully discrete methods is addressed. 
For explicit HTC schemes applied to the compressible Euler equations, a rigorous convergence analysis in the framework of dissipative weak solutions has recently been established in \cite{Dumbser2026}.

In this work, we extend the framework of HTC schemes to semi-implicit finite-volume discretizations  applied to the compressible Euler equations of gasdynamics, which is of particular interest for low-Mach-number flows. 
The Mach number is defined as $\mathrm{Ma}=\lVert \mathbf{u}\rVert/c$, where $\mathbf{u}$ denotes the velocity vector of the fluid and $c$ the speed of sound. It measures the relative importance of compressibility effects in the flow. In the low-Mach-number regime ($\mathrm{Ma} \ll 1$), the characteristic speeds of the compressible Euler equations are dominated by the acoustic component $c$, which scales as $\mathcal{O}(\mathrm{Ma}^{-1})$ relative to the advective speed $\lVert \mathbf{u}\rVert$. 
This disparity between the fast acoustic waves and the slow material transport is the source of the severe time-step restriction faced by explicit schemes under a classical CFL condition, and it is also responsible for the loss of accuracy of standard density-based Godunov-type schemes as $\mathrm{Ma}$ goes to zero, see e.g.~\cite{DELLACHERIE2010978,GUILLARD2004655,GUILLARD199963,KleinMach}.
This observation motivates the splitting of the flux into an advective part and an acoustic part, so that only the former is constrained by the CFL condition, while the latter is treated implicitly.
The first semi-implicit pressure-based method for the compressible Euler equations was proposed by Casulli and Greenspan \cite{CasulliGreenspan1984}, but that scheme was neither fully conservative, nor thermodynamically compatible. 
In the meantime a variety of strategies have been proposed in the literature to cope with the low-Mach-number regime, ranging from flux-preconditioning techniques \cite{GUILLARD199963,TURKEL1987277} to semi-implicit \cite{KleinMach,klein,munzMPV,MunzPark,DegondTang2011} and IMEX time discretizations \cite{ThomannZenk2019,THOMANN2020109723,BoscheriAllMach,AnandanLukacova2025}, the latter being designed to be asymptotic preserving (AP) in the sense that, in the limit as $\mathrm{Ma}$  goes to zero, they recover a consistent discretization of the incompressible Euler equations. Among these, the implicit treatment of acoustic waves can substantially relax the severe time-step restriction of explicit schemes in such regimes. 

However, the incorporation of a semi-implicit splitting into the discrete HTC framework is not straightforward, since it introduces additional challenges in the consistent evaluation of the terms responsible for preserving thermodynamic compatibility and thus total energy conservation.
For the core of the low-Mach-number scheme, we follow the approach introduced in \cite{DumbserCasulli2016}. They proposed a conservative semi-implicit finite volume scheme for the compressible Euler and Navier-Stokes equations with general equation of state in which the total energy conservation law is discretized and the discrete equations yield a mildly nonlinear pressure system with diagonal nonlinearity. 
This approach was subsequently extended to staggered discontinuous Galerkin discretizations \cite{TavelliDumbser2017} and hybrid finite volume/finite element formulations on unstructured meshes in \cite{Hybrid2}. 
A key property of these schemes is that of being asymptotic-preserving, i.e.\ ensuring that the discretization of the compressible model reduces to a consistent discretization of the incompressible Euler equations as the Mach number tends to zero.  

Motivated by these considerations, the present work aims at extending the HTC framework presented in \cite{HTCAbgrall,HTCAbgrall2} to pressure-based semi-implicit finite-volume schemes in order to recover an HTC scheme that is also well suited for low-Mach-number applications. 
The main novelties of this paper are the following. 

First, a semi-implicit discretization is applied to the compressible Euler equations with the entropy as primary evolution variable and not the total energy. This considerably simplifies the structure of the final pressure system to be solved, compared to semi-implicit discretizations of the Euler equations that are based on the total energy as primary evolution variable. This requires the terms responsible for thermodynamic compatibility to be evaluated consistently across the explicit and the implicit stage and, in particular, the introduction of an entropy production associated with the implicit time discretization. Second, the local cellwise flux correction of the explicit HTC schemes cannot be enforced implicitly without solving a coupled nonlinear system with one unknown per cell interface, and it is therefore replaced by a single global Abgrall-type flux correction, which conserves the total energy globally and not locally. Third, a pressure splitting is introduced that treats only a fraction of the pressure terms implicitly, so that accuracy and stability at larger CFL numbers are retained also in the presence of strong shocks, while the low-Mach-number limit remains unaffected. In the following, the new semi-implicit HTC scheme presented in this paper will be denoted by SIHTC scheme and its pressure-split variant by SIHTC p-split.

The structure of this paper is as follows. In Section~\ref{sec: PDEHTC} we recall the general framework of hyperbolic and thermodynamically compatible systems of PDE and give the HTC formulation of the compressible Euler equations. To make the paper self-contained, section~\ref{sec: semidiscHTC} briefly recalls the semi-discrete HTC finite volume framework in one and two space dimensions. In Section~\ref{sec: Numericalscheme} we introduce our new fully discrete semi-implicit HTC scheme, which discretizes the entropy inequality and conserves the total energy globally as a consequence, with the one-dimensional formulation given in detail in Section~\ref{sec: Numscheme1d} and the p-split version in Section~\ref{sec:allmachext}. Section~\ref{sec: Numres} presents numerical results in one and two space dimensions, covering the range from low-Mach-number flows to flows with shocks and discontinuities. Finally, Section~\ref{sec: conclusions} draws the conclusions and gives an outlook on future work.

\section{HTC formulation of the governing PDE system}
	\label{sec: PDEHTC}
	In what follows, we briefly recall the HTC framework for the class of hyperbolic PDE introduced by Godunov and Romenski in \cite{God1961,Godunov:1995a,Godunov:2003a,Rom1998}.
	In particular, this implies writing the governing system of equations together with its associated total energy conservation law and the resulting compatibility conditions between fluxes,  dissipative terms and the related entropy production, which will then serve as the continuous counterpart of the discrete scheme constructed in the following sections.
	The object of study in this paper is the system of Euler equations of gasdynamics, which admits an SHTC formulation that is compatible with an extra conservation law for the total energy.
	In general, an HTC system in Einstein summation notation takes the form
	\begin{equation}
	\label{eq:pdesystemnodissip}	\partial_t \mathbf{q} + \partial_k \mathbf{f}_k(\mathbf{q}) + \mathbf{B}_k(\mathbf{q})\partial_k \mathbf{q}=0,
	\end{equation}
	and it satisfies the extra conservation law for the total energy density
	\begin{equation}
		\label{eq:energyeqnodissip}
		\partial_t \mathcal{E} + \partial_k F_k(\mathbf{q})= 0.
	\end{equation}
	Therein, $\mathbf{q}$ is the vector of state variables, $\mathbf{f}_k$ are the conservative fluxes and $\mathbf{B}_k$ contains the non-conservative terms,
	 $\mathcal{E}$ is the total energy and $F_k$ the energy flux.
	 Then, a vanishing-viscosity term $\partial_k (\epsilon\partial_k \mathbf{q})$ is added as a regularization of the PDE system, together with the associated entropy production term, denoted by $\mathbf{P}$:
		\begin{equation}
		\label{eq:pdesystem}	\partial_t \mathbf{q} + \partial_k \mathbf{f}_k(\mathbf{q}) + \mathbf{B}_k(\mathbf{q})\partial_k \mathbf{q}
		-\partial_k (\epsilon\partial_k \mathbf{q})
		=\mathbf{P},
	\end{equation}
	and consequently the total energy conservation law is modified to 
	\begin{equation}
		\label{eq:energyeq}
		\partial_t \mathcal{E} + \partial_k F_k(\mathbf{q})
		-\partial_k (\epsilon\partial_k \mathcal{E})= 0.
	\end{equation}
 	To guarantee thermodynamic compatibility of the full system \eqref{eq:pdesystem} with the total energy conservation \eqref{eq:energyeq}, it must hold
	\begin{equation}
		\label{eq:compatibility}
		\mathbf{p}\cdot\partial_t \mathbf{q}=\partial_t \mathcal{E}
		\vphantom{\mathbf{p}\cdot\partial_t \mathbf{q}=\partial_t \mathcal{E}},
	\end{equation}
	where $\mathbf{p}=\partial_{\mathbf{q}} \mathcal{E}$ is the vector of dual variables or the so-called main field \cite{Ruggeri81}.
	Moreover, this translates into the compatibility between the fluxes and between the dissipation terms and the entropy production as
	\begin{subequations}\label{eq:continuous_comp}
		\begin{align}
			\mathbf{p}\cdot (\partial_k \mathbf{f}_k(\mathbf{q}) +  \mathbf{B}_k(\mathbf{q}) \partial_k \mathbf{q}) 
			&= \partial_k F_k(\mathbf{q}), \label{eq:continuous_comp_fluxes}\\
			\mathbf{p}\cdot (\mathbf{P} +  
			\partial_k(\epsilon \partial_k \mathbf{q})) 
			&= \partial_k(\epsilon \partial_k \mathcal{E}). \label{eq:continuous_comp_dissip}
		\end{align}
	\end{subequations}
	We recall that throughout this paper, repeated Latin indices $k,l$ denote summation over the spatial dimensions 
	(Einstein summation convention), whereas the dot product $\mathbf{p}\cdot(\cdot)$ is the Euclidean 
	product taken over the components of the state vector.
As a particular example of an HTC system, in this work, we consider the compressible Euler equations. Within the HTC framework, the Euler equations can be formulated in terms of the state variables $\mathbf{q}=(\rho,\rho u_l,\rho S)^\top$, where $\rho >0$ denotes the density, $u_l$ the velocity vector and $\rho S$ the entropy density, which is treated as an independent state variable.
The thermodynamics is specified by an equation of state for the specific internal energy $e=e(\rho,S)$, from which the pressure is obtained as
\begin{equation}
	\label{eq: pressure_eos}
	p = \rho^2 \partial_\rho e(\rho,S),
\end{equation}
where the equation of state is assumed to be such that the total energy $\mathcal{E}=\mathcal{E}(\mathbf{q})$ is a strictly convex function of the state variables $\mathbf{q}$.
The governing equations with vanishing-viscosity regularization are
\begin{subequations}
		\label{eq:eulerHTC}
	\begin{align}
		\partial_t \rho
		+\partial_k(\rho u_k)
		-\partial_k\left(\epsilon\,\partial_k \rho\right)
		&=0,\\
		\partial_t(\rho u_l)
		+\partial_k\left(\rho u_l u_k+p\,\delta_{lk}\right)
		-\partial_k\left(\epsilon\,\partial_k (\rho u_l)\right)
		&=0,\\
		\partial_t(\rho S)
		+\partial_k(\rho S\,u_k)
		-\partial_k\left(\epsilon\,\partial_k (\rho S)\right)
		&=\Pi \ge 0,
	\end{align}
		\end{subequations}
    where $\mathbf{f}_k =\left(\rho u_k,\,\rho u_l u_k + p\,\delta_{lk},\,\rho S u_k\right)^\top$ are the conservative fluxes in the \(k\)-th spatial direction and $\mathbf{P} =\left(0,\,0,\,\Pi\right)^\top$ denotes the entropy production term for $\Pi \ge 0$.
	The diffusive fluxes associated with the vanishing-viscosity regularization are of the simple form $\partial_k (\epsilon\partial_k \mathbf{q})$.
	This formulation satisfies the compatibility condition \eqref{eq:compatibility} with the additional total energy conservation law:
	\begin{equation}
		\partial_t \mathcal{E}
		+\partial_k\left[u_k(\mathcal{E}+p)\right]
		-\partial_k\left(\epsilon\,\partial_k\mathcal{E}\right)
		=0.
	\end{equation}
	The total energy density is given by
	\begin{equation}
		\mathcal{E}
		=
		\rho e(\rho,S)
		+\frac12 \rho u_l u_l,
	\end{equation}
	with corresponding energy flux
	\begin{equation}
		F_k
		=
		u_k(\mathcal{E}+p).
	\end{equation}
Recently, this compatibility structure has also been exploited to construct HTC schemes that preserve thermodynamic compatibility exactly at the discrete level for different models from continuum physics \cite{HTCAbgrall,HTCAbgrall2,HTCMHD,HTCGPR,thomann_2023_thermodynamically}.
The goal of HTC schemes is to ensure that the continuous compatibility conditions \eqref{eq:continuous_comp_fluxes} and \eqref{eq:continuous_comp_dissip} hold at the discrete level as well, so that discrete total energy conservation follows directly from the construction of the method. 
In the following section we briefly recall how semi-discrete HTC schemes are constructed for general HTC systems of the form \eqref{eq:pdesystem}.

	\section{Semi-discrete HTC finite volume scheme}
	\label{sec: semidiscHTC}
	To make this paper self-contained, in this section we briefly recall the general framework outlined in \cite{SWETurbulence,HTCGPR,HTCAbgrall,HTCMHD,HTCAbgrall2} for the construction of thermodynamically compatible schemes for overdetermined hyperbolic PDE systems endowed with an extra conservation law. The framework was mainly designed for explicit time discretizations. In the present work, however, we consider semi-implicit schemes,
	which are generally more efficient than fully explicit methods in the low-Mach-number regime, as they allow larger time steps by treating the fast acoustic dynamics implicitly. In this setting, additional difficulties arise in a consistent evaluation of the terms responsible for preserving thermodynamic compatibility.
	To address this difficulty, we develop a new approach to account for semi-implicit splitting within the HTC framework. The proposed approach is of the predictor-corrector type. 

Since the fundamental concepts are valid for all HTC schemes, semi-discrete or fully discrete, we first recall the semi-discrete framework here, before introducing the new semi-implicit extension for the Euler equations in the next section.
For ease of notation, we first consider the one-dimensional case on a uniform
Cartesian grid. 
Let $\Omega\subset\mathbb R$ be partitioned into $N$ non-overlapping cells
$\Omega_i=\left[x_{i-\frac12},x_{i+\frac12}\right],$ of constant size $\Delta x =|\Omega_i|= x_{i+\frac12}-x_{i-\frac12}$. The cell center is located at $x_i=\frac12\left(x_{i-\frac12}+x_{i+\frac12}\right),$ and the semi-discrete unknowns $\mathbf{q}_i \in \mathbb{R}^M$ represent cell averages over $\Omega_i$. Therein, the index $i$ plays a different role than $k$ above: since the discretization is one-dimensional, 
there is no spatial index to sum over, and $i$ simply labels the cell $\Omega_i$ over which the cell average is taken. The semi-discrete finite volume scheme reads
		\begin{equation}
		\partial_t \mathbf{q}_i =  - \frac{\mathbf{f}_{i+\frac{1}{2}}-\mathbf{f}_{i-\frac{1}{2}}}{\Delta x}- \frac{\boldsymbol{\mathcal{D}}^-_{i+\frac{1}{2}}+\boldsymbol{\mathcal{D}}^+_{i-\frac{1}{2}}}{\Delta x}+\frac{\boldsymbol{\mathcal{G}}_{i+\frac{1}{2}}-\boldsymbol{\mathcal{G}}_{i-\frac{1}{2}}}{\Delta x}+ \mathbf{P}_i,
		\label{eqn.sdhtc}
	\end{equation}
	where at the interface $x_{i+\frac{1}{2}}$ the numerical flux is denoted by
	$\mathbf{f}_{i+\frac{1}{2}}\in \mathbb{R}^M$, $\boldsymbol{\mathcal{D}}_{i+\frac{1}{2}}^\pm \in \mathbb{R}^M$
	is the numerical discretization of the non-conservative terms, and
	$\boldsymbol{\mathcal{G}}_{i+\frac{1}{2}} \in \mathbb{R}^M$ represents the numerical discretization of
	the dissipation terms. Further, $\mathbf{P}_i\in \mathbb{R}^M$ are the related entropy production
	terms, which are rates per unit volume and therefore do not carry the
	factor $1/\Delta x$ that multiplies the interface contributions.
	
	Following the classical Tadmor decomposition of a numerical flux into an inviscid 
	contribution and a dissipative correction \cite{Tadmor1}, we assume that
	$\mathbf{f}_{i+\frac12}$ is a consistent approximation of the physical flux 
	$\mathbf{f}(\mathbf{q})$ that introduces no numerical dissipation and therefore does not contribute to the entropy production of the scheme. All the dissipation 
	required for stability of the numerical scheme is instead contained in the separate term
	$\boldsymbol{\mathcal{G}}_{i+\frac12}$. Consequently, the effective numerical flux employed
	by the scheme,
	$\mathbf{f}_{i+\frac12}-\boldsymbol{\mathcal{G}}_{i+\frac12}$,
	contains an inviscid flux plus the dissipation necessary for monotone and entropy-stable
	methods, while keeping the conservative and dissipative contributions 
	separated at the discrete level.
	Since $\mathbf{f}_{i+\frac12}$ carries no dissipation, the
	thermodynamic compatibility correction contained in the inviscid flux and discussed below in Section~\ref{sec: Thermo_flux} can be
	applied to $\mathbf{f}_{i+\frac12}$ alone, without interfering with the dissipative mechanism responsible for entropy stability, which is controlled independently by
	$\boldsymbol{\mathcal{G}}_{i+\frac12}$ and the associated entropy production terms $\mathbf{P}_i$.
	
	The discrete energy conservation law that has to hold as a consequence has the form
	\begin{equation}
		\partial_t \mathcal{E}_i =  - \frac{\mathcal{F}_{i+\frac{1}{2}}-\mathcal{F}_{i-\frac{1}{2}}}{\Delta x}+\frac{\mathcal{G}^\mathcal{E}_{i+\frac{1}{2}}-\mathcal{G}^\mathcal{E}_{i-\frac{1}{2}}}{\Delta x},
	\end{equation}
	with $\mathcal{F}_{i+\frac{1}{2}}$ the inviscid numerical energy flux and $\mathcal{G}^\mathcal{E}_{i+\frac{1}{2}}$ the dissipative energy flux at the interface $x_{i+\frac{1}{2}}$.
	On a cell $\Omega_i$, the compatibility constraint on the time derivatives \eqref{eq:compatibility} naturally reads 
	\begin{equation}\label{eq:Time_compatibility}
	\mathbf{p}_i\cdot \partial_t \mathbf{q}_i=\partial_t \mathcal{E}_i.
	\end{equation}
	 Thus, to achieve compatibility at the discrete level one needs to verify that the numerical discretization satisfies
	\begin{subequations}
		\begin{align}
			\mathbf{p}_i\cdot \left(\frac{\mathbf{f}_{i+\frac{1}{2}}-\mathbf{f}_{i-\frac{1}{2}}}{\Delta x}+ \frac{\boldsymbol{\mathcal{D}}^-_{i+\frac{1}{2}}+\boldsymbol{\mathcal{D}}^+_{i-\frac{1}{2}}}{\Delta x}\right) &=\frac{\mathcal{F}_{i+\frac{1}{2}}-\mathcal{F}_{i-\frac{1}{2}}}{\Delta x}, \label{eq:flux_thermo_compatibility}\\
			\mathbf{p}_i\cdot \left(\frac{\boldsymbol{\mathcal{G}}_{i+\frac{1}{2}}-\boldsymbol{\mathcal{G}}_{i-\frac{1}{2}}}{\Delta x}+ \mathbf{P}_i \right) &=\frac{\mathcal{G}^\mathcal{E}_{i+\frac{1}{2}}-\mathcal{G}^\mathcal{E}_{i-\frac{1}{2}}}{\Delta x},\label{eq:dissipation_compatibility}
		\end{align}
	\end{subequations}
	which are the discrete counterparts of the continuous conditions \eqref{eq:continuous_comp}. In the discrete compatibility conditions \eqref{eq:flux_thermo_compatibility} and \eqref{eq:dissipation_compatibility}, 
		 the dot product $\mathbf{p}_i\cdot(\cdot)$ is, as before, the Euclidean 
		product taken over the components of the state vector on the cell $\Omega_i$.
		It remains therefore to enforce the compatibility conditions \eqref{eq:flux_thermo_compatibility} and \eqref{eq:dissipation_compatibility} at the discrete level by carefully designing numerical fluxes and entropy production terms in the following two subsections. 
	\subsection{Thermodynamic compatibility of the numerical flux}
	\label{sec: Thermo_flux}
We first focus on condition \eqref{eq:flux_thermo_compatibility}, which involves the numerical 
fluxes and the numerical discretization of the non-conservative terms.
The central flux $\widetilde{\mathbf{f}}_{i+\frac12}$ 
\begin{equation}
	\label{eq: central_flux}
	\widetilde{\mathbf{f}}_{i+\frac12}=\frac{1}{2}(\mathbf{f}(\mathbf{q}_i)+\mathbf{f}(\mathbf{q}_{i+1}))
\end{equation}
does not automatically satisfy \eqref{eq:flux_thermo_compatibility} in general, and has to be adjusted.
In \cite{HTCAbgrall,HTCAbgrall2,thomann_2023_thermodynamically}, this is achieved by means of a
local flux correction, computed independently at each interface, so that according to Abgrall \cite{Abgrall2018} the numerical flux $\mathbf{f}_{i+\frac12}$ takes the form
	\begin{equation}
	\mathbf{f}_{i+\frac12}
	: =
	\widetilde{\mathbf{f}}_{i+\frac12}
	-
	\alpha_{i+\frac12}
	\left(
	\mathbf{p}_{i+1}-\mathbf{p}_i
	\right),
	\label{eq:local_abgrall_flux}
\end{equation}
where $\alpha_{i+\frac12}$ is a local correction parameter which is motivated and derived in the following.

Since the total energy density $\mathcal{E}$ is not among the primary evolution variables of an HTC scheme, but is instead obtained as the consequence of a compatible discretization, similar to conservative time integrators for Hamiltonian ODE systems \cite{Brugnano1,Brugnano2}, it cannot be advanced directly through a conservative standard finite volume update.  
Nevertheless, condition
\eqref{eq:flux_thermo_compatibility} prescribes exactly this structure for the total energy: its
right-hand side is written as the difference of a numerical energy flux
$\mathcal{F}_{i+\frac12}$ across the two interfaces of cell $\Omega_i$. For this to be
meaningful, $\mathcal{F}_{i+\frac12}$ must be a well-defined single-valued quantity at each interface, which is
shared consistently by the two adjacent cells $\Omega_i$ and $\Omega_{i+1}$. 
Indeed, if a consistent numerical energy flux 
$\mathcal{F}_{i+\frac12}$ exists, local energy conservation is obvious and therefore also global energy conservation, see e.g.  \cite{Abgrall2018,HTCAbgrall,HTCAbgrall2}.

Following the classical approach used in the analysis of Riemann-solver-based methods \cite{LeVeque2002a}, we decompose the contribution of each
interface into a left-going and a
right-going fluctuation, acting on the cells $\Omega_i$ and $\Omega_{i+1}$,
respectively. 
Since the left- and right-going fluctuations are computed independently, they are, in general, not guaranteed to be consistent with a single, well-defined numerical energy flux $\mathcal{F}_{i+\frac12}$ at the interface.
Requiring their consistency is
therefore equivalent to enforcing the existence of a unique
numerical energy flux at the interface, which in turn guarantees discrete energy
conservation.
Imposing this consistency condition on the corrected numerical flux
\eqref{eq:local_abgrall_flux} leads to an explicit relation for the correction
parameter $\alpha_{i+\frac12}$.
The energy fluctuations are connected to the numerical energy fluxes by
\begin{equation}
	\mathcal{R}^{\mathcal{E},-}_{i+\frac12}
	=
	\mathcal{F}_{i+\frac12}-F(\mathbf{q}_i),
	\qquad
	\mathcal{R}^{\mathcal{E},+}_{i-\frac12}
	=
	F(\mathbf{q}_i)-\mathcal{F}_{i-\frac12},
	\label{eq:fluct_en}
\end{equation}
and thus the right-hand side of \eqref{eq:flux_thermo_compatibility} is
\begin{equation}
		\frac{
		\mathcal{F}_{i+\frac12}-\mathcal{F}_{i-\frac12}
	}{\Delta x}=\frac{
		\mathcal{R}^{\mathcal{E},-}_{i+\frac12}
		+
		\mathcal{R}^{\mathcal{E},+}_{i-\frac12}
	}{\Delta x}
	.
	\label{eq:fluctuation_ene_flux}
\end{equation}
We now recast the left-hand side of the compatibility requirement \eqref{eq:flux_thermo_compatibility} in terms of local total-energy fluctuations across each cell interface by adding and subtracting 
the physical flux $\mathbf{f}(\mathbf{q}_i)$. Thus, in order for \eqref{eq:flux_thermo_compatibility} to hold, the fluctuations need to satisfy
	\begin{align}
		\label{eq: leftfluct}
		\mathcal{R}^{\mathcal{E},-}_{i+\frac12}
		&=
		\mathbf{p}_i \cdot
		\left(\mathbf{f}_{i+\frac12}-\mathbf{f}(\mathbf{q}_i)\right)
		+
		\mathbf{p}_i \cdot \boldsymbol{\mathcal{D}}^{-}_{i+\frac12},
		\\
		\label{eq: rightfluct}
		\mathcal{R}^{\mathcal{E},+}_{i-\frac12}
		&=
		\mathbf{p}_i \cdot
		\left(\mathbf{f}(\mathbf{q}_i)-\mathbf{f}_{i-\frac12}\right)
		+
		\mathbf{p}_i \cdot \boldsymbol{\mathcal{D}}^{+}_{i-\frac12},
	\end{align}
	where the non-conservative terms $\boldsymbol{\mathcal{D}}^{\pm}_{i\mp\frac12}$ are discretized according to a path-conservative scheme \cite{Pares2006,Castro2006} based on a straight-line segment path and computed using the midpoint rule, i.e.
	\begin{equation}
		\boldsymbol{\mathcal{D}}^{\pm}_{i+\frac12}
		=
		\frac{1}{2}\,
		\mathbf{B}\left(\frac{\mathbf{q}_i+\mathbf{q}_{i+1}}{2}\right)
		\left(\mathbf{q}_{i+1}-\mathbf{q}_i\right).
		\label{eq:D_pm_midpoint}
	\end{equation}
	Therein $\mathbf{B}(\mathbf{q})$ denotes the matrix associated with the non-conservative products of the system.
	If the numerical energy flux $\mathcal F_{i+\frac12}$ is single-valued at the interface, summing the left and the right energy fluctuations at the same interface $i+\frac12$ yields the consistency requirement
	\begin{equation}
		\mathcal{R}^{\mathcal{E},+}_{i+\frac12}
		+
		\mathcal{R}^{\mathcal{E},-}_{i+\frac12}
		=
		F(\mathbf{q}_{i+1})-F(\mathbf{q}_i).
		\label{eq:interface_energy_compatibility}
	\end{equation}
	Thus, \eqref{eq:interface_energy_compatibility} must be enforced in the construction 
	of the scheme in order to guarantee that a single-valued $\mathcal F_{i+\frac12}$ exists at 
	each interface, and with it, conservation of the total energy. 
	The correction parameter $\alpha_{i+\frac12}$ is therefore computed by imposing the thermodynamic compatibility relations \eqref{eq: leftfluct}, \eqref{eq: rightfluct} and \eqref{eq:interface_energy_compatibility}, using the corrected numerical flux \eqref{eq:local_abgrall_flux} as $\mathbf{f}_{i+\frac12}$:
	\begin{align}
		\alpha_{i+\frac{1}{2}}
		=&
		\frac{
			F(\mathbf{q}_{i+1})-F(\mathbf{q}_i)
			+\widetilde{\mathbf{f}}_{i+\frac{1}{2}}\cdot\left(\mathbf{p}_{i+1}-\mathbf{p}_i\right)
			-\left(\mathbf{p}_{i+1}\cdot \mathbf{f}(\mathbf{q}_{i+1})-\mathbf{p}_{i}\cdot \mathbf{f}(\mathbf{q}_{i})\right)
		}{
			\left(\mathbf{p}_{i+1}-\mathbf{p}_{i}\right)\cdot\left(\mathbf{p}_{i+1}-\mathbf{p}_{i}\right)
		}\nonumber \\
		&-
		\frac{
			\mathbf{p}_{i+1}\cdot \boldsymbol{\mathcal{D}}_{\,i+\frac{1}{2}}^{+}
			+\mathbf{p}_{i}\cdot \boldsymbol{\mathcal{D}}_{\,i+\frac{1}{2}}^{-}
		}{
			\left(\mathbf{p}_{i+1}-\mathbf{p}_{i}\right)\cdot\left(\mathbf{p}_{i+1}-\mathbf{p}_{i}\right)
		}.
		\label{eq: alpha}
	\end{align}
	If $\mathbf{p}_{i+1}-\mathbf{p}_{i}=0$, the correction term vanishes and we set $\alpha_{i+\frac{1}{2}}=0$.
	The Abgrall flux \eqref{eq:local_abgrall_flux} with the correction parameter \eqref{eq: alpha} ensures the discrete thermodynamic compatibility \eqref{eq:flux_thermo_compatibility}.
\subsection{Discrete compatibility of dissipative terms}
\label{sec: Dissip_comp_1d}
Let us now consider the compatibility of the discrete dissipative operators as expressed in \eqref{eq:dissipation_compatibility}. 
As discussed above, at the discrete level, $\boldsymbol{\mathcal{G}}_{i+\frac12}$ represents the  dissipative flux introduced alongside the inviscid flux $\mathbf{f}_{i+\frac12}$ \eqref{eq:local_abgrall_flux}, in order to ensure the stability of the scheme.
It is defined as
\begin{align}\label{eq:dissip_terms}
	\boldsymbol{\mathcal{G}}_{i+\frac{1}{2}}
	=
	\epsilon_{i+\frac12}
	\frac{\Delta \mathbf{q}_{i+\frac12}}{\Delta x}.
\end{align}
Therein we have used the notation $\Delta (\cdot)_{i+\frac12} = (\cdot)_{i+1} - (\cdot)_i$ and $\epsilon_{i+\frac12}$ represents a non-negative scalar numerical viscosity coefficient at the interface $i+\frac12$. It can, in principle, be chosen as a constant; in our case, we instead use a Rusanov-type form given by
\begin{equation}
	\label{eq:dissipation_eps}
	\epsilon_{i+\frac12}
	=
	\frac12 \Delta x\, s^{\max}_{i+\frac12},
	\qquad
	s^{\max}_{i+\frac12}
	=
	\max_j
	\Big(
	|\lambda_j(\mathbf{q}_i)|,
	|\lambda_j(\mathbf{q}_{i+1})|
	\Big),
\end{equation}
where $\lambda_j(\mathbf{q})$ denotes the $j$-th eigenvalue of the flux Jacobian	$\frac{\partial \mathbf{f}}{\partial \mathbf{q}},$
and $s^{\max}_{i+\frac12}$ is the maximum absolute characteristic wave speed at the interface. 
Substituting the definition of the dissipation terms \eqref{eq:dissip_terms} into \eqref{eq:dissipation_compatibility}
and adding and subtracting the terms $\frac{1}{2}\mathbf{p}_{i+1}\cdot\boldsymbol{\mathcal{G}}_{i+\frac{1}{2}}$ and $\frac{1}{2}\mathbf{p}_{i-1}\cdot\boldsymbol{\mathcal{G}}_{i-\frac{1}{2}}$, we obtain
\begin{eqnarray}
	&& \phantom{-} \frac{1}{2}\frac{\mathbf{p}_i+\mathbf{p}_{i+1}}{\Delta x}
	\cdot \epsilon_{i+\frac12}\frac{\Delta \mathbf{q}_{i+\frac12}}{\Delta x}
	-\frac{1}{2}\frac{\mathbf{p}_i+\mathbf{p}_{i-1}}{\Delta x}
	\cdot \epsilon_{i-\frac12}\frac{\Delta \mathbf{q}_{i-\frac12}}{\Delta x}
	\nonumber \\[0.3em]
	&& -\frac{1}{2}\frac{\mathbf{p}_{i+1}-\mathbf{p}_i}{\Delta x}
	\cdot \epsilon_{i+\frac12}\frac{\Delta \mathbf{q}_{i+\frac12}}{\Delta x}
	-\frac{1}{2}\frac{\mathbf{p}_i-\mathbf{p}_{i-1}}{\Delta x}
	\cdot \epsilon_{i-\frac12}\frac{\Delta \mathbf{q}_{i-\frac12}}{\Delta x}
	+\mathbf{p}_i\cdot \mathbf{P}_i
	=\frac{\mathcal{G}^\mathcal{E}_{i+\frac{1}{2}}-\mathcal{G}^\mathcal{E}_{i-\frac{1}{2}}}{\Delta x}.
	\label{eq:explicit_dissipation}
\end{eqnarray}
From \eqref{eq:explicit_dissipation} we can immediately identify the conservative dissipative flux of the total energy as
\begin{equation}
	\mathcal{G}^\mathcal{E}_{i+\frac{1}{2}} = 
	\frac{1}{2} \left( \mathbf{p}_i+\mathbf{p}_{i+1} \right)
	\cdot \epsilon_{i+\frac12}\frac{\Delta \mathbf{q}_{i+\frac12}}{\Delta x}
	\approx \epsilon_{i+\frac12} \frac{\mathcal{E}_{i+1}-\mathcal{E}_i}{\Delta x},
\end{equation}
since the term $\frac{1}{2} \left( \mathbf{p}_i+\mathbf{p}_{i+1} \right) \cdot \Delta \mathbf{q}_{i+\frac12}$ is a second-order approximation of the path integral 
\begin{equation}
	\label{eq: diss_comp_ene}
	\Delta \mathcal{E}_{i+\frac{1}{2}}=\mathcal{E}_{i+1}-\mathcal{E}_i
	=
	\int_{\mathbf{q}_i}^{\mathbf{q}_{i+1}}
	\partial_{\mathbf{q}} \mathcal{E}\cdot d\mathbf{q}
	=
	\int_{\mathbf{q}_i}^{\mathbf{q}_{i+1}}
	\mathbf{p}\cdot d\mathbf{q}
	\approx
	\frac12
	\left(
	\mathbf{p}_{i+1}+\mathbf{p}_i
	\right)
	\cdot
	\Delta \mathbf{q}_{i+\frac12}
\end{equation}
via the trapezoidal rule. Hence, it follows that the production term $\mathbf{P}_i$ has to compensate the remaining terms in \eqref{eq:explicit_dissipation} and it is therefore defined as
\begin{equation}
	\mathbf{p}_i\cdot \mathbf{P}_i
	=
	\frac12
	\frac{\Delta \mathbf{p}_{i+\frac12}}{\Delta x}
	\cdot
	\epsilon_{i+\frac12}
	\frac{\Delta \mathbf{q}_{i+\frac12}}{\Delta x}
	+
	\frac12
	\frac{\Delta \mathbf{p}_{i-\frac12}}{\Delta x}
	\cdot
	\epsilon_{i-\frac12}
	\frac{\Delta \mathbf{q}_{i-\frac12}}{\Delta x}.
	\label{eq:production_term}
\end{equation}
To show that \eqref{eq:production_term} leads to a consistent entropy production, we make use of the discrete jump relation
\begin{equation}
	\mathbf p_{i+1}-\mathbf p_i
	=\tilde{\mathcal{E}}_{qq,i+\frac12}\left(\mathbf{q}_{i+1}-\mathbf{q}_i\right),
\end{equation}
where $\tilde{\mathcal{E}}_{qq,i+\frac12}$ is a Roe-type matrix defined as
\begin{equation}
	\tilde{\mathcal{E}}_{qq,i+\frac12}
	=
	\int_0^1
	\mathcal{E}_{qq}\left(\boldsymbol\phi(\xi)\right)\,d \xi.
	\label{eq:comp_dissip_1d}
\end{equation}
Therein $\mathcal{E}_{qq}$ is the Hessian of the total energy and $\boldsymbol\phi(\xi) = \mathbf{q}_i + \xi \left(\mathbf{q}_{i+1} - \mathbf{q}_i\right)$ is the linear segment path between the states $ \mathbf{q}_i$ and $ \mathbf{q}_{i+1}$.
Substituting this relation into \eqref{eq:production_term}, we can write
\begin{equation}
	\label{eq: production_roe_1d}
	\mathbf{p}_i\cdot \mathbf{P}_i
	=
	\frac12
	\epsilon_{i+\frac12}
	\,
	\frac{\Delta \mathbf{q}_{i+\frac12}}{\Delta x} \cdot
	\,\tilde{\mathcal{E}}_{qq,i+\frac12}\,
	\frac{\Delta \mathbf{q}_{i+\frac12}}{\Delta x}+ 
	\frac12
	\epsilon_{i-\frac12}
	\,
	\frac{\Delta \mathbf{q}_{i-\frac12}}{\Delta x} \cdot
	\,\tilde{\mathcal{E}}_{qq,i-\frac12}\,
	\frac{\Delta \mathbf{q}_{i-\frac12}}{\Delta x} \geq 0.
\end{equation}
If the energy potential $\mathcal{E}$ is strictly convex, the Hessian $\mathcal{E}_{qq}$ is positive definite by construction, and consequently so is the Roe-type matrix $\tilde{\mathcal{E}}_{qq,i+\frac12}$. 
Since the production term is allowed only in the entropy equation because mass, momentum and total energy must be conserved, we define 
$\mathbf{P}_i = \left( 0, 0, \Pi_i \right)^\top$ with the non-negative entropy production 
\begin{equation}
	\Pi_i = \frac{1}{T_i} \left( 
	\frac12
	\epsilon_{i+\frac12}
	\,
	\frac{\Delta \mathbf{q}_{i+\frac12}}{\Delta x} \cdot
	\,\tilde{\mathcal{E}}_{qq,i+\frac12}\,
	\frac{\Delta \mathbf{q}_{i+\frac12}}{\Delta x}+ 
	\frac12
	\epsilon_{i-\frac12}
	\,
	\frac{\Delta \mathbf{q}_{i-\frac12}}{\Delta x} \cdot
	\,\tilde{\mathcal{E}}_{qq,i-\frac12}\,
	\frac{\Delta \mathbf{q}_{i-\frac12}}{\Delta x}
	\right) \geq 0. 
\end{equation}

Together with the flux compatibility condition \eqref{eq:flux_thermo_compatibility} discussed above, this provides the last ingredient 
needed to obtain a semi-discrete thermodynamically compatible finite volume scheme. It is important to stress that 
thermodynamic compatibility is established here only at the semi-discrete level and is not automatically inherited 
by the fully discrete scheme once a time integrator is applied. 
It must be enforced separately, taking into account the fully discrete nature of the numerical scheme.
This is addressed in Section~\ref{sec: Numericalscheme}, where fully 
discrete compatibility is established for our new semi-implicit scheme.	

\subsection{Thermodynamic compatibility in multiple space dimensions}
For completeness, we also provide a short formulation in multiple space dimensions.
Let $\Omega_\ell$ be a control volume and let $\mathcal N_\ell$ denote the set
of neighboring cells. The common interface between two neighboring control
volumes $\Omega_\ell$ and $\Omega_r$ is denoted by
$\partial\Omega_{\ell r}$ and the associated outward unit normal vector by
$\mathbf n_{\ell r}$.
The semi-discrete finite volume scheme reads
\begin{equation}
	\partial_t \mathbf{q}_\ell
	=
	-\sum_{r\in\mathcal N_\ell}
	\frac{|\partial\Omega_{\ell r}|}{|\Omega_\ell|}
	\left(
	\mathbf{f}_{\ell r}
	+
	\boldsymbol{\mathcal{D}}_{\ell r}
	\right)\cdot \mathbf n_{\ell r}
	+
	\sum_{r\in\mathcal N_\ell}
	\frac{|\partial\Omega_{\ell r}|}{|\Omega_\ell|}
	\boldsymbol{\mathcal{G}}_{\ell r}\cdot \mathbf n_{\ell r}
	+
	\sum_{r\in\mathcal N_\ell}
	\frac{|\partial\Omega_{\ell r}|}{|\Omega_\ell|}
	\mathbf{P}_{\ell r},
	\label{eq:fv_scheme_2d}
\end{equation}
where $\mathbf{f}_{\ell r}$ are the inviscid numerical flux of the PDE system, $\boldsymbol{\mathcal{D}}_{\ell r}$ the non-conservative contributions, $\boldsymbol{\mathcal{G}}_{\ell r}$ the dissipative numerical flux and $\mathbf{P}_{\ell r}$ the associated entropy production term at the cell interface. 
The corresponding discrete total-energy conservation law is
\begin{equation}
	\partial_t \mathcal E_\ell
	=
	-\sum_{r\in\mathcal N_\ell}
	\frac{|\partial\Omega_{\ell r}|}{|\Omega_\ell|}
	\mathcal F_{\ell r} \cdot \mathbf n_{\ell r}
	+
	\sum_{r\in\mathcal N_\ell}
	\frac{|\partial\Omega_{\ell r}|}{|\Omega_\ell|}
	\mathcal G^{\mathcal E}_{\ell r}\cdot \mathbf n_{\ell r}.
	\label{eq:energy_scheme_2d}
\end{equation}
The relation  
$\mathbf{p}_\ell\cdot \partial_t \mathbf{q}_\ell=\partial_t \mathcal E_\ell$
obviously holds by construction. Furthermore, the numerical discretization must satisfy the compatibility conditions 
\begin{subequations}
	\begin{align}
		\mathbf p_\ell \cdot
		\left[
		\sum_{r\in\mathcal N_\ell}
		\frac{|\partial\Omega_{\ell r}|}{|\Omega_\ell|}
		\left(
		\mathbf{f}_{\ell r}
		+
		\boldsymbol{\mathcal{D}}_{\ell r}
		\right)\cdot \mathbf n_{\ell r}
		\right]
		&=
		\sum_{r\in\mathcal N_\ell}
		\frac{|\partial\Omega_{\ell r}|}{|\Omega_\ell|}
		\mathcal F_{\ell r}\cdot \mathbf n_{\ell r},\label{eq:comp_flux_multid}
		\\
		\mathbf p_\ell \cdot
		\left[
		\sum_{r\in\mathcal N_\ell}
		\frac{|\partial\Omega_{\ell r}|}{|\Omega_\ell|}\left(
		\boldsymbol{\mathcal{G}}_{\ell r}\cdot \mathbf n_{\ell r}
		+
		\mathbf{P}_{\ell r}\right)
		\right]
		&=
		\sum_{r\in\mathcal N_\ell}
		\frac{|\partial\Omega_{\ell r}|}{|\Omega_\ell|}
		\mathcal G^{\mathcal E}_{\ell r}\cdot \mathbf n_{\ell r}.\label{eq:comp_dissip_multid}
	\end{align}
\end{subequations}
Compatibility of the fluxes \eqref{eq:comp_flux_multid} is ensured by defining a local Abgrall correction on the baseline central
numerical flux $\widetilde{\mathbf{f}}_{\ell r}\in\mathbb R^{M\times d}$
\begin{equation}
	\mathbf{f}_{\ell r}
	=
	\widetilde{\mathbf{f}}_{\ell r}
	-
	\alpha_{\ell r}
	\,
	(\mathbf p_r-\mathbf p_\ell)
	\otimes
	\mathbf n_{\ell r},
	\qquad 
	\textnormal{ with } \qquad 
	\widetilde{\mathbf{f}}_{\ell r} = \frac12 
	\left( 
	\mathbf{f}(\mathbf q_r) + \mathbf{f}(\mathbf q_\ell) 
	\right), 
	\label{eq:abgrall_flux_2d}
\end{equation}
where $\otimes$ denotes the tensor product and the local Abgrall correction factor is defined as 
\begin{align}
	\alpha_{\ell r}
	&=
	\frac{
		\big(F(\mathbf{q}_r)-F(\mathbf{q}_\ell)\big)\cdot\mathbf n_{\ell r}
		+
		\big(
		\widetilde{\mathbf{f}}_{\ell r}\cdot\mathbf n_{\ell r}
		\big)
		\cdot
		(\mathbf p_r-\mathbf p_\ell)
	}
	{
		(\mathbf p_r-\mathbf p_\ell)\cdot(\mathbf p_r-\mathbf p_\ell)
	}
	\nonumber\\
	&\quad
	-
	\frac{
		\Big[
		\mathbf p_r\cdot
		\big(
		\mathbf{f}(\mathbf{q}_r)\cdot\mathbf n_{\ell r}
		\big)
		-
		\mathbf p_\ell\cdot
		\big(
		\mathbf{f}(\mathbf{q}_\ell)\cdot\mathbf n_{\ell r}
		\big)
		\Big]
	}
	{
		(\mathbf p_r-\mathbf p_\ell)\cdot(\mathbf p_r-\mathbf p_\ell)	}
	-
	\frac{
		\mathbf p_\ell\cdot
		(\boldsymbol{\mathcal{D}}_{\ell r}\cdot\mathbf n_{\ell r})
		+
		\mathbf p_r\cdot
		(\boldsymbol{\mathcal{D}}_{r\ell}\cdot\mathbf n_{r\ell})
	}
	{
		(\mathbf p_r-\mathbf p_\ell)\cdot(\mathbf p_r-\mathbf p_\ell)
	}.
	\label{eq:alpha_2d}
\end{align}
As in the one-dimensional case we set $\alpha_{\ell r}=0$ when the denominator vanishes.
Let the dissipative numerical flux be
\begin{equation}
	\boldsymbol{\mathcal{G}}_{\ell r}\cdot\mathbf n_{\ell r}
	=
	\epsilon_{\ell r}
	\frac{\mathbf{q}_r-\mathbf{q}_\ell}{\delta_{\ell r}},
\end{equation}
where $\delta_{\ell r}=\left\lVert \mathbf x_r-\mathbf x_\ell \right\rVert$ denotes the distance between neighboring cell
barycenters and 
\begin{equation}
	\label{eq: eps_multid}
	\epsilon_{\ell r}
	=
	\frac12 \delta_{\ell r}
	\,s^{\max}_{\ell r},
	\qquad
	s^{\max}_{\ell r}
	=
	\max \limits_j
	\Big(
	|\lambda_j(\mathbf{q}_\ell)|,
	|\lambda_j(\mathbf{q}_r)|
	\Big).
\end{equation}
Following the one-dimensional reasoning, which is not repeated here because the operations are essentially the same, thermodynamic compatibility of the dissipation terms in \eqref{eq:comp_dissip_multid} is ensured if at each interface the production term contribution is chosen such that
\begin{equation}
	\mathbf p_\ell\cdot \mathbf{P}_{\ell r}
	=
	\frac12
	(\mathbf p_r-\mathbf p_\ell)
	\cdot
	\left(
	\boldsymbol{\mathcal{G}}_{\ell r}\cdot\mathbf n_{\ell r}
	\right)
	= 
	\frac12  \frac{\epsilon_{\ell r}}{\delta_{\ell r}} 
	\left( \mathbf p_r-\mathbf p_\ell \right)  
	\cdot 
	\left( \mathbf{q}_r-\mathbf{q}_\ell \right), 
	\label{eq:prd_term_multiD}
\end{equation}
with the dissipative numerical energy flux given by
\begin{equation}
	\mathcal G^{\mathcal E}_{\ell r}\cdot \mathbf n_{\ell r} = \frac{1}{2} \frac{\epsilon_{\ell r}}{\delta_{\ell r}} 
	\left( \mathbf{p}_\ell+\mathbf{p}_{r} \right)
	\cdot \left( \mathbf{q}_r-\mathbf{q}_\ell \right). 
\end{equation}
As in the one-dimensional case, with the aid of the Roe matrix
\begin{equation}
\widetilde{\mathcal E}_{qq,\ell r}  
=
\int_0^1
\mathcal{E}_{qq}\left(\boldsymbol\phi(\xi)\right)\,d \xi, \qquad \textnormal{ with } \qquad \boldsymbol\phi(\xi) = \mathbf{q}_\ell + \xi \left( \mathbf{q}_{r} - \mathbf{q}_\ell \right),  
\end{equation}
we can write the compatibility constraint 
\eqref{eq:prd_term_multiD} as
\begin{equation}
	\mathbf p_\ell \cdot \mathbf{P}_{\ell r}
	=
	\frac12 \frac{\epsilon_{\ell r}}{\delta_{\ell r}} 
	\left( \mathbf p_r - \mathbf p_\ell \right) 
	\cdot
	\left( \mathbf{q}_r-\mathbf{q}_\ell \right) 
	=
	\frac12
	\frac{\epsilon_{\ell r}}{\delta_{\ell r}} 
	\,
	\left(\mathbf{q}_r-\mathbf{q}_\ell\right)
	\cdot
	\widetilde{\mathcal{E}}_{qq,\ell r}
	\,
	\left( \mathbf{q}_r-\mathbf{q}_\ell \right),
	\label{eq:comp_prod_multid}
\end{equation}   
with $\mathbf{P}_{\ell r} = \left( 0, 0, \Pi_{\ell r} \right)^\top$ and 
the non-negative entropy production term at the interface 
\begin{equation}
	\Pi_{\ell r} = \frac{1}{2 T_\ell}   \frac{\epsilon_{\ell r}}{\delta_{\ell r}} 
	\,
	\left(\mathbf{q}_r-\mathbf{q}_\ell\right)
	\cdot
	\widetilde{\mathcal{E}}_{qq,\ell r}
	\,
	\left( \mathbf{q}_r-\mathbf{q}_\ell \right) \geq 0,
\end{equation}
that is needed to guarantee total energy conservation.

	\section{Semi-implicit HTC scheme}
	\label{sec: Numericalscheme}
	In this section we present a new thermodynamically compatible semi-implicit finite volume scheme for the Euler equations.
	The semi-implicit splitting follows the approach originally proposed by Casulli and Greenspan \cite{CasulliGreenspan1984} and by Park and Munz \cite{MunzPark} for compressible flows and by Casulli \cite{Casulli1990} for the shallow water equations, in which the governing system is split into a slow, nonlinear convective subsystem, discretized explicitly, and a fast, acoustic (pressure) subsystem, discretized implicitly. In this way, the acoustic terms, which are responsible for the most restrictive stability constraint of a fully explicit scheme, are removed from the time step restriction, while the nonlinear convective terms retain their explicit discretization.
	This property makes semi-implicit schemes particularly well suited for low-Mach-number flows, where the disparity between the fast acoustic speed and the slow convective waves would otherwise force a fully explicit scheme to take prohibitively small time steps. 
	
	In Section~\ref{sec:splitting} we present the splitting of the continuous system into a convective and a pressure subsystem. Subsequently, we derive the fully discrete, semi-implicit HTC scheme in one dimension in Section~\ref{sec: Numscheme1d}. The extension to multiple dimensions on a Cartesian grid is straightforward. In Section~\ref{sec:allmachext} we provide a pressure splitting which has better performance at higher Mach numbers.
\subsection{Semi-implicit splitting}
\label{sec:splitting}

We consider the Euler equations \eqref{eq:eulerHTC} written in the variables $\mathbf{q}=(\rho, \rho u_l, \rho S)^\top$.
They can be recast in the compact form \eqref{eq:pdesystem} where the flux is decomposed into explicit and implicit contributions
\begin{equation}
	\mathbf{f}_k(\mathbf{q})=\mathbf{f}_k^{E}(\mathbf{q})+\mathbf{f}_k^{I}(\mathbf{q}),
\end{equation}
with
\begin{equation}
	\mathbf{f}_k^{E}(\mathbf{q})=(0,\rho u_l u_k,\rho S u_k)^\top,
	\qquad
	\mathbf{f}_k^{I}(\mathbf{q})=(\rho u_k,p\delta_{lk},0)^\top.
\end{equation}
The explicit subsystem is hyperbolic and its characteristic speeds in direction $x_k$
\begin{equation}\label{eq:expl_charrhosspeeds}
	\lambda^E = (0,\,2u_k,\,u_k)
\end{equation} 
are independent of the Mach number.
The non-conservative contribution of the system is null, so the matrix $\mathbf{B}_k$ vanishes.
As before, the terms $\partial_k (\epsilon\partial_k \mathbf{q})$ denote the vanishing-viscosity regularization terms which can be interpreted as the numerical viscosity at the continuous level and $\mathbf{P}=(0,0,\Pi)^\top$ is the entropy production term.
We consider the ideal gas equation of state (EOS), where the specific internal energy $e$ is given by
\begin{equation}
	\label{eq: energy_eos}
	e(\rho,S)=\frac{p(\rho,S)}{(\gamma-1)\rho} =\frac{\rho^{(\gamma-1)}\exp(S/c_v)}{(\gamma-1)},
\end{equation}
where $\gamma$ is the ratio of specific heats and $c_v$ is the specific heat at constant volume.
The dual variables are $\mathbf{p}=\partial_{\mathbf{q}} \mathcal{E}=(\partial_{\rho} \mathcal{E}, u_l , T)$ where $T=e/c_v$ is the temperature and $\partial_{\rho} \mathcal{E}=-\frac{ u_k u_k}{2}+ \gamma e - T S$.
The Hessian of the total energy $\mathcal{E}_{\mathbf{q} \mathbf{q}}$ is given by
\begin{align}\label{eq: hessian_q}\mathcal{E}_{\mathbf{q} \mathbf{q}}=
	\begin{pmatrix}
		\displaystyle
		\frac{u_k u_k}{\rho}
		+\frac{e}{\rho}
		\left[
		\left(\gamma-\frac{S}{c_v}\right)^2
		-\gamma
		+2\frac{S}{c_v}
		\right] \qquad
		&
		\displaystyle
		-\frac{u_k}{\rho} \qquad
		&
		\displaystyle
		\frac{e}{c_v\rho}
		\left(
		\gamma-1-\frac{S}{c_v}
		\right)
		\\[3ex]
		\displaystyle
		-\frac{u_l}{\rho}
		&
		\displaystyle
		\frac{\delta_{lk}}{\rho}
		&
		0
		\\[3ex]
		\displaystyle
		\frac{e}{c_v\rho}
		\left(
		\gamma-1-\frac{S}{c_v}
		\right)
		&
		0
		&
		\displaystyle
		\frac{e}{c_v^2\rho}
	\end{pmatrix}.	
\end{align}
In~\cite{Lukacova2025,Godunov:2003a} it is demonstrated that, under the thermodynamic stability conditions, the Hessian of the total energy with respect to the conservative variables $(\rho,\rho u_l,\rho S)$ is symmetric and positive definite.
For the ideal gas equation of state considered here this can also be verified directly on the explicit expression \eqref{eq: hessian_q}, whose leading principal minors are all positive whenever $\rho>0$, $e>0$ and $\gamma>1$, the determinant being
\begin{equation}
	\label{eq: det_hessian_q}
	\det\left(\mathcal{E}_{\mathbf{q} \mathbf{q}}\right)
	=
	\frac{(\gamma-1) e^2}{c_v^2\rho^{d+2}}
	>0.
\end{equation}
This guarantees, at the discrete level, the correct sign of the production term as seen in Section~\ref{sec: Dissip_comp_1d}. 

Note that in the EOS \eqref{eq: energy_eos}, as well as in the dual variables, the specific entropy $S$ is required while the system \eqref{eq:eulerHTC} gives the evolution of the entropy density $\rho S$.
In the proposed flux splitting, the density is evolved implicitly, while the entropy density is treated explicitly. 
Thus the information to compute $S=(\rho S)/\rho$ is available at different time levels throughout the numerical scheme yielding an inconsistent entropy $S$, which has a direct and undesirable consequence for the numerical solution, e.g.\ spurious oscillations across a contact discontinuity. 

To ensure the contact-preserving property, it is therefore preferable to evolve the specific entropy $S$ directly, rather than the entropy density $\rho S$, avoiding the reconstruction $S=(\rho S)/\rho$ within the implicit pressure. 
This leads to an alternative formulation of the governing equations in terms of the state variables $\mathbf{Q} = (\rho, \rho u_l, S)^\top$, which obey the non-conservative system 
\begin{subequations}
	\label{eq:Euler-S}
	\begin{align}
		\partial_t \rho + \partial_k(\rho u_k)
		-\partial_k(\epsilon \partial_k \rho)&=0,\\
		\partial_t(\rho u_l)
		+\partial_k(\rho u_l u_k)
		+\partial_k (p\,\delta_{lk})
		-\partial_k(\epsilon \partial_k (\rho u_l))&=0,\\
		\partial_t S
		+u_k\,\partial_k S
		-\partial_k(\epsilon \partial_k S)
		&=\Pi/\rho.
	\end{align}
\end{subequations}
The vanishing-viscosity term in the entropy equation is chosen directly as a regularization of the specific entropy, which agrees with the regularization of the entropy density in \eqref{eq:eulerHTC} in the inviscid limit $\epsilon \to 0$.
It can be recast in the compact non-conservative form \eqref{eq:pdesystem} where the flux is decomposed into explicit and implicit contributions
\begin{equation}
	\label{eq: flux_impl_exp_split}
	\mathbf{f}_k(\mathbf{Q})=\mathbf{f}_k^{E}(\mathbf{Q})+\mathbf{f}_k^{I}(\mathbf{Q}),
\end{equation}
with
\begin{equation}
	\label{eq: flux_impl_exp}
	\mathbf{f}_k^{E}(\mathbf{Q})=(0,\rho u_l u_k,0)^\top,
	\qquad
	\mathbf{f}_k^{I}(\mathbf{Q})=(\rho u_k,p\delta_{lk},0)^\top.
\end{equation}
The entropy equation is treated explicitly and contains the only non-conservative contribution of the system, which can be represented by the matrix
\begin{equation}
	\mathbf{B}_k(\mathbf{Q})=
	\begin{pmatrix}
		0 & 0 & 0\\
		0 & 0 & 0\\
		0 & 0 & u_k
	\end{pmatrix},
\end{equation}
so that $\mathbf{B}_k(\mathbf{Q})\partial_k \mathbf{Q}=(0,	0, u_k\partial_k S)^\top$.
The explicit subsystem is hyperbolic and its characteristic speeds in direction $x_k$
\begin{equation}\label{eq:expl_charSpeeds}
	\lambda^E = (0,\,2u_k,\,u_k)
\end{equation} 
are independent of the Mach number.
As before, $\partial_k (\epsilon\partial_k \mathbf{Q})$ denote the vanishing-viscosity regularization terms and $\mathbf{P}=(0,0,\Pi/\rho)^\top$ an entropy production term. 
The dual variables now read $\mathbf{p}=\partial_{\mathbf{Q}} \mathcal{E}=(\partial_{\rho} \mathcal{E}, u_l , \rho T)$, where again $T=e/c_v$ is the temperature and $\partial_{\rho} \mathcal{E}=-\frac{ u_k u_k}{2}+ \gamma e $. In this case the Hessian of the total energy reads
\begin{align}\mathcal{E}_{\mathbf{Q} \mathbf{Q}}=
	\begin{pmatrix}
		\displaystyle
		\frac{u_k u_k}{\rho}
		+\gamma(\gamma-1)\frac{e}{\rho}\qquad
		&
		\displaystyle
		-\frac{u_k}{\rho} \qquad
		&
		\displaystyle
		\gamma\frac{e}{c_v}
		\\[3ex]
		\displaystyle
		-\frac{u_l}{\rho}
		&
		\displaystyle
		\frac{\delta_{lk}}{\rho}
		&
		0
		\\[3ex]
		\displaystyle
		\gamma\frac{e}{c_v}
		&
		0
		&
		\displaystyle
		\frac{\rho e}{c_v^2}
	\end{pmatrix}.	
\end{align}
Unfortunately, it is not positive definite by Sylvester's criterion, since
\begin{equation}
\det\left(\mathcal{E}_{\mathbf{Q}\mathbf{Q}}\right)
=
-\frac{\gamma e^2}{c_v^2\rho^d}
<0.
\end{equation}

 This property has important consequences for the construction of the discrete production term. As shown in Section~\ref{sec: Dissip_comp_1d}, working directly with these variables may in fact yield negative entropy production terms, precisely because of the lack of positive definiteness of the Hessian. 
With this in mind, the evolution will be carried out in terms of the new variables $\mathbf{Q}$, whereas the entropy production terms will be evaluated using the original variables $\mathbf{q}$ to ensure that the entropy production term is non-negative, and hence that the discrete entropy inequality is ensured. 
The same set of evolution variables $\mathbf{Q}$ has also been recently used in the FEEC scheme proposed 
in \cite{Zampa2} for the isentropic Euler equations, which is asymptotic-preserving in the low-Mach-number limit, but which is not thermodynamically compatible with the full Euler system. 
In what follows, we present the numerical scheme in detail. For clarity, its derivation is carried out in one spatial dimension; the extension to multiple dimensions on a Cartesian grid is straightforward direction by direction and is therefore omitted.

\subsection{The one-dimensional scheme}
\label{sec: Numscheme1d}
We apply the same discretization and notation as in Section~\ref{sec: semidiscHTC}.
In addition to the cells on the primal grid $\Omega_i=\left[x_{i-\frac12},x_{i+\frac12}\right]$ of constant size $\Delta x =|\Omega_i|= x_{i+\frac12}-x_{i-\frac12}$ we define the cells on a staggered dual grid 
$\Omega_{i+\frac{1}{2}}=[x_i,x_{i+1}]$ on which we evolve additionally the staggered discrete momentum $m_{i+\frac{1}{2}}$. 
It is well known, see e.g.\ \cite{DumbserCasulli2016}, that a staggered treatment of the momentum 
provides asymptotic consistency and accuracy in the low-Mach-number limit of \eqref{eq:Euler-S}, which is given by the incompressible Euler equations \cite{KlaMaj}. 
Accordingly, we define the staggered velocity as
\begin{equation}
	\label{eq: staggered_velocity}
	u_{i+\frac{1}{2}}=\frac{m_{i+\frac{1}{2}}}{\rho_{i+\frac{1}{2}}}.
\end{equation}
Quantities that are defined in the cells but are needed on the dual grid are transferred by simple arithmetic averaging of the two adjacent cell values, so that the staggered density is given by
\begin{equation}
	\label{eq: staggered_density}
	\rho_{i+\frac{1}{2}} = \frac{1}{2}\left( \rho_i + \rho_{i+1} \right),
\end{equation}
and the same averaging is used throughout for the transfer between the primal and the dual grid. The staggered momentum plays the role of a numerical flux rather than of a state variable, since all HTC quantities are defined in the cell centers.
The development of the numerical scheme is based on an explicit implicit strategy as e.g.\ done in \cite{MunzPark,DumbserCasulli2016}. 
Its main steps are as follows 
\begin{enumerate}
\item Given the cell averages $\mathbf{Q}^n_i$ at time $t^n$, we perform first the explicit update which we denote by $\mathbf{Q}^*$.
\item Using the explicit data $\mathbf{Q}^*$, we perform an implicit step obtaining the update $\mathbf{Q}^{**}$, by solving an implicit pressure equation.
\item After this update we recover a thermodynamically compatible formulation of the entropy density $\rho S$, thus obtaining the approximated values of the conserved variables $\mathbf{q}^{**}=(\rho^{**},m^{**},\rho S^{**})$.
To ensure thermodynamic compatibility of the fully discrete scheme, the entropy density $\rho S^{***}$ is then updated with the correct entropy production term. 
\item In the end, to ensure global energy conservation, we will correct our variables $\mathbf{q}^{***}$ with a global correction parameter $\alpha$. With this final step we have computed $\mathbf{q}^{n+1}$ and thus also $\mathbf{Q}^{n+1}$ by defining $S^{n+1}=(\rho S)^{n+1}/\rho^{n+1}$.
\end{enumerate}
\subsubsection{Explicit step}
Following a standard finite volume framework, the explicit update is given by
\begin{equation}
	\label{eq:expl_update_1d}
\mathbf{Q}^{*}_i=\mathbf{Q}^n_i-\frac{\Delta t}{\Delta x}\left(\mathbf{f}^{E,n}_{i+\frac12}-\mathbf{f}^{E,n}_{i-\frac12}\right)-\frac{\Delta t}{\Delta x}\left(\boldsymbol{\mathcal{D}}^{-,n}_{i+\frac{1}{2}}+\boldsymbol{\mathcal{D}}^{+,n}_{i-\frac{1}{2}}\right)+\frac{\Delta t}{\Delta x}\left(\boldsymbol{\mathcal{G}}_{i+\frac{1}{2}}^n-\boldsymbol{\mathcal{G}}_{i-\frac{1}{2}}^n\right),
\end{equation}
where the explicit numerical flux at the interface is given by
\begin{equation}
	\label{eq: fluxes_expl_n}
	\mathbf{f}^{E,n}_{i+\frac12}=\left(0,f^{\rho u^2,n}_{i+\frac12},0\right)^\top, \qquad f^{\rho u^2,n}_{i+\frac12}
	=
	\frac12
	u^{n}_{i+\frac12}
	\left(
	m^{n}_i
	+
	m^{n}_{i+1}
	\right).
\end{equation}
The non-conservative contributions are 
\begin{align}
	\label{eq: ncp_expl}
	\boldsymbol{\mathcal{D}}^{-,n}_{i+\frac{1}{2}}=\left(0,0,\mathcal{D}^{S,-,n}_{i+\frac12}\right)^\top \qquad	\mathcal{D}^{S,-,n}_{i+\frac12}
	=
	\frac12\,u^n_{i+\frac12}
	\left(S^n_{i+1}-S^n_i \right), \\
		\label{eq: ncm_expl}
		\boldsymbol{\mathcal{D}}^{+,n}_{i-\frac{1}{2}}=\left(0,0,\mathcal{D}^{S,+,n}_{i-\frac12}\right)^\top \qquad \mathcal{D}^{S,+,n}_{i-\frac12}
		=
		\frac12\,u^n_{i-\frac12}
		\left(S^n_i-S^n_{i-1}\right),
\end{align}
where the non-conservative products are evaluated with the staggered velocity $u^n_{i+\frac12}$ of \eqref{eq: staggered_velocity}, while the dissipation terms are defined as
\begin{equation}
		\label{eq: dissip_expl}
		\boldsymbol{\mathcal{G}}^{n}_{i+\frac{1}{2}} =\left(\mathcal{G}^{\rho,n}_{i+\frac12},\mathcal{G}^{m,n}_{i+\frac12},\mathcal{G}^{S,n}_{i+\frac12}\right)^\top=
	\epsilon_{i+\frac12}
	\frac{\Delta \mathbf{Q}^{n}_{i+\frac12}}{\Delta x}.
\end{equation}
Similarly to \eqref{eq:dissipation_eps}, we choose a Rusanov-type dissipation which only depends on the velocity according to \eqref{eq:expl_charSpeeds}. These eigenvalues not only determine the value of $\epsilon_{i+\frac12}$ appearing in the numerical dissipation terms, but, more importantly, they also control the stability of the explicit step via the CFL condition
\begin{equation}
	\Delta t
	\le
	\mathrm{CFL}\,
	\frac{\Delta x}{\max \limits_i (\max \limits_r \left|\lambda^E_r(\mathbf{Q}^n_i)\right|)}
	=
	\mathrm{CFL}\,
	\frac{\Delta x}{2\max \limits_i \left|u_i\right|},
	\qquad
	\mathrm{CFL}\in(0,1],
	\label{eq:cfl_condition}
\end{equation}
which is considerably less restrictive than the CFL condition of a fully explicit scheme, in which the acoustic eigenvalues $u\pm c$ would also enter the time step restriction.
This concludes the explicit step. 
\subsubsection{Implicit step}
According to the flux splitting, the implicit update reads
\begin{equation}
	\label{eq:impl_update_1d}
\mathbf{Q}^{**}_i=\mathbf{Q}^*_i-\frac{\Delta t}{\Delta x}\left(\mathbf{f}^{I,**}_{i+\frac12}-\mathbf{f}^{I,**}_{i-\frac12}\right),
\end{equation}
where $\mathbf{Q}_i^*$ denotes the state vector after the explicit step.
The numerical flux depends implicitly on $\mathbf{Q}^{**}$ and is given  by
\begin{equation}
	\label{eq:implicit_up_1d}
\mathbf{f}_{i+\frac12}^{I,**}=(m^{**}_{i+\frac12},p^{**}_{i+\frac12},0)^\top.
\end{equation}
Note that the specific entropy $S$ is not affected by the implicit step and therefore it holds $S^{**}=S^{*}$. 
In particular, mass and momentum obey
\begin{align}
	\rho_i^{**} &=\rho_i^{*}-\frac{\Delta t}{\Delta x}\left(m_{i+\frac12}^{**}-m_{i-\frac12}^{**}\right),\label{eq:rho_n+1_star}\\
	m_i^{**}&=m_i^{*} - \frac{
		\Delta t
	}{\Delta x}\left(p^{**}_{i+\frac12}
	-
	p^{**}_{i-\frac12}\right).\label{eq:m_n+1_star}
\end{align}
In addition, the momentum is evolved on the staggered grid applying the same splitting as on the primal grid as follows 
\begin{equation}
	\label{eq:mom_stagg_n+1_star}
	m^{**}_{i+\frac12}
	=m^{*}_{i+\frac12}-
	\frac{\Delta t}{\Delta x}\left(p^{**}_{i+1}-p^{**}_i\right), \quad m^{*}_{i+\frac12}
	=
	\frac12\left(m^{*} _{i+1}+m^{*}_i\right).
\end{equation}
This allows us to write one implicit scalar pressure wave equation, unlike the coupled system \eqref{eq:rho_n+1_star}, \eqref{eq:m_n+1_star}, by substituting \eqref{eq:mom_stagg_n+1_star} in \eqref{eq:rho_n+1_star}. 
It reads 
\begin{equation}\label{eq:pressure_implicit}	
\rho\left(p_i^{**},S_i^{**}\right)
-
\left(
\frac{\Delta t}{\Delta x}
\right)^2
\left(
p^{**}_{i+1}
-
2p^{**}_i
+
p^{**}_{i-1}
\right)
	=
	\rho_i^{*}
	-
	\frac{\Delta t}{\Delta x}
	\left(
	m^{*}_{i+\frac12}
	-
	m^{*}_{i-\frac12}
	\right),
\end{equation}
since it is possible to write $\rho=\rho(p,S)$ from the chosen equation of state, i.e.\ on the discrete level $\rho^{**}_i=\rho\left(p_i^{**},S_i^{**}\right)$.
Thus, \eqref{eq:pressure_implicit} is a mildly nonlinear implicit system for the pressure $p_i^{**}$ 
in the sense of Casulli and Zanolli, see \cite{casulli_2010_a,casulli_2012_iterative}, where the nonlinearity is only on the diagonal.
We can write the discrete pressure system more compactly as 
\begin{equation}
	\boldsymbol{\rho}\left(\mathbf{p}^{**},\mathbf{S}^{**}\right)
	+
	\mathbf{A} \, \mathbf{p}^{**}
	= \mathbf{b}^*,
	\label{eq:pressure_system}
\end{equation}
where $\mathbf{p}^{**}=(p_1^{**},p_2^{**},\ldots)$ is the vector of the unknown pressure, $\mathbf{S}^{**}=(S_1^{**},S_2^{**},\ldots)$ is the known vector of specific entropies and the matrix $\mathbf{A}$ corresponds to the linear part of the system, which is a symmetric, tridiagonal and positive semi-definite coefficient matrix. The right-hand side vector is denoted by $\mathbf{b}^{*}$ and is given by the right-hand side of \eqref{eq:pressure_implicit}. Since system \eqref{eq:pressure_system} is mildly nonlinear, it can be solved with the aid of the nested Newton method of Casulli and Zanolli, see \cite{casulli_2010_a,casulli_2012_iterative},
which is proven to converge for a suitable initial guess.
Note that the Jacobian of the pressure system, which is needed in the nested Newton method, is 
\begin{equation}
	\mathbf{J} = \mathbf{D} + \mathbf{A},
	\label{eqn.J.def} 
\end{equation}
with the diagonal matrix 
\begin{equation}
\mathbf{D} = \frac{\partial \boldsymbol{\rho}(\mathbf{p}^{**},\mathbf{S}^{**})}{\partial \mathbf{p}^{**} } = \textnormal{diag} \left( \frac{1}{(c^{**}_1)^2}, \frac{1}{(c^{**}_2)^2}, \ldots \right),
\label{eqn.D.def} 
\end{equation}
which is simply composed of the reciprocal values of the square of the sound speed in each cell, with $(c^{**}_i)^2 = c^2(\rho_i^{**},S_i^{**})$. Since $\mathbf{D}$ is a positive diagonal matrix and $\mathbf{A}$ is symmetric, positive semi-definite and has non-positive off-diagonal entries, the Jacobian $\mathbf{J}$ is symmetric and positive definite and is moreover an M-matrix. Note that by its definition the square of the adiabatic sound speed is precisely given by $c^2 = c^2(\rho,S) = \partial p(\rho,S) / \partial \rho$, i.e.\ the partial derivative of the pressure with respect to density at constant entropy. The low-Mach compatibility of the scheme is therefore immediate, as the pressure system \eqref{eq:pressure_implicit} reduces to the classical pressure Poisson equation of the incompressible Euler equations in the case when $\mathrm{Ma} \to 0$, i.e.\ $c \to \infty$. 

In the one-dimensional case, the associated linear system appearing in the Newton method is solved with the aid of the Thomas algorithm, while in the two-dimensional case a matrix-free conjugate gradient method with a tolerance of $10^{-15}$ is used. 
Once the pressure field $p^{**}$ has been obtained, the momentum on the primal and on the staggered grid is updated according to \eqref{eq:m_n+1_star} and \eqref{eq:mom_stagg_n+1_star}, respectively, where the pressure at the interface is defined by
\begin{equation}
	p^{**}_{i+\frac12}
	=
	\frac12
	\left(
	p^{**}_{i+1}
	+
	p^{**}_i
	\right).
\end{equation}
Last, the density is updated through the relation \eqref{eq:rho_n+1_star}.

\subsubsection{Restoring thermodynamic compatibility}
We have now obtained a predictor $\mathbf{Q}_i^{**}$ for the state variables which is in general not thermodynamically compatible, but which is asymptotic-preserving in the low-Mach-number limit since the discrete pressure wave equation tends to the pressure projection equation of an incompressible flow solver, as also confirmed numerically in Section~\ref{sec:tgv}.
As discussed above, ensuring the discrete compatibility of the dissipative terms with the total energy conservation law and with entropy production requires a positive definite total energy Hessian, a property that fails to hold in the non-conservative variables $\mathbf{Q}$, while it is satisfied in the original variables $\mathbf{q}$, see e.g.~\cite{Feireisl2021book}. For this reason, the entropy production term must be evaluated using $\mathbf{q} = (\rho, m, \rho S)$ as evolution variables. As a third step in the numerical scheme, we therefore first recover an updated predictor for $\rho S$, consistent with the viscous fluxes already applied to $\rho$ and to $S$, and only then add the correct, thermodynamically compatible entropy production. 
Since the mass and momentum equations are not affected by this change of variables, we simply set
\begin{equation}
	\rho^{***}_i=\rho^{**}_i, \qquad m^{***}_i=m^{**}_i.
\end{equation}
Let us define $\rho S^{**}$ as
\begin{equation}\label{eq:rhoS_pred}
	\rho S_i^{**}
	=
	\biggl(\rho_i^{**}
	-\frac{\Delta t}{\Delta x}
	\left(\mathcal{G}^{\rho,n}_{i+\frac{1}{2}}-\mathcal{G}^{\rho,n}_{i-\frac{1}{2}}\right)
	\biggr)
	\biggl(S_i^{**}
	-\frac{\Delta t}{\Delta x}
	\left(\mathcal{G}^{S,n}_{i+\frac{1}{2}}-\mathcal{G}^{S,n}_{i-\frac{1}{2}}\right)\biggr)
	+
	\frac{\Delta t}{\Delta x}
	\left(\mathcal{G}^{\rho S,n}_{i+\frac{1}{2}}-\mathcal{G}^{\rho S,n}_{i-\frac{1}{2}}\right)
	.
\end{equation}
In this expression, the first term subtracts the dissipation already embedded in the predictors $\rho_i^{**}$ and $S_i^{**}$ from the explicit step, while the last term adds back the correct dissipation in $\rho S$:
\begin{equation}
	\mathcal{G}^{\rho S,n}_{i+\frac12}
	=
	\epsilon_{i+\frac12}
	\frac{\Delta (\rho S)_{i+\frac12}^n}{\Delta x}.
\end{equation}
Formulation \eqref{eq:rhoS_pred} together with setting 
\begin{equation}
	\label{eq:entr_prod}
	\rho S_i^{***}
	=\rho S_i^{**}+ \Delta t \, \Pi_i,
\end{equation}
ensures discrete thermodynamic compatibility with the viscous terms in the discrete energy conservation law as detailed in Section~\ref{sec: Dissip_comp_1d}.
Note that $\Pi_i$ denotes a total entropy production term, which is obtained as the sum of a spatial and a temporal entropy production 
\begin{equation}
	\Pi_i=\Pi^{s}_i+\Pi^{t}_i,
\end{equation}
where $\Pi^{s}_i$ denotes the entropy production from the spatial numerical viscosity as in \eqref{eq:production_term} in the semi-discrete scheme. 
Moreover, we obtain an additional production term $\Pi^{t}_i$ due to the time discretization. 
Before its derivation, we first give the formulation of $\Pi_i^s$ specific to the case of the Euler equations. 
Since the production term only appears in the entropy, whose dual is the temperature $T$, we have
\begin{equation}
	\label{eq: pi_s_dual}
	\Pi^{s}_i=\frac{\epsilon_{i+\frac12}}{T_i^{**}}\frac{1}{2}	\frac{\Delta \mathbf{p}^{**}_{i+\frac12}}{\Delta x}
	\cdot	
	\frac{\Delta \mathbf{q}^{**}_{i+\frac12}}{\Delta x}
	+\frac{\epsilon_{i-\frac12}}{T_i^{**}}
	\frac12
	\frac{\Delta \mathbf{p}^{**}_{i-\frac12}}{\Delta x}
	\cdot	
	\frac{\Delta \mathbf{q}^{**}_{i-\frac12}}{\Delta x} \geq 0.
\end{equation}
Note that this formulation of the production term does not necessarily require the implementation or evaluation of the total energy Hessian and is non-negative provided the total energy is convex, and therefore ensures a non-negative entropy production. Indeed, one can also write it in terms of the Roe matrix as in \eqref{eq: production_roe_1d}
\begin{equation}
	\label{eq: pi_s_roe}
	\Pi^{s}_i
	=
	\frac{1}{T_i^{**}}\frac{1}{2\Delta x^2}\left(
	\epsilon_{i+\frac12}
	\left(\Delta \mathbf{q}^{**}_{i+\frac12}\right)^\top
	\,\tilde{\mathcal{E}}^{**}_{qq,i+\frac12}\,
	\Delta \mathbf{q}^{**}_{i+\frac12}
	+
	\epsilon_{i-\frac12}
	\left(\Delta \mathbf{q}^{**}_{i-\frac12}\right)^\top
	\,\tilde{\mathcal{E}}^{**}_{qq,i-\frac12}\,
	\Delta \mathbf{q}^{**}_{i-\frac12}
	\right) \geq 0. 
\end{equation}
Although mathematically equivalent, the formulation \eqref{eq: pi_s_dual} is preferred since it is computationally more efficient and easier to implement than \eqref{eq: pi_s_roe}, as it only involves vector dot products between quantities that are already available in the numerical scheme.
We turn now to the evaluation of the additional production term $\Pi^{t}_i$.
\subsubsection{Fully discrete thermodynamic compatibility}
A fully discrete scheme must satisfy the discrete analogue of the compatibility condition \eqref{eq:Time_compatibility}, which is given on an interval $[t^n,t^{n+1})$ by 
\begin{equation}
	\widetilde{\mathbf{p}}_i
	\cdot
	\left(
	\mathbf{q}_i^{n+1}-\mathbf{q}_i^n
	\right)
	=
	\mathcal{E}_i^{n+1}
	-
	\mathcal{E}_i^n,
	\label{eq:fully_discrete_compatibility} \quad \widetilde{\mathbf{p}}_i
	=
	\int_0^1
	\mathbf{p}(\psi_i(\xi))\,d\xi.
\end{equation}
In the derivation of the Roe-type discrete main field we employ a simple segment path $	
\psi_i(\xi)=\mathbf{q}_i^n+\xi\left(\mathbf{q}_i^{n+1}-\mathbf{q}_i^n\right)$ connecting the states $\mathbf{q}_i^n$ and $\mathbf{q}_i^{n+1}$, as in \cite{HTCAbgrall,HTCAbgrall2,HTCMHD,SWETurbulence,HTCGPR,thomann_2023_thermodynamically}.
Integration by parts along the path yields
\begin{equation}
	\widetilde{\mathbf{p}}_i
	=
	\mathbf{p}_i^{n+1}
	-
	\int_0^1
	\xi\,
	\mathcal{E}_{qq}(\psi_i(\xi))
	\left(
	\mathbf{q}_i^{n+1}-\mathbf{q}_i^n
	\right)
	\,d\xi.
\end{equation}
	Substituting this expression into \eqref{eq:fully_discrete_compatibility} and defining 
	\begin{equation}
		\mathcal{V}_i\left(\mathbf{q}_i^{n},\mathbf{q}_i^{n+1}\right)
		:=
		\int_0^1
		\xi\,
		\left(
		\mathbf{q}_i^{n+1}-\mathbf{q}_i^n
		\right)^\top
		\mathcal{E}_{qq}(\psi_i(\xi))
		\left(
		\mathbf{q}_i^{n+1}-\mathbf{q}_i^n
		\right)
		\,d\xi \geq 0,
		\label{eq:time_viscosity}
	\end{equation}
	 the compatibility condition becomes
\begin{equation}
	\mathbf{p}_i^{n+1}\cdot
	\left(
	\mathbf{q}_i^{n+1}-\mathbf{q}_i^n
	\right)
	-
	\mathcal{V}_i\left(\mathbf{q}_i^{n},\mathbf{q}_i^{n+1}\right)
	=
	\mathcal{E}_i^{n+1}
	-
	\mathcal{E}_i^n.
\end{equation}
The Hessian of the total energy with respect to the conservative variables is symmetric positive definite for convex potentials, hence the production term defined in \eqref{eq:time_viscosity} is non-negative and adds to the total entropy production. 
Since the entropy production is added after the implicit step, we define the entropy production due to the implicit time discretization as
\begin{equation}
	\label{eq: pi_t}
	\Pi_i^t = \frac{\mathcal{V}_i\left(\mathbf{q}_i^{n},\mathbf{q}_i^{**}\right)}{\Delta t \, T_i^{**}} \geq 0.
\end{equation}
From the above given derivation, strictly speaking, the temporal production has to be evaluated as $\mathcal{V}_i(\mathbf{q}_i^{n},\mathbf{q}_i^{***})$ instead of $\mathcal{V}_i(\mathbf{q}_i^{n},\mathbf{q}_i^{**})$. 
Since the production term is nonlinear, this would lead to a large coupled nonlinear implicit system. 
To compensate for this, we place steps 1 to 3 of the numerical scheme in a Picard loop as follows
\begin{align}
	\mathbf{Q}^{*,(k+1)}_i&=\mathbf{Q}^n_i-\frac{\Delta t}{\Delta x}\left(\mathbf{f}^{E,(k)}_{i+\frac12}-\mathbf{f}^{E,(k)}_{i-\frac12}\right)-\frac{\Delta t}{\Delta x}\left(\boldsymbol{\mathcal{D}}^{-,(k)}_{i+\frac{1}{2}}+\boldsymbol{\mathcal{D}}^{+,(k)}_{i-\frac{1}{2}}\right)+\frac{\Delta t}{\Delta x}\left(\boldsymbol{\mathcal{G}}_{i+\frac{1}{2}}^{(k)}-\boldsymbol{\mathcal{G}}_{i-\frac{1}{2}}^{(k)}\right) \label{eq: picardi_expl_upd_1d}\\
	\mathbf{Q}^{**,(k+1)}_i&=\mathbf{Q}^{*,(k+1)}_i-\frac{\Delta t}{\Delta x}\left(\mathbf{f}^{I,**,(k+1)}_{i+\frac12}-\mathbf{f}^{I,**,(k+1)}_{i-\frac12}\right)
\end{align}
where the terms on the right-hand side of \eqref{eq: picardi_expl_upd_1d} denoted by $^{(k)}$ are evaluated at the state $\mathbf{Q}^{***,(k)}$ given by \eqref{eq: fluxes_expl_n}--\eqref{eq: dissip_expl}.
Analogously to \eqref{eq:rhoS_pred}, we obtain $\rho S^{**,(k+1)}$ as
\begin{align}\label{eq:rhoS_pred_pic}
	\rho S_i^{**,(k+1)}
	=&
	\biggl(\rho_i^{**,(k+1)}
	-\frac{\Delta t}{\Delta x}
	\left(\mathcal{G}^{\rho,(k)}_{i+\frac{1}{2}}-\mathcal{G}^{\rho,(k)}_{i-\frac{1}{2}}\right)
	\biggr)
	\biggl(S_i^{**,(k+1)}
	-\frac{\Delta t}{\Delta x}
	\left(\mathcal{G}^{S,(k)}_{i+\frac{1}{2}}-\mathcal{G}^{S,(k)}_{i-\frac{1}{2}}\right)\biggr)
	\nonumber \\ 
	&+
	\frac{\Delta t}{\Delta x}
	\left(\mathcal{G}^{\rho S,(k)}_{i+\frac{1}{2}}-\mathcal{G}^{\rho S,(k)}_{i-\frac{1}{2}}\right)
\end{align}
and thus we have obtained $\mathbf{q}^{**,(k+1)}$ with which the entropy production is computed. Following \eqref{eq:entr_prod}, we finally obtain $\rho S^{***,(k+1)}$ and thus $ \mathbf{q}^{***,(k+1)}$. 
Therein, $k \in \mathbb{N}_0$ denotes the Picard index, and we set the start value as $\mathbf{q}_i^{***,(0)} = \mathbf{q}_i^n$. 
The state $\mathbf{Q}^{***,(k)}$ is obtained from $\mathbf{q}^{***,(k)}$ by $S^{***,(k)} = (\rho S)^{***,(k)}/\rho^{***,(k)}$, in the same way as at the end of the time step. 
The loop is carried out for a fixed number of iterations, and in all computations presented in this work two Picard iterations were performed, which we found to be sufficient to ensure the discrete compatibility; the result of the last iteration is then taken as $\mathbf{q}^{***}_i = \mathbf{q}^{***,(k+1)}$.
\subsubsection{Total energy conservation}
Once all Picard iterations are completed, the only thermodynamic compatibility condition still to be satisfied is the conservation of total energy. 
In the semi-discrete scheme \eqref{eqn.sdhtc}, this has been done by the local flux correction of Abgrall \eqref{eq:local_abgrall_flux}. 
Enforcing this correction implicitly leads in 1D to $N+1$ unknowns, one per cell interface, which form a large coupled nonlinear implicit system. 
Since resolving such a system numerically is quite challenging and computationally expensive, we deliberately choose to ensure global energy conservation instead, by enforcing the compatibility by means of a single global correction parameter $\alpha$ yielding the final numerical approximation $\mathbf{q}_i^{n+1}$ from which the non-conservative variables $\mathbf{Q}_i^{n+1}$ are computed.  

Similarly to the flux correction \eqref{eq:local_abgrall_flux} in the semi-discrete scheme, we define the flux-corrected update
\begin{equation}
	\label{eq:abgrall_q_corr}
	\mathbf{q}_i^{n+1}(\alpha)
	=
	\mathbf{q}_i^{***}
	-
	\frac{\Delta t}{\Delta x}
	\left[
	\widehat{f}_{i+\frac12}(\alpha)
	-
	\widehat{f}_{i-\frac12}(\alpha)
	\right],	\end{equation}
	where
\begin{equation}
	\label{eq:abgrall_flux}
	\widehat{f}_{i+\frac12}(\alpha)
	=
	-\alpha \mathbf{L}
	\left(
	\mathbf{p}^{***}_{i+1}-\mathbf{p}^{***}_i
	\right).
\end{equation}
Therein the correction factor $\alpha$ is shared by all interfaces. 
Unlike in \eqref{eq:local_abgrall_flux}, the correction is not necessarily applied to all variables in $\mathbf{q}$. 
This is managed by the binary diagonal matrix $L$. More precisely,
\begin{equation}
	\mathbf{L}=\operatorname{diag}(L_{11},\ldots,L_{MM}),
	\qquad
	L_{ll}\in\{0,1\}, 
	\label{eq:Loperator}
\end{equation}
with $L_{ll}=1$ if the $l$-th component of the state vector
$\mathbf{q}\in\mathbb R^M$ is corrected by the global Abgrall flux and
$L_{ll}=0$ otherwise. For instance, choosing $\mathbf{L}=\mathbf{I}$ applies the correction to
all components of $\mathbf{q}$, while setting selected diagonal entries to zero leaves
the corresponding variables unchanged.
Note that $\alpha \in \mathbb{R}$ has no fixed sign, as it modifies the numerical fluxes of the scheme in such a way that global total energy conservation is achieved. Since we do not want this correction to affect the 
discrete entropy inequality, given that we have already ensured a positive entropy production term, we 
choose to not apply it to the entropy equation, setting $L_{MM} = 0$, while all remaining components are corrected, so that $\mathbf{L}=\operatorname{diag}(1,\ldots,1,0)$ in all computations presented in this work.
The value of $\alpha$ is determined by imposing global conservation of total energy, hence 
let
\begin{equation}
	\label{eq:global_abgrall_constraint}
	g(\alpha)
	:=
	\sum_i \mathcal E\left(\mathbf{q}_i^{n+1}(\alpha)\right)\Delta x
	-
	\sum_i \mathcal E(\mathbf{q}_i^n)\Delta x
	+
	\Delta t\left(F(\mathbf{q}_R)-F(\mathbf{q}_L)\right)=0,
\end{equation}
where $F(\mathbf{q}_R)-F(\mathbf{q}_L)$ is the net outward
energy flux through the computational boundary, with $\mathbf{q}_L$ and $\mathbf{q}_R$ the fixed Dirichlet boundary values prescribed at the left and at the right boundary of the computational domain, which are assumed constant in time.  
The nonlinear system given by \eqref{eq:abgrall_q_corr} and \eqref{eq:global_abgrall_constraint} is solved with the Newton method, using $\alpha = 0$ as initial guess and the relative total energy error below $\varepsilon = 10^{-14}$ as stopping criterion.
After the computation of $\alpha$ we set $\mathbf{q}^{n+1}$ from \eqref{eq:abgrall_q_corr} and $\mathbf{Q}^{n+1}=(\rho^{n+1},m^{n+1},S^{n+1})$ with $S^{n+1}=\dfrac{\rho S^{n+1}}{\rho^{n+1}}$.
The semi-implicit scheme conserves mass and momentum locally, since all its stages are in flux form, and it also satisfies the entropy inequality for the entropy density $\rho S$ in each cell. However, it conserves total energy only globally, the local cellwise compatibility of Section~\ref{sec: semidiscHTC} being replaced by the single global constraint \eqref{eq:global_abgrall_constraint}. The numerical results of Section~\ref{sec: Numres} show that despite the resulting local energy defect the position and the strength of the computed shocks are still in good agreement with the exact solution.
\subsection{Pressure splitting for higher Mach numbers}
\label{sec:allmachext}

The implicit-explicit flux splitting \eqref{eq: flux_impl_exp_split} applied in the above scheme treats only convective terms explicitly and the entire pressure subsystem implicitly. 
This choice is well suited for the low-Mach-number regime, where the fast acoustic waves are handled implicitly and thus neither contribute to the CFL condition nor cause a scale-dependent dissipation. 
However, as the Mach number increases and the flow approaches the sonic regime, the resolution of all wave phenomena is of interest.
Moreover, material and acoustic effects become comparable in magnitude and speed.  
Consequently, treating the entire pressure contribution implicitly is no longer necessary as it does not negatively influence the CFL condition and the dissipation. 
On the contrary, additional dissipation stabilizes the numerical scheme in the presence of strong shocks and in practice improves the convergence of the Newton algorithms, which otherwise would require a reduction of the CFL number to obtain numerical results. 

Therefore, we propose a pressure splitting approach, see e.g.~\cite{THOMANN2020109723,DegondTang2011} and references therein, which is based on a blending coefficient $\beta=\beta(\mathrm{Ma})$ depending on the local Mach number $\mathrm{Ma}$.
The new flux splitting, depending on $\beta$, reads 
\begin{equation}
	\label{eq: flux_impl_exp_allmach}
	\mathbf{f}_k^{E}(\mathbf{Q})=\bigl((1-\beta)\rho u_k,\,
	\rho u_lu_k+(1-\beta)p\,\delta_{lk},\,0\bigr)^\top,
	\qquad
	\mathbf{f}_k^{I}(\mathbf{Q})=\bigl(\beta\rho u_k,\,
	\beta p\,\delta_{lk},\,0\bigr)^\top.
\end{equation}
The explicit subsystem is strictly hyperbolic for $\beta<1$, with characteristic speeds in direction $x_k$
\begin{equation}
	\label{eq:expl_charSpeeds_allmach}
	\lambda^E=
	\left(
	u_k,\,
	u_k-\sqrt{\beta u_k^2+(1-\beta)^2c^2},\,
	u_k+\sqrt{\beta u_k^2+(1-\beta)^2c^2}
	\right),
\end{equation}
where $c^2=\partial p(\rho,S)/\partial\rho$ denotes the square of the adiabatic sound speed. The coefficient $\beta\in[0,1]$ is chosen such that in the low-Mach-number limit the original splitting \eqref{eq: flux_impl_exp} is recovered for $\beta = 1$. 
As the Mach number increases beyond $\mathrm{Ma} = 0.1$ and the flow enters the compressible regime, $\beta$ decreases and the fraction of explicitly treated acoustic waves grows accordingly.
Here, we define the blending coefficient as
\begin{equation}
	\label{eq: beta_mach}
	\beta(\mathrm{Ma})=\min\left(1,\ \frac{1}{10\,\max_{\Omega}\mathrm{Ma}}\right),
\end{equation}
where $\max_{\Omega}\mathrm{Ma}$ is the maximal local Mach number at time $t^n$ over the computational domain. 
This definition provides a continuous transition between the low-Mach-number and compressible regimes. Other choices of $\beta$ are possible, see e.g.~\cite{AbateIolloPuppo}. 
Repeating the derivation of the implicit pressure equation \eqref{eq:pressure_implicit} with the blended fluxes \eqref{eq: flux_impl_exp_allmach}, the blending coefficient is inherited by the implicit mass flux and by the implicit pressure flux, so that the resulting mildly nonlinear system for the pressure reads
\begin{equation}
	\label{eq:pressure_implicit_allmach}
	\rho\left(p_i^{**},S_i^{**}\right)
	-
	\beta^2
	\left(
	\frac{\Delta t}{\Delta x}
	\right)^2
	\left(
	p^{**}_{i+1}
	-
	2p^{**}_i
	+
	p^{**}_{i-1}
	\right)
	=
	\rho_i^{*}
	-
	\beta
	\frac{\Delta t}{\Delta x}
	\left(
	m^{*}_{i+\frac12}
	-
	m^{*}_{i-\frac12}
	\right),
\end{equation}
which reduces to \eqref{eq:pressure_implicit} for $\beta=1$. Since $\beta$ is constant in space at each time step, the linear part of the system is simply multiplied by $\beta^2$ and therefore retains the properties discussed in Section~\ref{sec: Numscheme1d}.
Moreover, this setting allows one to apply the total energy correction directly to the $\mathbf{q}$ variables, which facilitates the interpretation of the correction parameter $\alpha$ as dissipation or anti-dissipation depending on its sign. 
Thus, when $\beta<1$, rather than using \eqref{eq:abgrall_flux} as correction flux, we employ
\begin{equation}
	\label{eq:abgrall_flux2}
	\widehat{f}_{i+\frac12}(\alpha)
	=
	-\alpha L
	\left(
	\mathbf{q}^{***}_{i+1}-\mathbf{q}^{***}_i
	\right).
\end{equation}
With the flux splitting \eqref{eq: flux_impl_exp_allmach}, the pressure equation \eqref{eq:pressure_implicit_allmach} and the correction flux \eqref{eq:abgrall_flux2}, all remaining steps of the semi-implicit HTC scheme are the same as described above.

\subsection{Summary of the algorithm}
\label{sec: summary}
The individual stages of the scheme have been derived one after the other in the preceding subsections. For the convenience of the reader, and of anyone who wants to implement the method, they are collected here in the order in which they are executed. All formulas are repeated without a number and each of them is accompanied by the number of the equation where it was derived, so that every step can be traced back to its origin.

\smallskip
\noindent \textbf{Explicit step.} Everything that is not tied to the acoustic waves is advanced first, with the convective flux, the entropy fluctuations and the numerical dissipation all evaluated at time $t^n$,
\begin{equation*}
\mathbf{Q}^{*}_i=\mathbf{Q}^n_i-\frac{\Delta t}{\Delta x}\left(\mathbf{f}^{E,n}_{i+\frac12}-\mathbf{f}^{E,n}_{i-\frac12}\right)-\frac{\Delta t}{\Delta x}\left(\boldsymbol{\mathcal{D}}^{-,n}_{i+\frac{1}{2}}+\boldsymbol{\mathcal{D}}^{+,n}_{i-\frac{1}{2}}\right)+\frac{\Delta t}{\Delta x}\left(\boldsymbol{\mathcal{G}}_{i+\frac{1}{2}}^n-\boldsymbol{\mathcal{G}}_{i-\frac{1}{2}}^n\right)
\end{equation*}
see \eqref{eq:expl_update_1d}, with the convective flux \eqref{eq: fluxes_expl_n}, the non-conservative fluctuations \eqref{eq: ncp_expl} and \eqref{eq: ncm_expl} and the dissipation terms \eqref{eq: dissip_expl}. The acoustic eigenvalues are absent from this stage, so that the admissible time step is governed by the material velocity alone, see \eqref{eq:cfl_condition}.

\smallskip
\noindent \textbf{Implicit step.} Inserting the discrete momentum equation into the mass balance leaves one mildly nonlinear system for the new pressure,
\begin{equation*}
\rho\left(p_i^{**},S_i^{**}\right)
-
\left(
\frac{\Delta t}{\Delta x}
\right)^2
\left(
p^{**}_{i+1}
-
2p^{**}_i
+
p^{**}_{i-1}
\right)
	=
	\rho_i^{*}
	-
	\frac{\Delta t}{\Delta x}
	\left(
	m^{*}_{i+\frac12}
	-
	m^{*}_{i-\frac12}
	\right),
\end{equation*}
see \eqref{eq:pressure_implicit}, whose nonlinearity sits on the diagonal and enters only through the equation of state, since the specific entropy has already been computed in the previous explicit step. The mildly nonlinear pressure system is solved with the nested Newton method of Casulli and Zanolli \cite{casulli_2010_a,casulli_2012_iterative}. Density and momentum on the primal grid follow from \eqref{eq:rho_n+1_star} and \eqref{eq:m_n+1_star}, the staggered momentum from \eqref{eq:mom_stagg_n+1_star}.

\smallskip
\noindent \textbf{Entropy production and Picard iteration.} The predictor obtained so far is not yet thermodynamically compatible. The entropy density assembled from the two previous stages, see \eqref{eq:rhoS_pred}, is therefore corrected by a total entropy production term,
\begin{equation*}
	\rho S_i^{***}
	=\rho S_i^{**}+ \Delta t \, \Pi_i,
	\qquad
	\Pi_i=\Pi^{s}_i+\Pi^{t}_i,
\end{equation*}
see \eqref{eq:entr_prod}, which contains a spatial entropy production term \eqref{eq: pi_s_dual} and temporal entropy production \eqref{eq: pi_t}. Both are non-negative, so that the entropy $S^{***}$ satisfies a cell entropy inequality. Since the production depends on the state that it corrects, the three stages are repeated within a Picard loop with a fixed number of iterations, two in every computation reported below, similar to \cite{casulli_2010_a}.

\smallskip
\noindent \textbf{Global energy correction.} To avoid the solution of a large nonlinear system, the cellwise compatibility of the semi-discrete scheme is replaced here by one single global constraint. The corrected state
\begin{equation*}
	\mathbf{q}_i^{n+1}(\alpha)
	=
	\mathbf{q}_i^{***}
	-
	\frac{\Delta t}{\Delta x}
	\left[
	\widehat{f}_{i+\frac12}(\alpha)
	-
	\widehat{f}_{i-\frac12}(\alpha)
	\right],
\end{equation*}
see \eqref{eq:abgrall_q_corr} with the correction flux \eqref{eq:abgrall_flux}, contains the single scalar $\alpha$, which is fixed by demanding that the discrete total energy balance \eqref{eq:global_abgrall_constraint} be satisfied globally,
\begin{equation*}
	g(\alpha)
	:=
	\sum_i \mathcal E\left(\mathbf{q}_i^{n+1}(\alpha)\right)\Delta x
	-
	\sum_i \mathcal E(\mathbf{q}_i^n)\Delta x
	+
	\Delta t\left(F(\mathbf{q}_R)-F(\mathbf{q}_L)\right)=0.
\end{equation*}
 This scalar equation is solved by the Newton method with the initial guess $\alpha = 0$.

\smallskip
\noindent \textbf{P-split variant.} The flux is split with a blending coefficient $0 \leq \beta \leq 1$,
\begin{equation*}
	\mathbf{f}_k^{E}(\mathbf{Q})=\bigl((1-\beta)\rho u_k,\,
	\rho u_lu_k+(1-\beta)p\,\delta_{lk},\,0\bigr)^\top,
	\qquad
	\mathbf{f}_k^{I}(\mathbf{Q})=\bigl(\beta\rho u_k,\,
	\beta p\,\delta_{lk},\,0\bigr)^\top,
\end{equation*}
see \eqref{eq: flux_impl_exp_allmach}, with the blending coefficient $\beta$ as chosen in \eqref{eq: beta_mach}. The discrete pressure wave equation of the p-split scheme becomes
\begin{equation*}
	\rho\left(p_i^{**},S_i^{**}\right)
	-
	\beta^2
	\left(
	\frac{\Delta t}{\Delta x}
	\right)^2
	\left(
	p^{**}_{i+1}
	-
	2p^{**}_i
	+
	p^{**}_{i-1}
	\right)
	=
	\rho_i^{*}
	-
	\beta
	\frac{\Delta t}{\Delta x}
	\left(
	m^{*}_{i+\frac12}
	-
	m^{*}_{i-\frac12}
	\right),
\end{equation*}
see \eqref{eq:pressure_implicit_allmach}. The four stages are then executed as above, with the characteristic speeds \eqref{eq:expl_charSpeeds_allmach} entering the CFL condition and with the global corrector applied to the conservative variables, see \eqref{eq:abgrall_flux2}.

\section{Numerical results}
\label{sec: Numres}
In this section, we present several numerical test cases aiming at assessing the behavior of the proposed new semi-implicit HTC scheme for the Euler equations. All tests are performed on uniform Cartesian grids.
The time step is computed according to the CFL conditions \eqref{eq:cfl_condition} and
\begin{equation}
	\label{eq:cfl_condition2} 
	\Delta t=
	\mathrm{CFL} \min \limits_{i,j} 
	\left(
	\frac{ \max \limits_r |\lambda_{x,r}^E(\mathbf{Q}_{i,j}^n)|}{\Delta x}
	+
	\frac{\displaystyle\max_r |\lambda_{y,r}^E(\mathbf{Q}_{i,j}^n)|}{\Delta y}
	\right)^{-1},
\end{equation}
for the one- and the two-dimensional test cases, respectively. 
 The numerical scheme described in Section~\ref{sec: Numscheme1d} is denoted by SIHTC, while the p-split version of Section~\ref{sec:allmachext} is denoted by SIHTC p-split. 
 For the SIHTC scheme, the CFL number depends on the test case and for accuracy reasons it is reduced in the presence of shocks. By contrast, the SIHTC p-split formulation, which treats part of the pressure terms explicitly, mitigates this restriction and remains stable and accurate with $\mathrm{CFL}=0.9$ for all simulations considered in this work. 

 In both formulations, the explicit stability restriction depends only on the characteristic speeds of the explicit subsystem and is therefore significantly less restrictive than in a fully explicit discretization.
For all test cases, we employ the ideal gas equation of state \eqref{eq: energy_eos} with $\gamma=1.4$ and $c_v=1.0$. The initial data are prescribed in terms of the primitive variables, from which the specific entropy evolved by the scheme follows as
\begin{equation}
	\label{eq: entropy_from_primitives}
	S = c_v \log\left( \frac{p}{\rho^\gamma} \right).
\end{equation}
\subsection{One-dimensional Riemann problems}
\label{sec:1D_RP}
Here we consider three Riemann problems (RPs) for the 
Euler equations: first the Sod shock tube (RP1)~\cite{SOD19781} and the Lax problem (RP2)~\cite{lax2} which have both become standard benchmarks for 
numerical methods for the Euler equations in general. The solution in both cases consists of three waves: a left-moving rarefaction fan, an intermediate contact discontinuity, and a 
right-moving shock wave. 
Additionally, we consider a low-Mach-number Riemann problem (RP3) on which we test the behavior of the numerical scheme in case of a slowly moving contact at low Mach number.
The computational domain for all the Riemann problems is $\Omega = [-0.5, 0.5]$, 
with the initial discontinuity located at $x_0 = 0$ and the fixed Dirichlet boundary values $\mathbf{q}_L$ and $\mathbf{q}_R$, which are constant in time, imposed at the left and at the right boundary. 
In this section, we compare the numerical results of the SIHTC and SIHTC p-split schemes to the exact solution of the Riemann problem, the explicit local Lax--Friedrichs (LLF) or Rusanov scheme, and the HTC scheme introduced in Section~\ref{sec: semidiscHTC} equipped with a fifth-order explicit Runge--Kutta time integration. 

\subsubsection{RP1: Sod shock tube}
\label{sec:sod}

\begin{figure}[htpb]
	\centering
	\begin{subfigure}[b]{0.48\textwidth}
		\includegraphics[width=\textwidth]{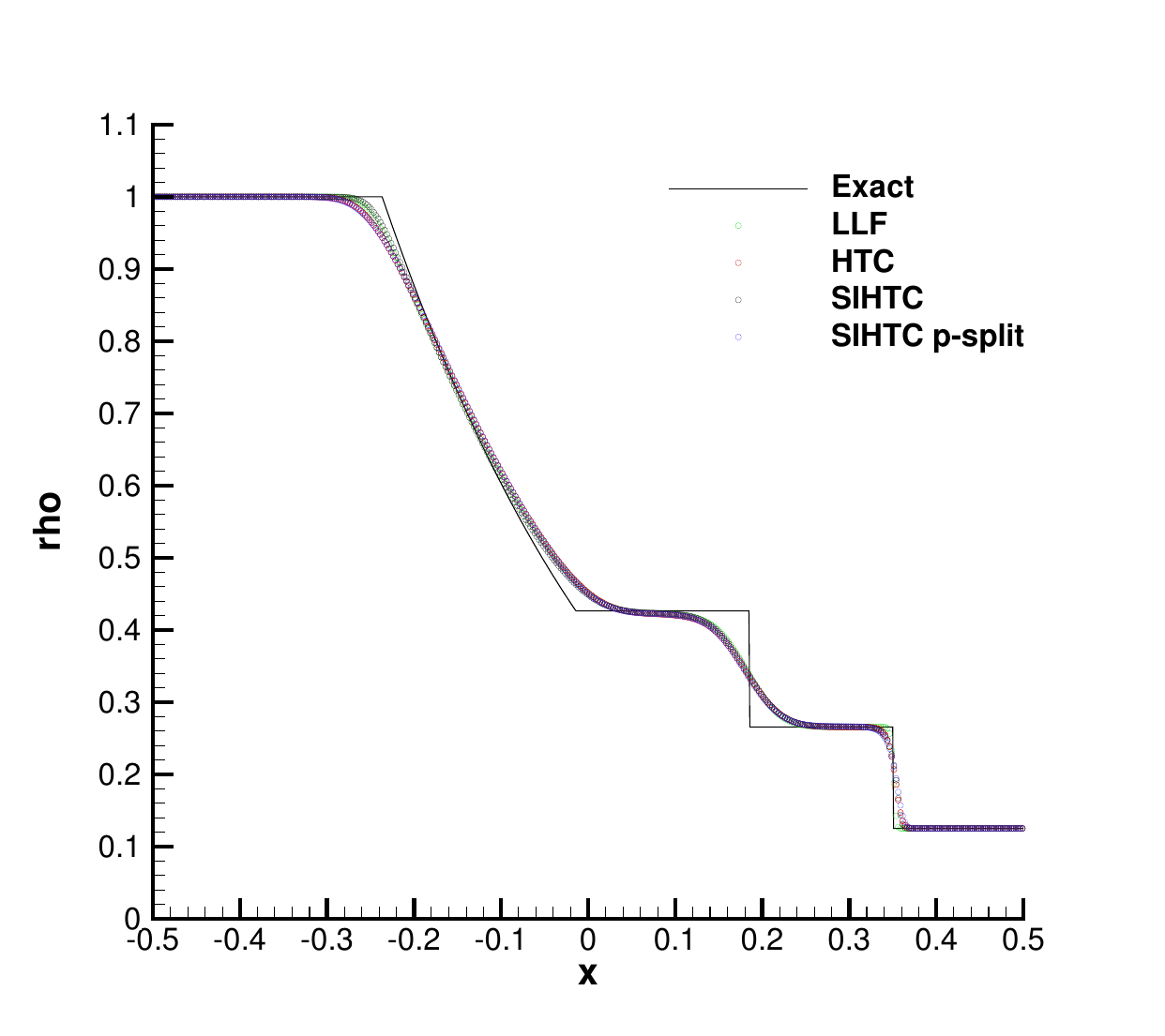}
		\caption{Density $\rho$}
	\end{subfigure}
	\hfill
	\begin{subfigure}[b]{0.48\textwidth}
		\includegraphics[width=\textwidth]{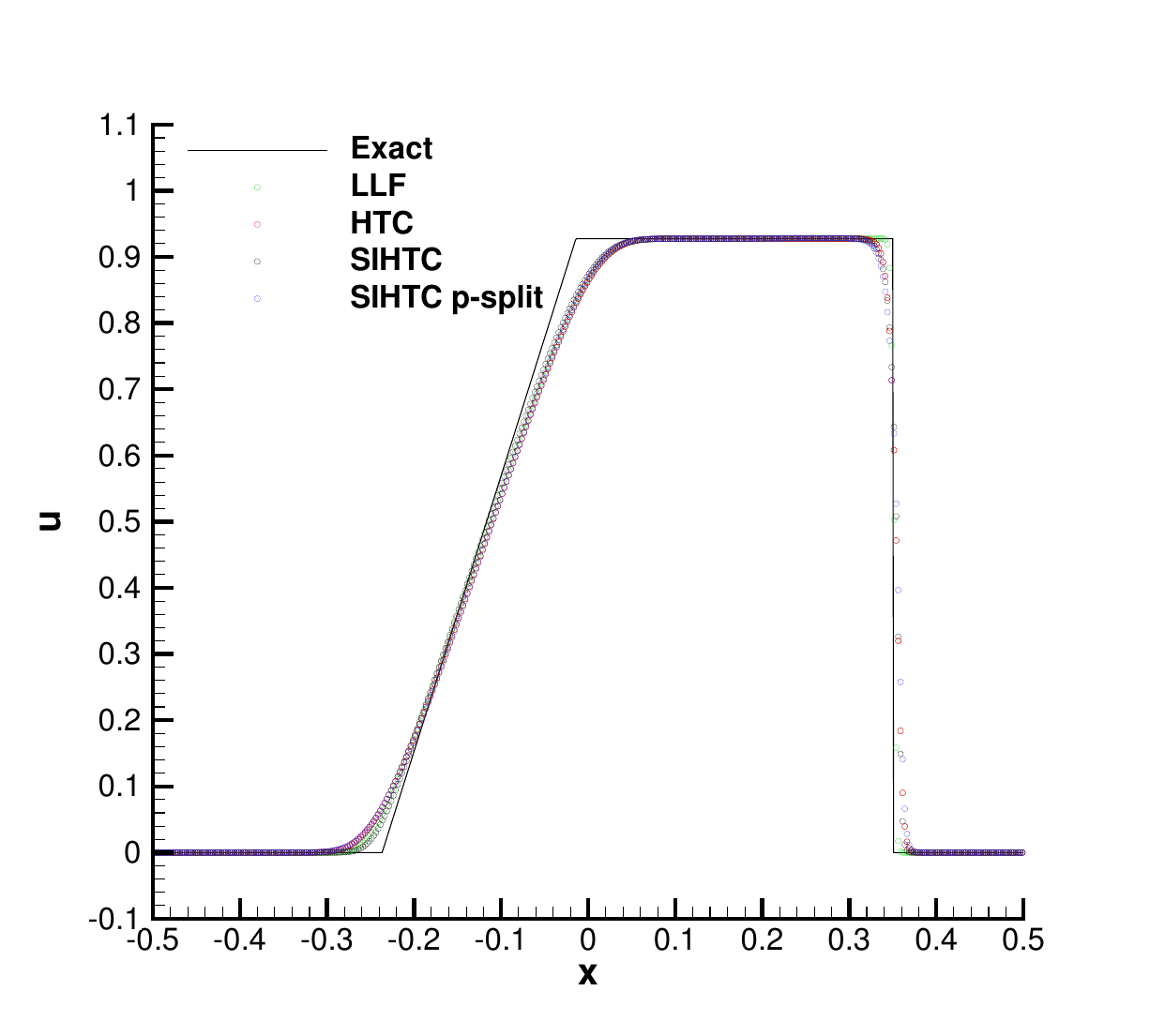}
		\caption{Velocity $u$}
	\end{subfigure}
	
	\vspace{1ex}
	
	\begin{subfigure}[b]{0.48\textwidth}
		\includegraphics[width=\textwidth]{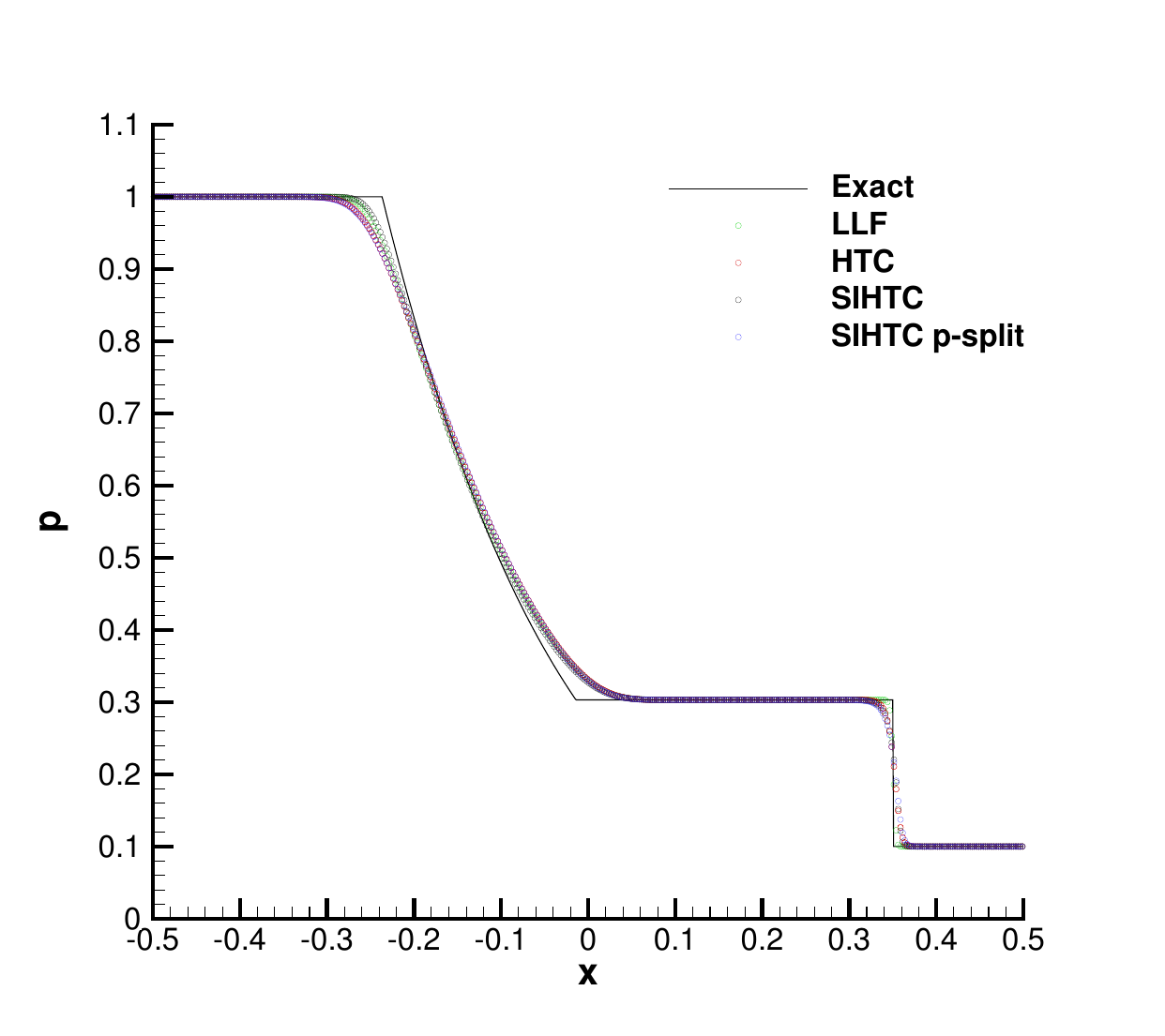}
		\caption{Pressure $p$}
	\end{subfigure}
	\hfill
	\begin{subfigure}[b]{0.48\textwidth}
		\includegraphics[width=\textwidth]{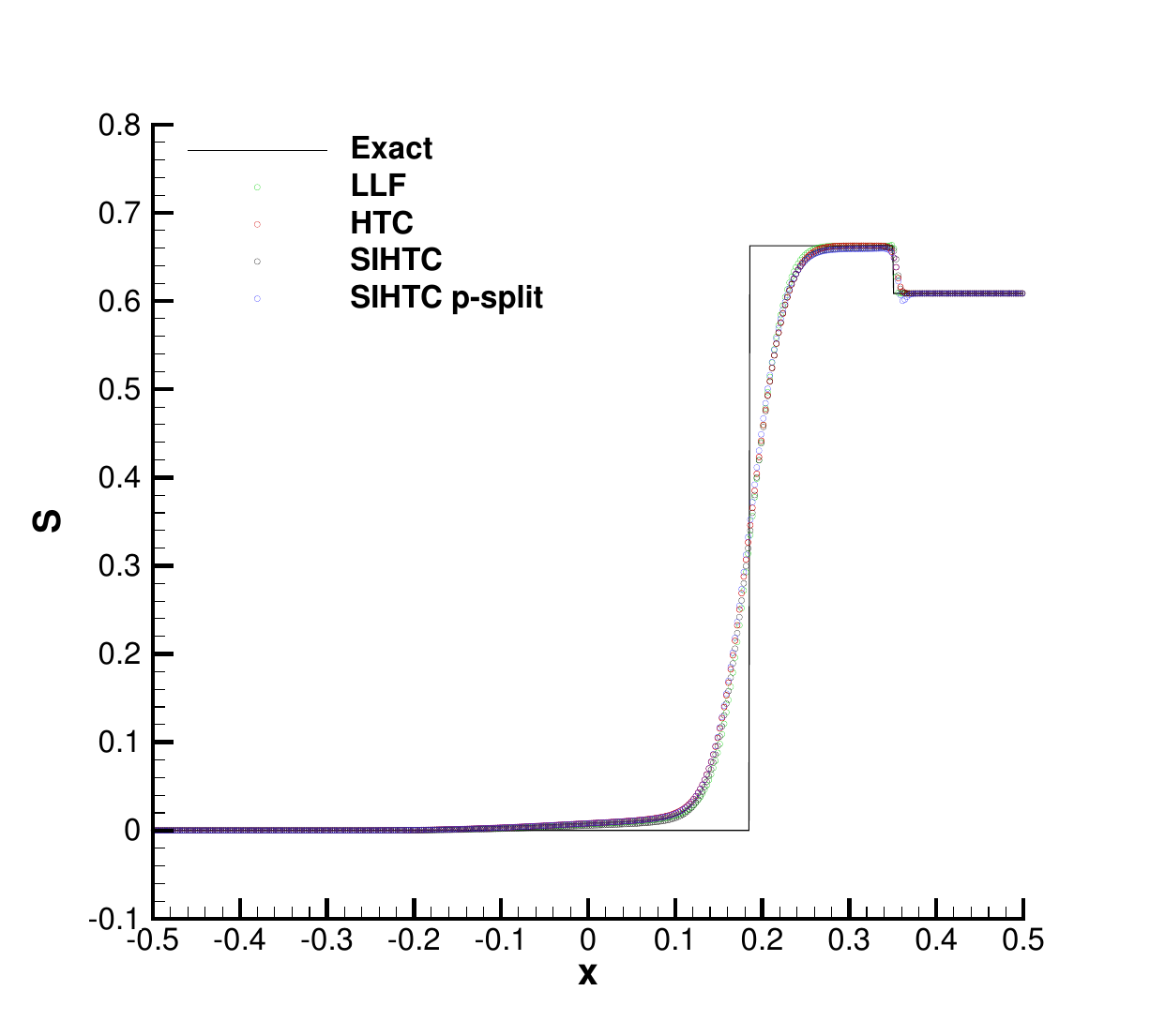}
		\caption{Entropy $S$}
	\end{subfigure}
	\caption{RP1 -- Sod: numerical solutions at $t = 0.2$, $N = 400$ cells. The CFL number is $0.45$ for SIHTC and $0.9$ for all other numerical schemes.}
	\label{fig:sod_solution}
\end{figure}

The first test case that we show here is the classical Sod shock tube problem~\cite{SOD19781}, whose initial condition reads
\begin{equation}
	(\rho, u, p)(x,0) =
	\left\{
	\begin{array}{@{}lll@{\quad}l@{}}
		(1,     & 0, & 1)   & \text{if } x < 0,\\
		(0.125, & 0, & 0.1) & \text{if } x \geq 0.
	\end{array}
	\right.
\end{equation}
The solution is computed up to the final time $t = 0.2$ with $\mathrm{CFL} = 0.45$ for the SIHTC scheme.
Since this test is situated in the compressible regime, the SIHTC p-split scheme is also considered with $\mathrm{CFL}=0.9$.
Table~\ref{tab:sod_convergence} reports the $L^1$ errors and convergence rates for the SIHTC scheme with respect to the exact solution, computed with the Riemann solver from \cite{ToroBook}, for the density $\rho$, the velocity $u$, the pressure $p$, and the specific entropy $S$
obtained by refining the mesh from $N=50$ to $N=3200$ cells. 
Table~\ref{tab:sod_convergence_psplit} reports the corresponding convergence results for the SIHTC p-split scheme.
Both formulations achieve the expected experimental order of convergence in the $L^1$ norm, with rates stabilizing around $0.5$ for the specific entropy, around $0.65$ for the density and around $0.8$ for the velocity and the pressure, 
consistent with the theoretical expectations for first-order schemes in the presence of discontinuities \cite{ToroBook}. 
The SIHTC p-split formulation exhibits convergence rates comparable to those of the SIHTC scheme.
In Figure~\ref{fig:sod_solution} the numerical results of both new semi-implicit HTC schemes are compared with a first-order finite volume scheme based on the LLF flux, with the semi-discrete HTC scheme and with the exact solution.
We observe a good agreement with respect to shock positions and amplitudes, as well as the increase of the specific entropy at the shock.

\begin{table}[htpb]
	\centering
	\caption{RP1 -- Sod: $L^1$ errors and convergence rates at $t=0.2$ for the SIHTC scheme.}
	\label{tab:sod_convergence}
	\begin{tabular}{c cc cc cc cc}
		\hline
		$N$
		& $L^1(\rho)$ & rate
		& $L^1(u)$    & rate
		& $L^1(p)$    & rate
		& $L^1(S)$    & rate \\
		\hline
		50    & 2.967e-02 & --   & 7.161e-02 & --   & 3.081e-02 & --   & 5.074e-02 & --   \\
		100   & 2.225e-02 & 0.41 & 4.434e-02 & 0.69 & 1.921e-02 & 0.68 & 3.845e-02 & 0.40 \\
		200   & 1.499e-02 & 0.57 & 2.530e-02 & 0.81 & 1.138e-02 & 0.76 & 2.757e-02 & 0.48 \\
		400   & 9.747e-03 & 0.62 & 1.418e-02 & 0.84 & 6.596e-03 & 0.79 & 1.959e-02 & 0.49 \\
		800   & 6.265e-03 & 0.64 & 7.922e-03 & 0.84 & 3.771e-03 & 0.81 & 1.379e-02 & 0.51 \\
		1600  & 4.012e-03 & 0.64 & 4.360e-03 & 0.86 & 2.136e-03 & 0.82 & 9.618e-03 & 0.52 \\
		3200  & 2.601e-03 & 0.63 & 2.536e-03 & 0.78 & 1.215e-03 & 0.81 & 6.838e-03 & 0.49 \\
		\hline
	\end{tabular}
\end{table}
\begin{table}[htpb]
	\centering
	\caption{RP1 -- Sod: $L^1$ errors and convergence rates at $t=0.2$ for the SIHTC p-split scheme.}
	\label{tab:sod_convergence_psplit}
	\begin{tabular}{c cc cc cc cc}
		\hline
		$N$
		& $L^1(\rho)$ & rate
		& $L^1(u)$    & rate
		& $L^1(p)$    & rate
		& $L^1(S)$    & rate \\
		\hline
		50    & 3.731e-02 & --   & 9.252e-02 & --   & 4.161e-02 & --   & 5.331e-02 & --   \\
		100   & 2.777e-02 & 0.43 & 5.756e-02 & 0.68 & 2.658e-02 & 0.65 & 4.070e-02 & 0.39 \\
		200   & 1.868e-02 & 0.57 & 3.318e-02 & 0.79 & 1.611e-02 & 0.72 & 2.914e-02 & 0.48 \\
		400   & 1.213e-02 & 0.62 & 1.879e-02 & 0.82 & 9.519e-03 & 0.76 & 2.072e-02 & 0.49 \\
		800   & 7.755e-03 & 0.65 & 1.056e-02 & 0.83 & 5.533e-03 & 0.78 & 1.455e-02 & 0.51 \\
		1600  & 4.907e-03 & 0.66 & 5.705e-03 & 0.89 & 3.148e-03 & 0.81 & 1.011e-02 & 0.52 \\
		3200  & 3.111e-03 & 0.66 & 3.232e-03 & 0.82 & 1.785e-03 & 0.82 & 7.022e-03 & 0.53 \\
		\hline
	\end{tabular}
\end{table}
\subsubsection{RP2: Lax shock tube}
\label{sec:lax}
In this test we compute the Lax problem~\cite{lax2} with initial data
\begin{equation}
	(\rho, u, p)(x,0) =
	\left\{
	\begin{array}{@{}lll@{\quad}l@{}}
		(0.445, & 0.698, & 3.528) & \text{if } x < 0,\\
		(0.5,   & 0,     & 0.571) & \text{if } x \geq 0.
	\end{array}
	\right.
\end{equation}
The final time is $t = 0.14$ with $\mathrm{CFL} = 0.1$ for the SIHTC scheme, while for the SIHTC p-split scheme a nine times larger CFL number is used. 
Tables~\ref{tab:lax_convergence} and \ref{tab:lax_convergence_psplit} report the $L^1$ errors and convergence rates for the SIHTC and the SIHTC p-split scheme, respectively, with respect to the exact solution of RP2, for the density $\rho$, the velocity $u$, the pressure $p$ and the specific entropy $S$. 
Also for RP2, the experimental order of convergence is consistent with the theoretical expectations for a first-order scheme on a solution containing strong shocks. The SIHTC p-split formulation yields errors and convergence rates comparable to those of the SIHTC scheme, even though larger time steps are applied.  

\begin{figure}[htpb]
	\centering
	\begin{subfigure}[b]{0.48\textwidth}
		\includegraphics[width=\textwidth]{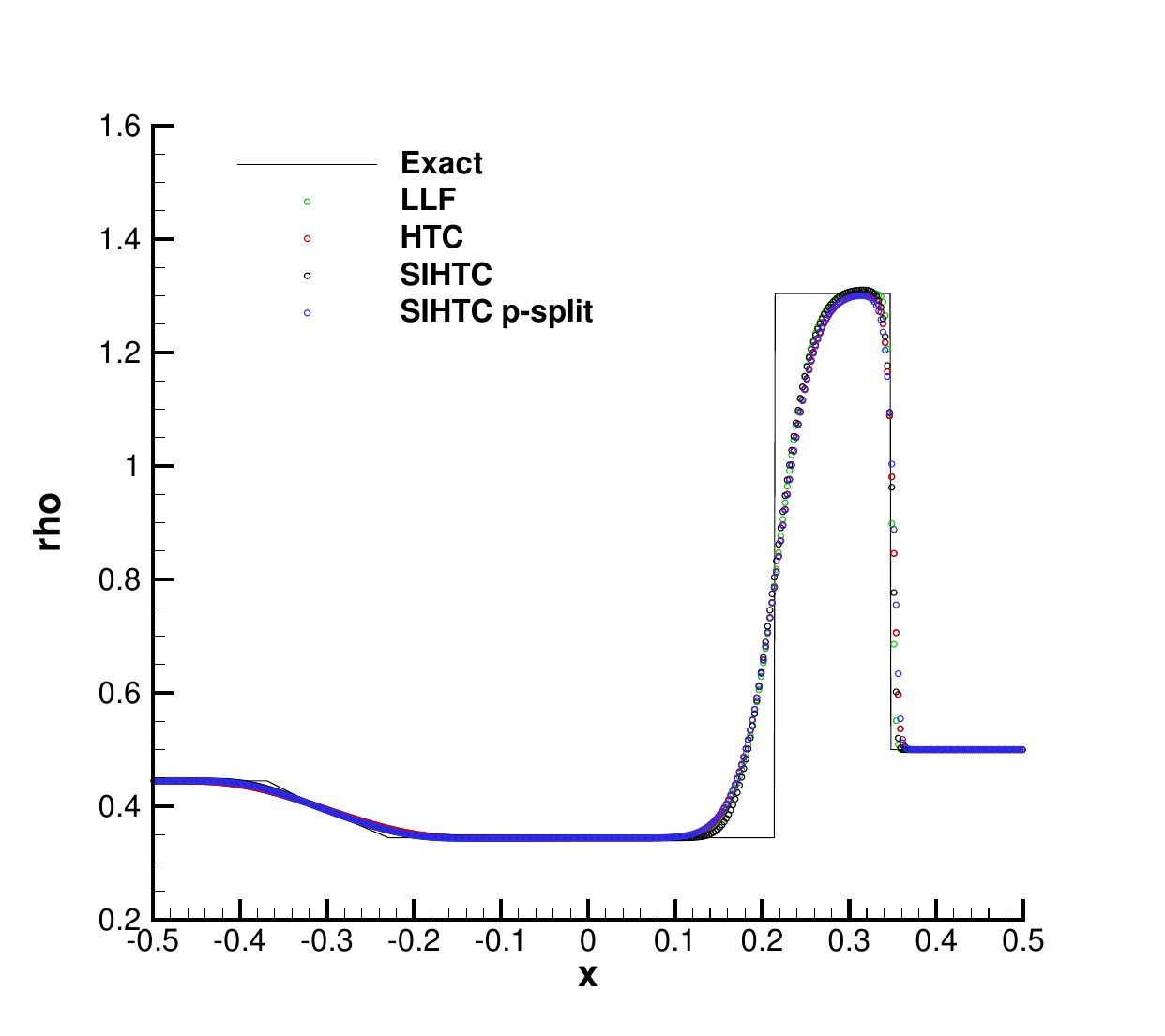}
		\caption{Density $\rho$}
	\end{subfigure}
	\hfill
	\begin{subfigure}[b]{0.48\textwidth}
		\includegraphics[width=\textwidth]{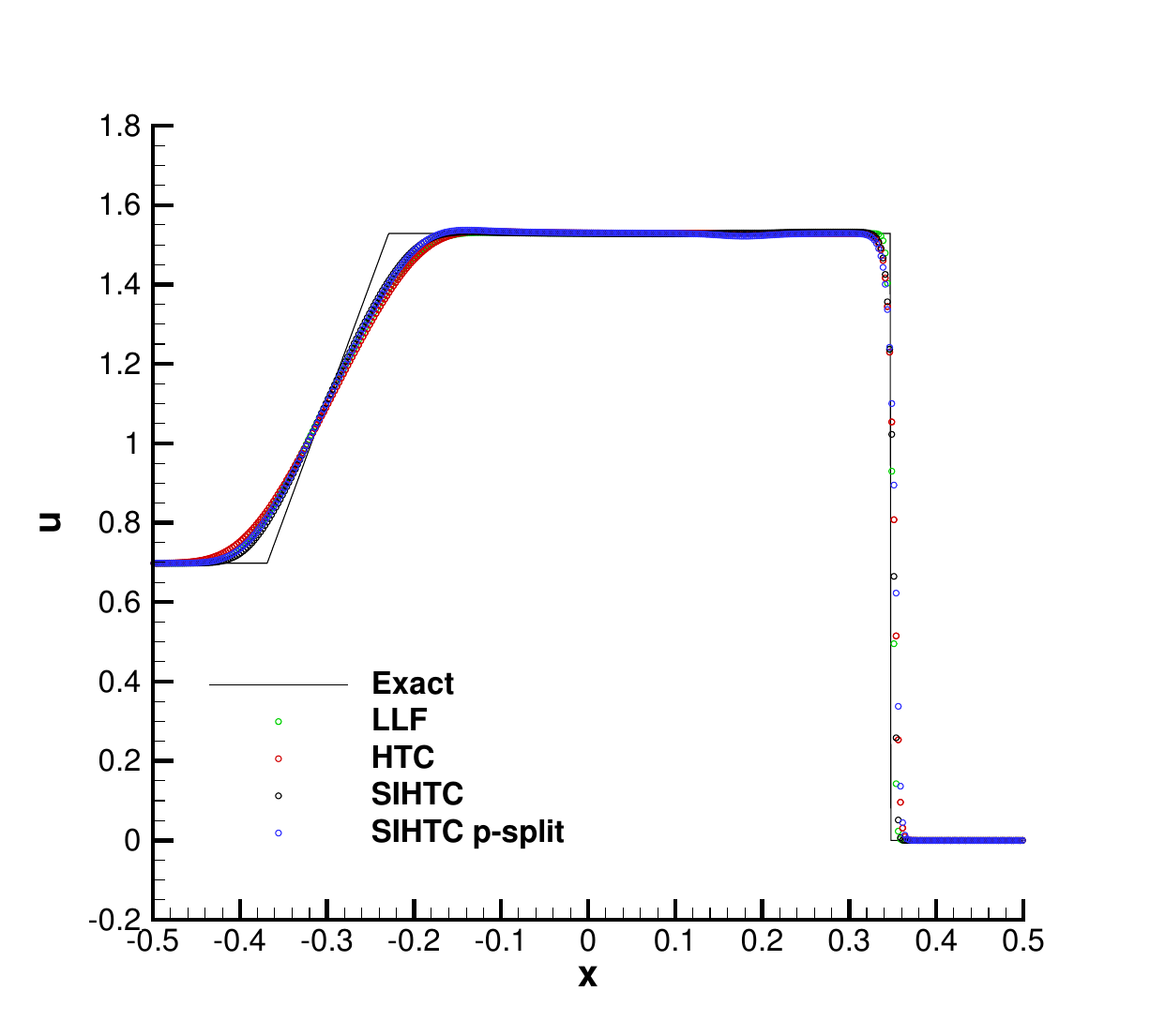}
		\caption{Velocity $u$}
	\end{subfigure}
	
	\vspace{1ex}
	
	\begin{subfigure}[b]{0.48\textwidth}
		\includegraphics[width=\textwidth]{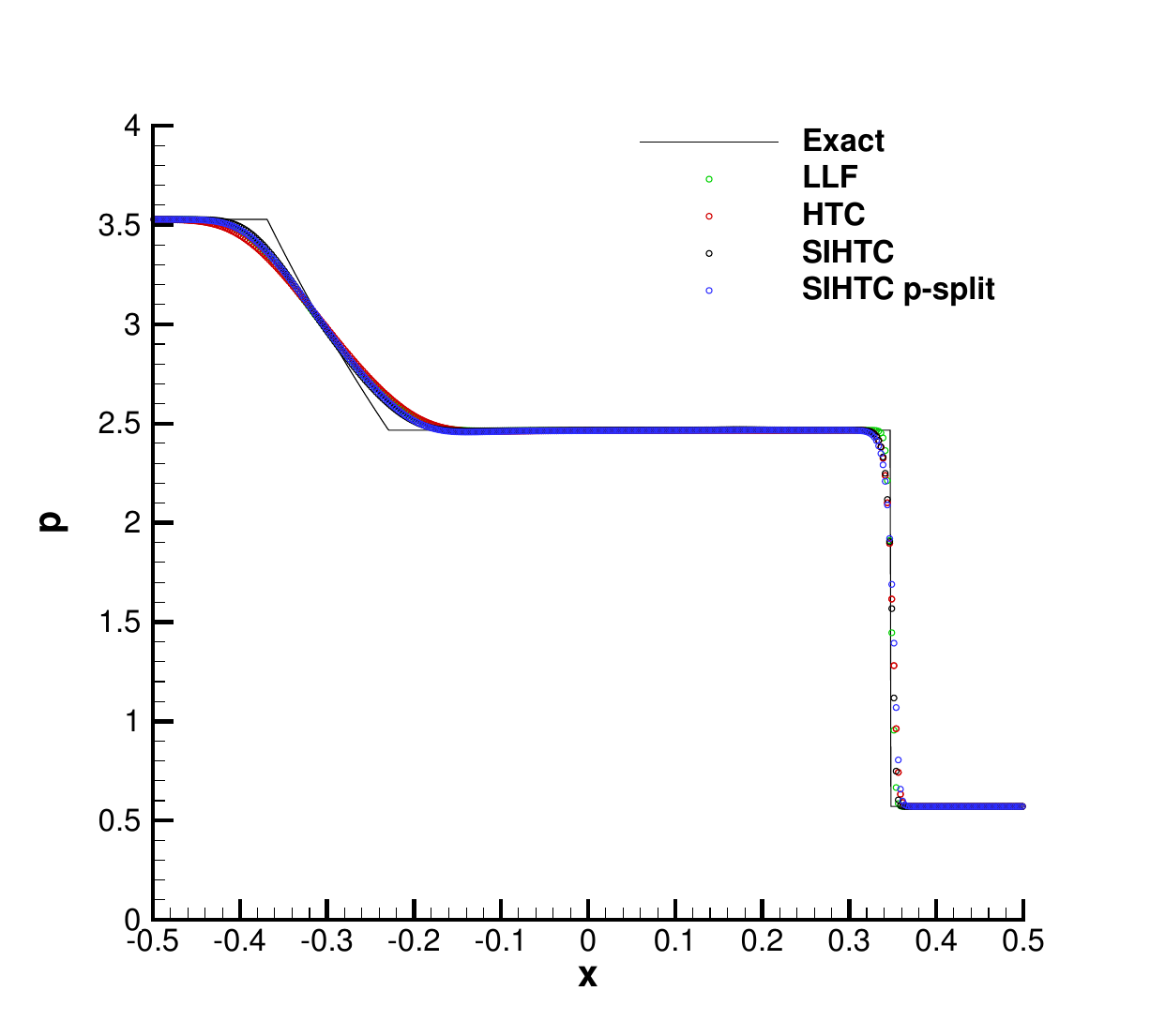}
		\caption{Pressure $p$}
	\end{subfigure}
	\hfill
	\begin{subfigure}[b]{0.48\textwidth}
		\includegraphics[width=\textwidth]{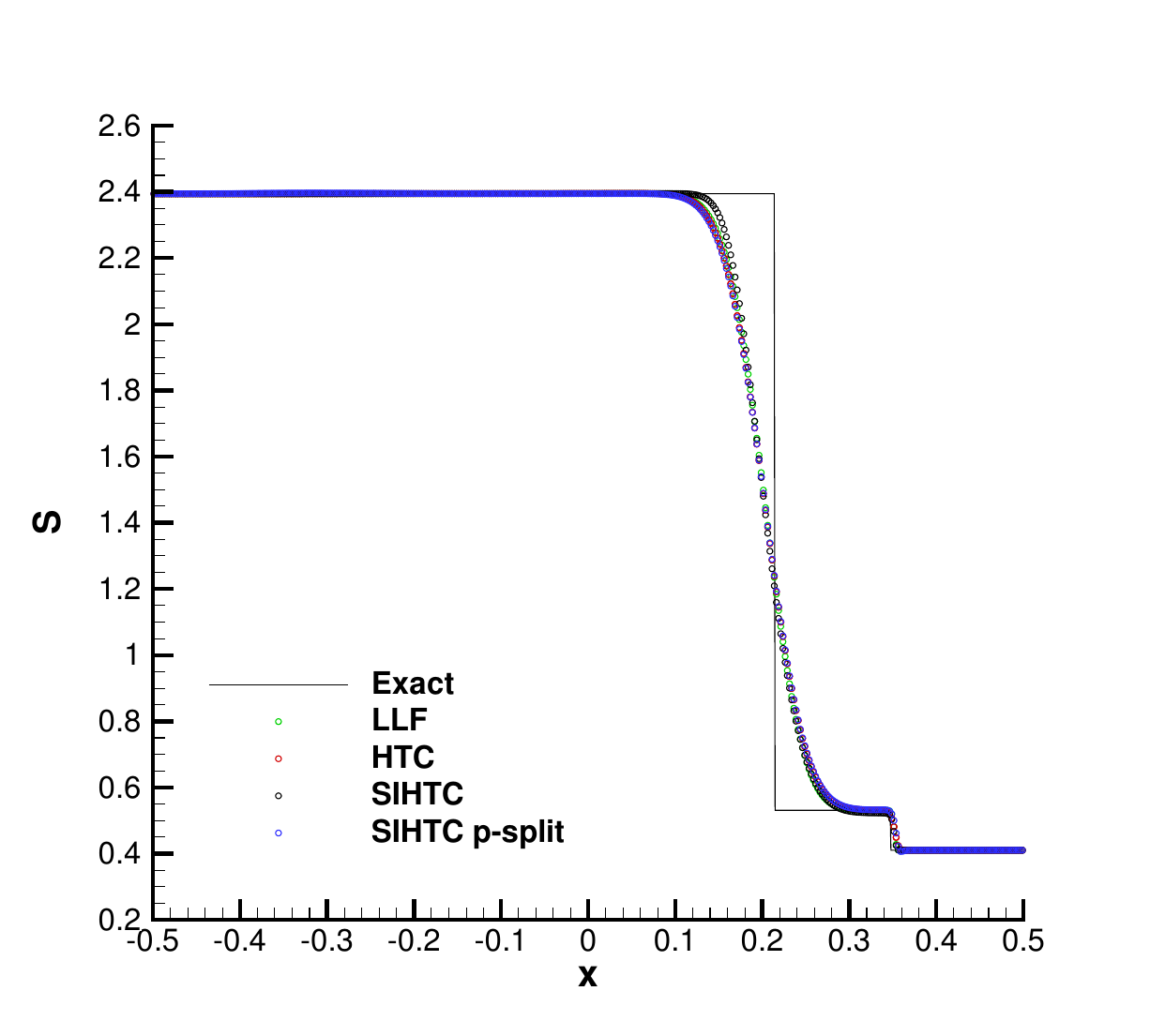}
		\caption{Entropy $S$}
	\end{subfigure}
	\caption{RP2 -- Lax: numerical solutions
		at $t = 0.14$, $N = 400$ cells. The CFL number is $0.1$ for SIHTC and $0.9$ for all other numerical schemes.}
	\label{fig:lax_solution}
\end{figure}

\begin{figure}[htbp]
	\centering
	\includegraphics[width=0.48\textwidth]{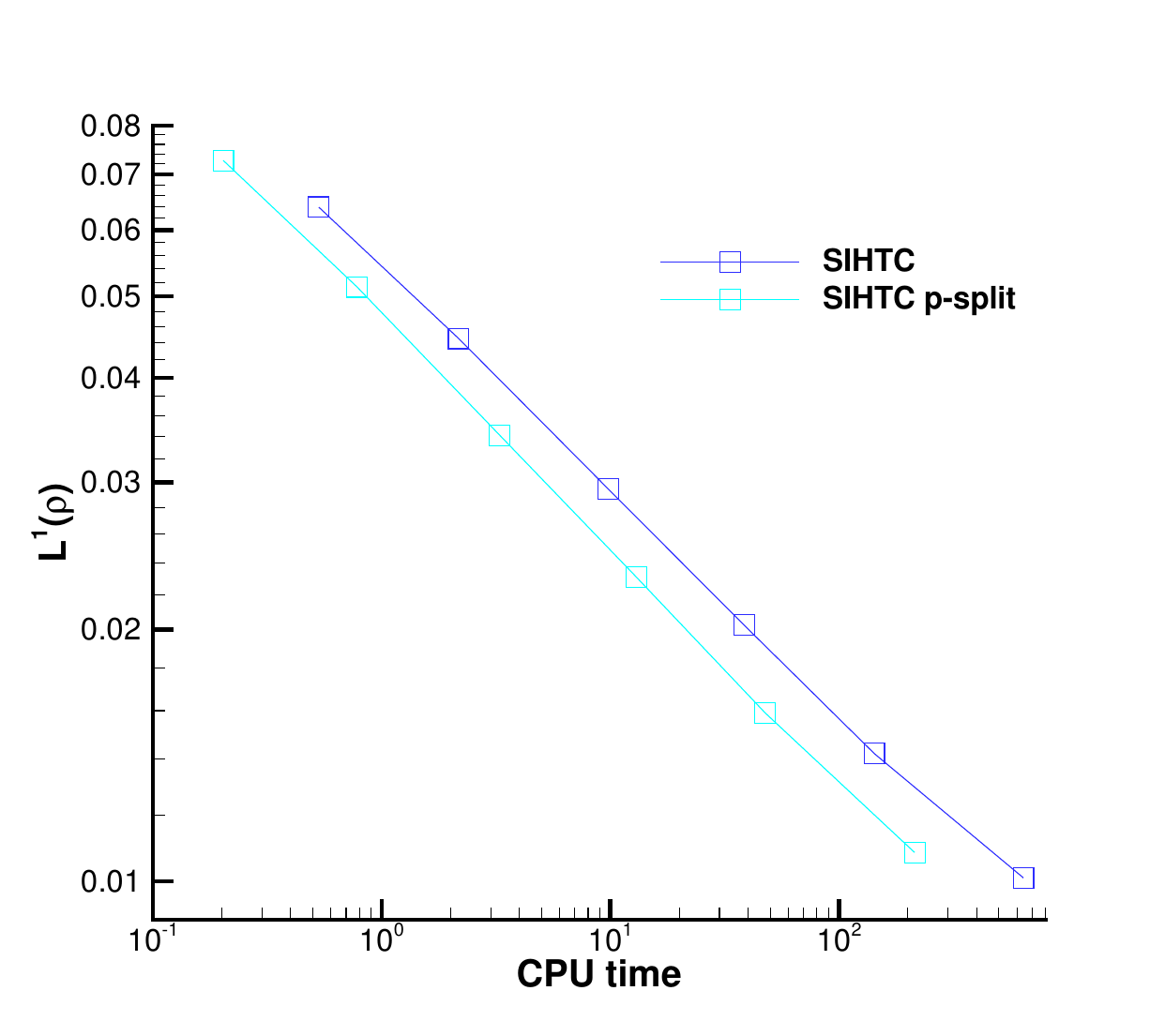}
	\caption{RP2 -- Lax: comparison of the $L^1(\rho)$ errors of the SIHTC and the SIHTC p-split scheme as a function of the CPU time in seconds.}
	\label{fig:cputime}
\end{figure}

In Figure~\ref{fig:lax_solution}, the numerical results are shown. The two formulations achieve the right shock position and amplitudes before and after the shock for all variables. 
\begin{table}[htpb]
	\centering
	\caption{RP2 -- Lax: $L^1$ errors and convergence rates at $t=0.14$ for the SIHTC scheme.}
	\label{tab:lax_convergence}
	\begin{tabular}{c cc cc cc cc}
		\hline
		$N$
		& $L^1(\rho)$ & rate
		& $L^1(u)$    & rate
		& $L^1(p)$    & rate
		& $L^1(S)$    & rate \\
		\hline
		50    & 9.178e-02 & --   & 9.751e-02 & --   & 1.128e-01 & --   & 1.396e-01 & --   \\
		100   & 6.394e-02 & 0.52 & 5.495e-02 & 0.83 & 6.691e-02 & 0.75 & 9.790e-02 & 0.51 \\
		200   & 4.455e-02 & 0.52 & 3.590e-02 & 0.61 & 4.237e-02 & 0.66 & 7.304e-02 & 0.42 \\
		400   & 2.942e-02 & 0.60 & 1.926e-02 & 0.90 & 2.359e-02 & 0.85 & 5.236e-02 & 0.48 \\
		800   & 2.025e-02 & 0.54 & 1.117e-02 & 0.79 & 1.375e-02 & 0.78 & 3.672e-02 & 0.51 \\
		1600  & 1.422e-02 & 0.51 & 6.978e-03 & 0.68 & 8.302e-03 & 0.73 & 2.594e-02 & 0.50 \\
		3200  & 1.010e-02 & 0.49 & 4.022e-03 & 0.79 & 4.953e-03 & 0.75 & 1.870e-02 & 0.47 \\
		\hline
	\end{tabular}
\end{table}
\begin{table}[htpb]
	\centering
	\caption{RP2 -- Lax: $L^1$ errors and convergence rates at $t=0.14$ for the SIHTC p-split scheme.}
	\label{tab:lax_convergence_psplit}
	\begin{tabular}{c cc cc cc cc}
		\hline
		$N$
		& $L^1(\rho)$ & rate
		& $L^1(u)$    & rate
		& $L^1(p)$    & rate
		& $L^1(S)$    & rate \\
		\hline
		50    & 9.766e-02 & --   & 1.054e-01 & --   & 1.203e-01 & --   & 1.548e-01 & --   \\
		100   & 7.271e-02 & 0.43 & 6.725e-02 & 0.65 & 8.043e-02 & 0.58 & 1.109e-01 & 0.48 \\
		200   & 5.133e-02 & 0.50 & 4.514e-02 & 0.58 & 5.226e-02 & 0.62 & 8.282e-02 & 0.42 \\
		400   & 3.415e-02 & 0.59 & 2.467e-02 & 0.87 & 2.917e-02 & 0.84 & 5.954e-02 & 0.48 \\
		800   & 2.310e-02 & 0.56 & 1.375e-02 & 0.84 & 1.645e-02 & 0.83 & 4.198e-02 & 0.50 \\
		1600  & 1.590e-02 & 0.54 & 8.093e-03 & 0.76 & 9.564e-03 & 0.78 & 2.982e-02 & 0.49 \\
		3200  & 1.083e-02 & 0.55 & 4.337e-03 & 0.90 & 5.229e-03 & 0.87 & 2.121e-02 & 0.49 \\
		\hline
	\end{tabular}
\end{table}
While the SIHTC p-split scheme allows a higher CFL number, its explicit eigenvalues are also larger in magnitude than those of the SIHTC scheme, so that the two effects partially compensate each other and the gain in the admissible time step is smaller than the ratio of the CFL numbers alone would suggest.
To quantify the advantage of applying the SIHTC p-split, Table~\ref{tab:sihtc_psplit_comparison} compares the computational performance of the SIHTC and SIHTC p-split schemes for the Lax Riemann problem. The SIHTC p-split admits an average time step approximately three times larger than the one of the SIHTC scheme. As a result, the total CPU time is reduced by roughly a factor of three across all grid resolutions, with this computational advantage remaining essentially independent of the mesh size.
Comparing the error versus CPU times in Figure~\ref{fig:cputime} shows that the proposed p-split version significantly improves the computational efficiency of the SIHTC formulation in compressible regimes without compromising its accuracy.

\begin{table}[htbp]
	\centering
	\caption{RP2 -- Lax: computational cost comparison between SIHTC and SIHTC p-split. The CPU times were measured with a timer resolution of $1/64$ s.}
	\label{tab:sihtc_psplit_comparison}
	\begin{tabular}{rrrrrrr}
		\toprule
		$N$ & $T_{\mathrm{SIHTC}}$ [s] & $T_{\mathrm{p-split}}$ [s] & $T_{\mathrm{p-split}}/T_{\mathrm{SIHTC}}$ & $\overline{\Delta t}_{\mathrm{SIHTC}}$ & $\overline{\Delta t}_{\mathrm{p-split}}$ & $\overline{\Delta t}_{\mathrm{p-split}}/\overline{\Delta t}_{\mathrm{SIHTC}}$ \\
		\midrule
		
		100  & 0.53   & 0.20  & 0.38 & $3.2887 \cdot 10^{-4}$ & $1.0000 \cdot 10^{-3}$ & 3.04 \\
		200  & 2.17   & 0.78  & 0.36 & $1.6361 \cdot 10^{-4}$ & $5.1664 \cdot 10^{-4}$ & 3.16 \\
		400  & 9.83   & 3.28  & 0.33 & $8.1667 \cdot 10^{-5}$ & $2.5673 \cdot 10^{-4}$ & 3.14 \\
		800  & 38.45  & 13.02 & 0.34 & $4.0811 \cdot 10^{-5}$ & $1.2800 \cdot 10^{-4}$ & 3.14 \\
		1600 & 142.59 & 47.58 & 0.33 & $2.0400 \cdot 10^{-5}$ & $6.3930 \cdot 10^{-5}$ & 3.13 \\
		3200 & 639.95 & 214.13& 0.33 & $1.0198 \cdot 10^{-5}$ & $3.1953 \cdot 10^{-5}$ & 3.13 \\
		\bottomrule
	\end{tabular}
\end{table}

\subsubsection{RP3: Low-Mach Riemann problem}
\label{sec:rp_lowmach}
To assess the behavior of the scheme in the low-Mach-number regime,
we consider a one-dimensional Riemann problem in which the initial discontinuity involves only velocity and pressure perturbations of order $\mathrm{Ma}$,
while the density is uniform. The initial data read
\begin{equation}
	(\rho, u, p)(x,0) =
	\left\{
	\begin{array}{@{}lll@{\quad}l@{}}
		(1, & 0,           & 1)        & \text{if } x < 0,\\
		(1, & 8\,\mathrm{Ma}, & 1-\mathrm{Ma}) & \text{if } x \geq 0,
	\end{array}
	\right.
\end{equation}
with $\mathrm{Ma} = 10^{-3}$. This configuration generates two acoustic waves of
amplitude $\mathcal{O}(\mathrm{Ma})$ propagating approximately at speeds $\pm c_0$, with
$c_0 = \sqrt{\gamma p_0/\rho_0}$, $\rho_0 = 1$ and $p_0 = 1$. 
Since this problem belongs entirely to the low-Mach-number regime, the blending coefficient satisfies $\beta=1$ throughout the simulation. Therefore, the p-split version coincides with the original SIHTC scheme and no separate SIHTC p-split results are reported.

\begin{figure}[htpb]
	\centering
	\begin{subfigure}[b]{0.48\textwidth}
		\includegraphics[width=\textwidth]{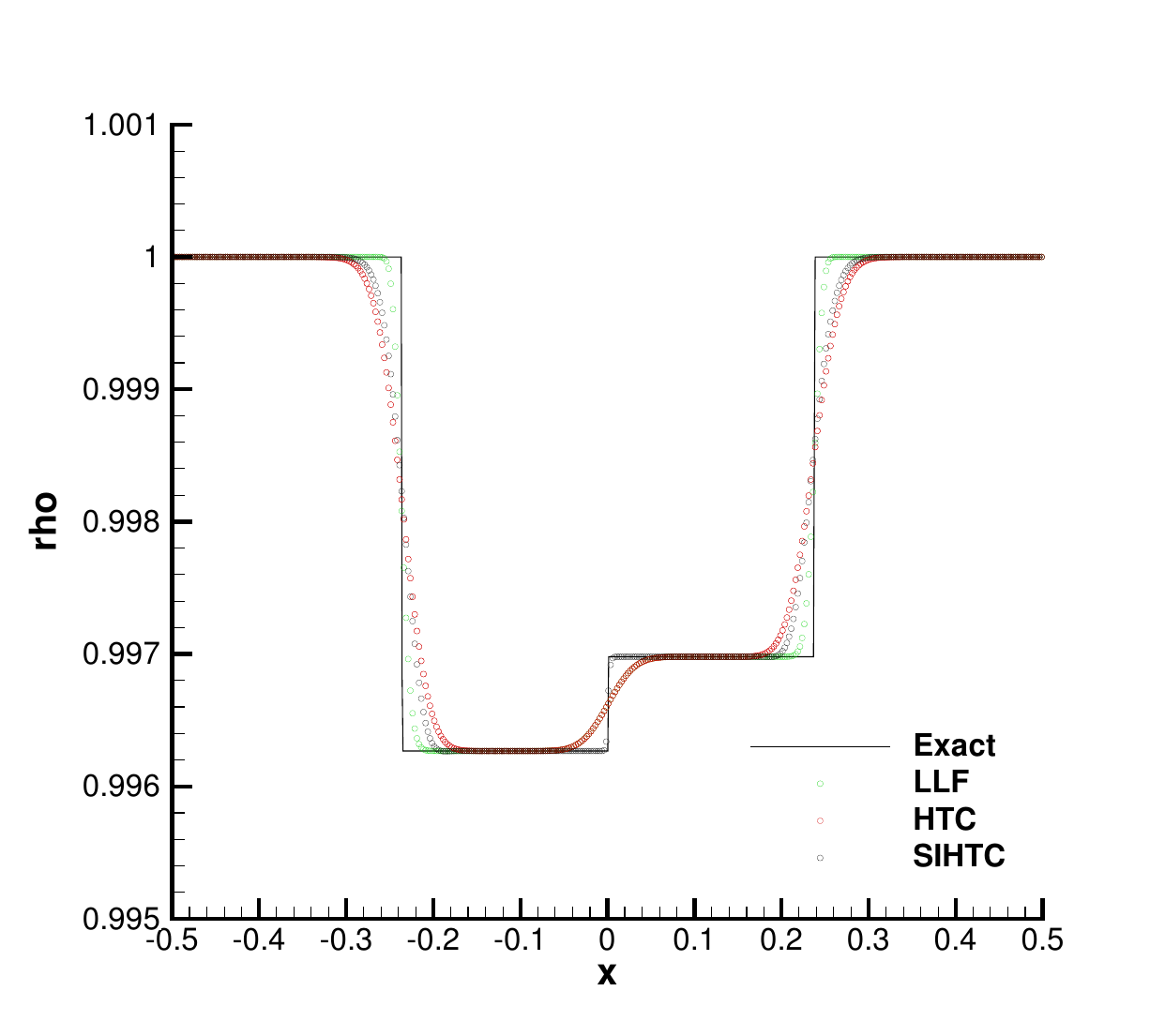}
		\caption{Density $\rho$}
	\end{subfigure}
	\hfill
	\begin{subfigure}[b]{0.48\textwidth}
		\includegraphics[width=\textwidth]{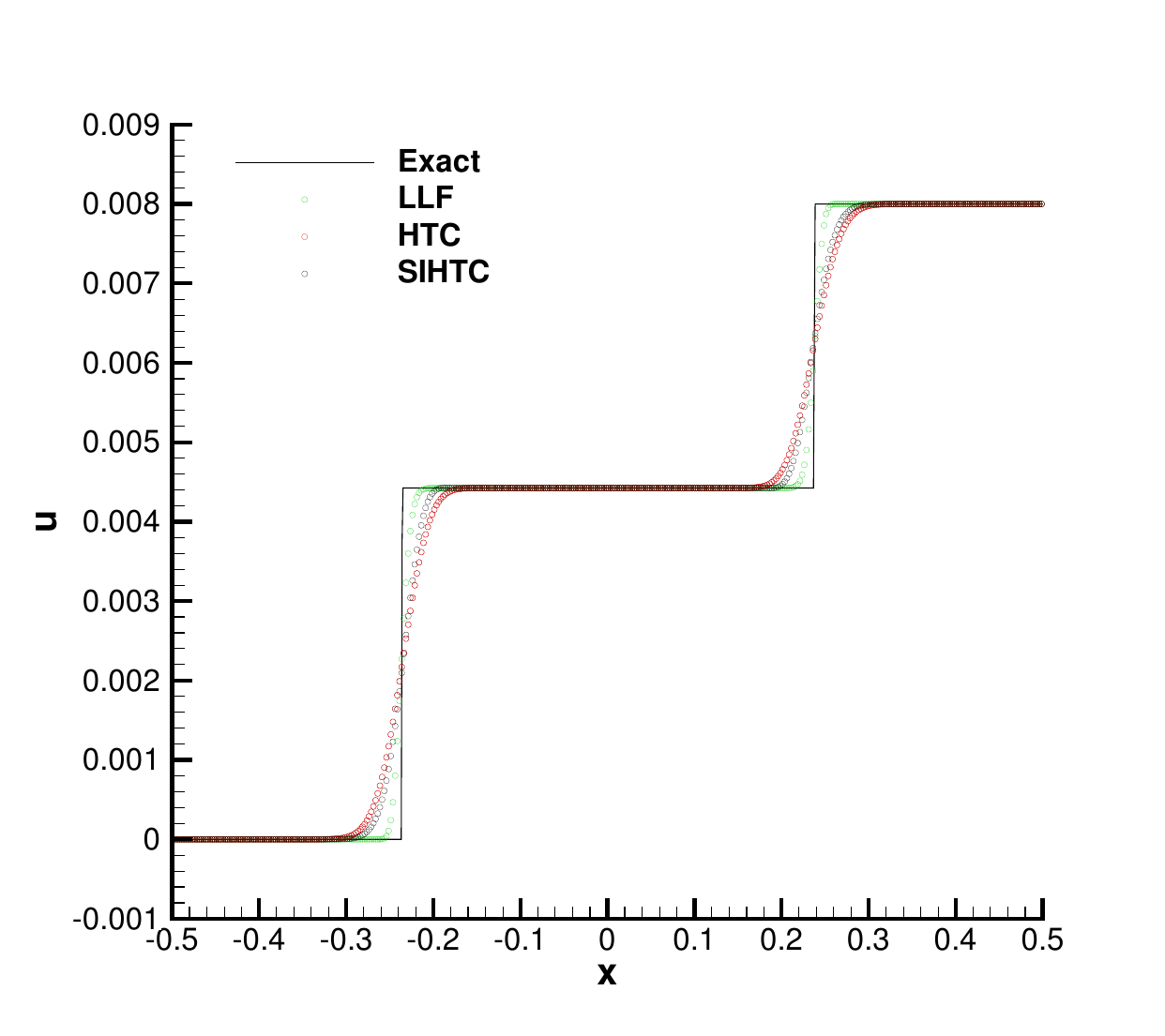}
		\caption{Velocity $u$}
	\end{subfigure}
	
	\vspace{1ex}
	
	\begin{subfigure}[b]{0.48\textwidth}
		\includegraphics[width=\textwidth]{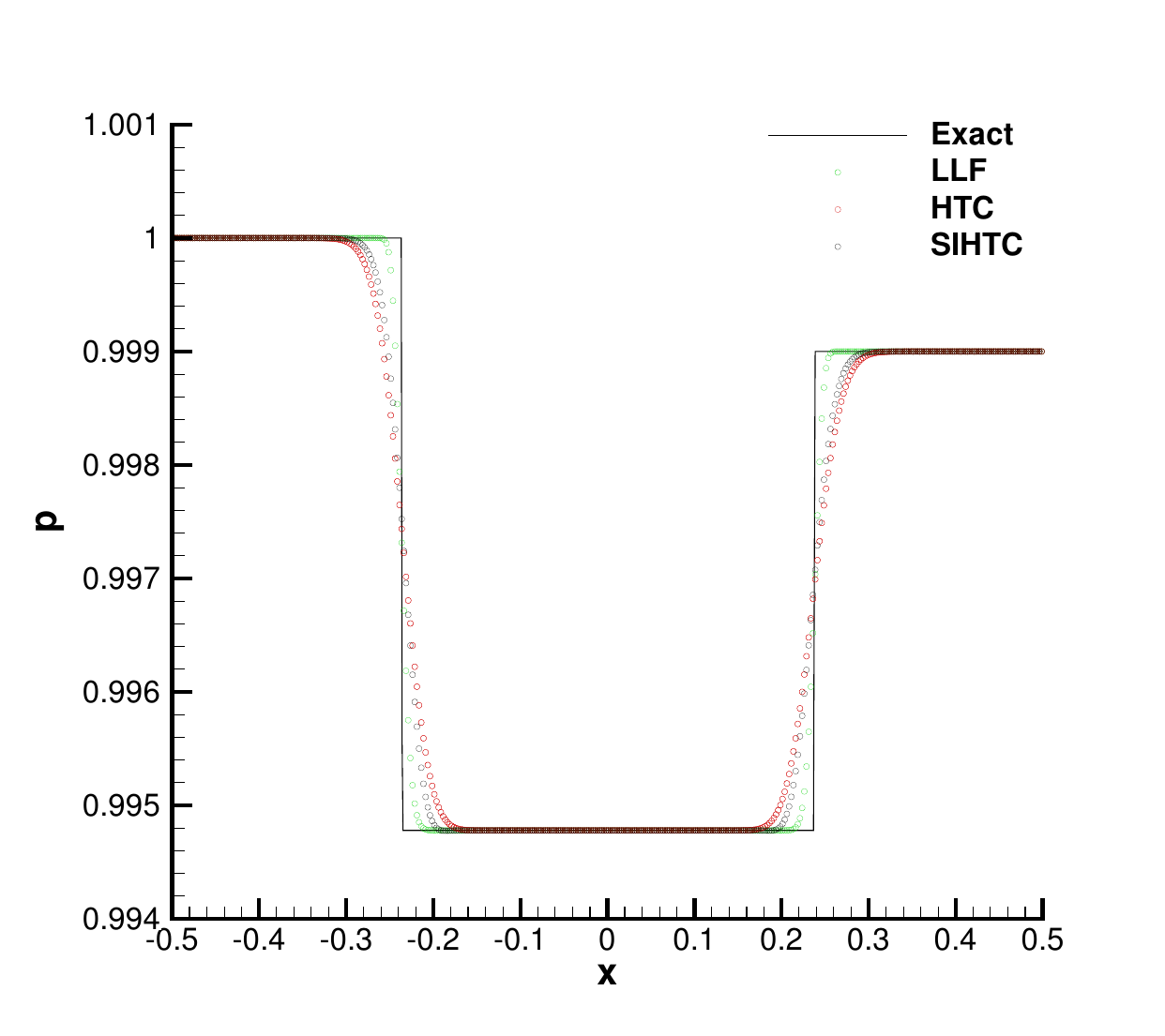}
		\caption{Pressure $p$}
	\end{subfigure}
	\hfill
	\begin{subfigure}[b]{0.48\textwidth}
		\includegraphics[width=\textwidth]{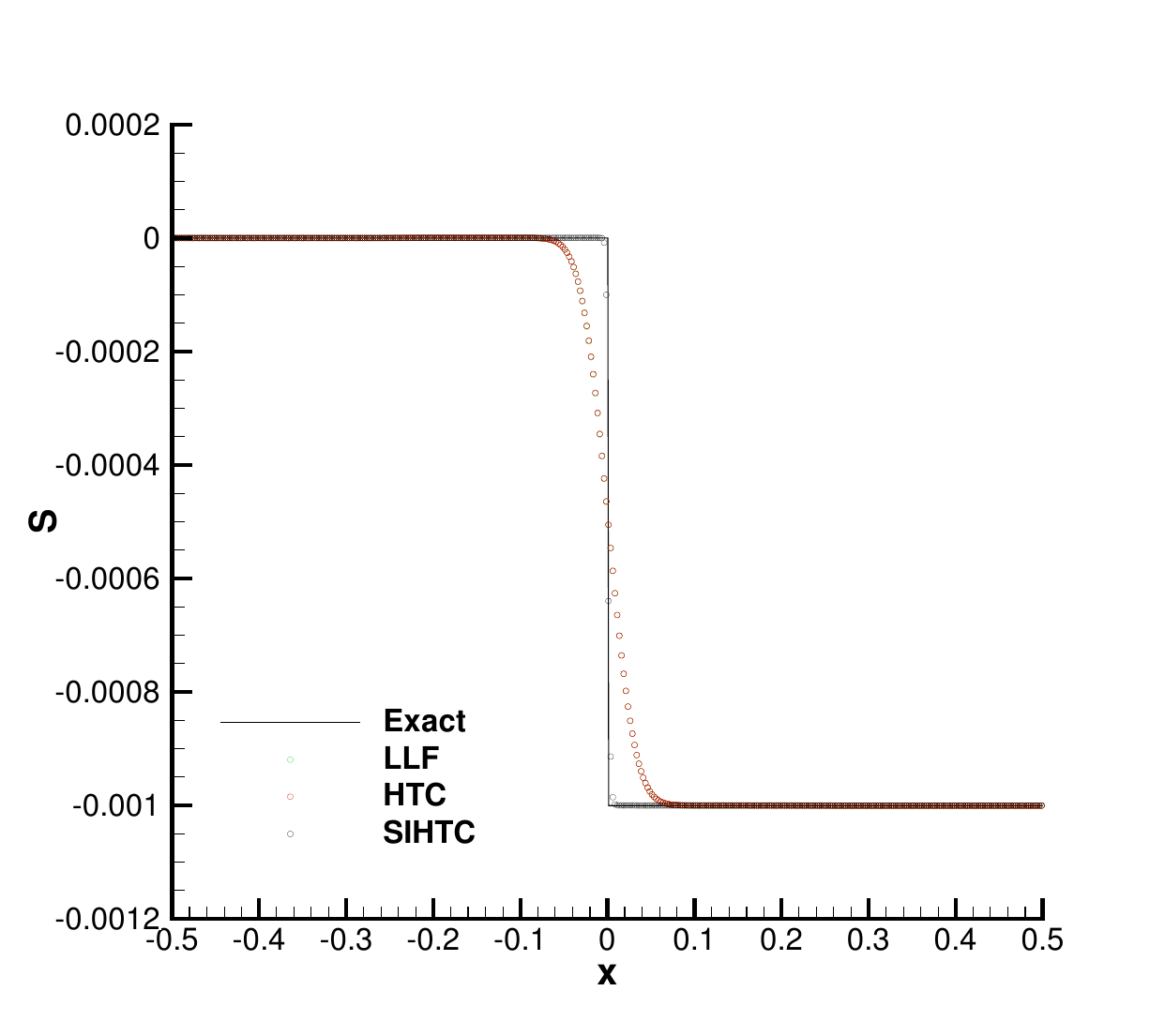}
		\caption{Entropy $S$}
	\end{subfigure}
	\caption{RP3 -- Low-Mach Riemann problem: numerical solutions at $t = 0.2$,
		$\mathrm{Ma} = 10^{-3}$, $N = 400$ cells.}
	\label{fig:rp3_lowmach}
\end{figure}

Figure~\ref{fig:rp3_lowmach} shows the numerical results at the final time $t = 0.2$ on a grid of $N=400$ cells and $\mathrm{CFL} = 0.5$. The exact solution to this Riemann problem consists of a steep left-moving rarefaction wave, a slowly moving contact discontinuity and a steep right-moving rarefaction. Since the Mach number regime is quite low, the acoustic waves travel at speed $c \approx \sqrt{\gamma p / \rho} \approx 1.18$, 
while the bulk flow velocity is $u = 8 \cdot 10^{-3}$. 
As a consequence, the pressure and density perturbations carried by 
the acoustic waves are of order $\mathcal{O}(\mathrm{Ma})$.
This test is particularly challenging even for low-Mach-number schemes since the contact wave moves slowly and is easily smeared out by the numerical viscosity, as is the case for the LLF and the HTC scheme.
However, the SIHTC scheme is able to capture the slow-moving wave accurately with only a few data points on the contact discontinuity, while the acoustic waves are, by design, dissipated. 
This is also reflected in the jump in the specific entropy, which is resolved accurately without any spurious oscillations. 

\subsubsection{Behavior of the global correction parameter}

\begin{figure}[htbp]
	\centering
	\begin{subfigure}[b]{0.48\textwidth}
		\includegraphics[width=\textwidth]{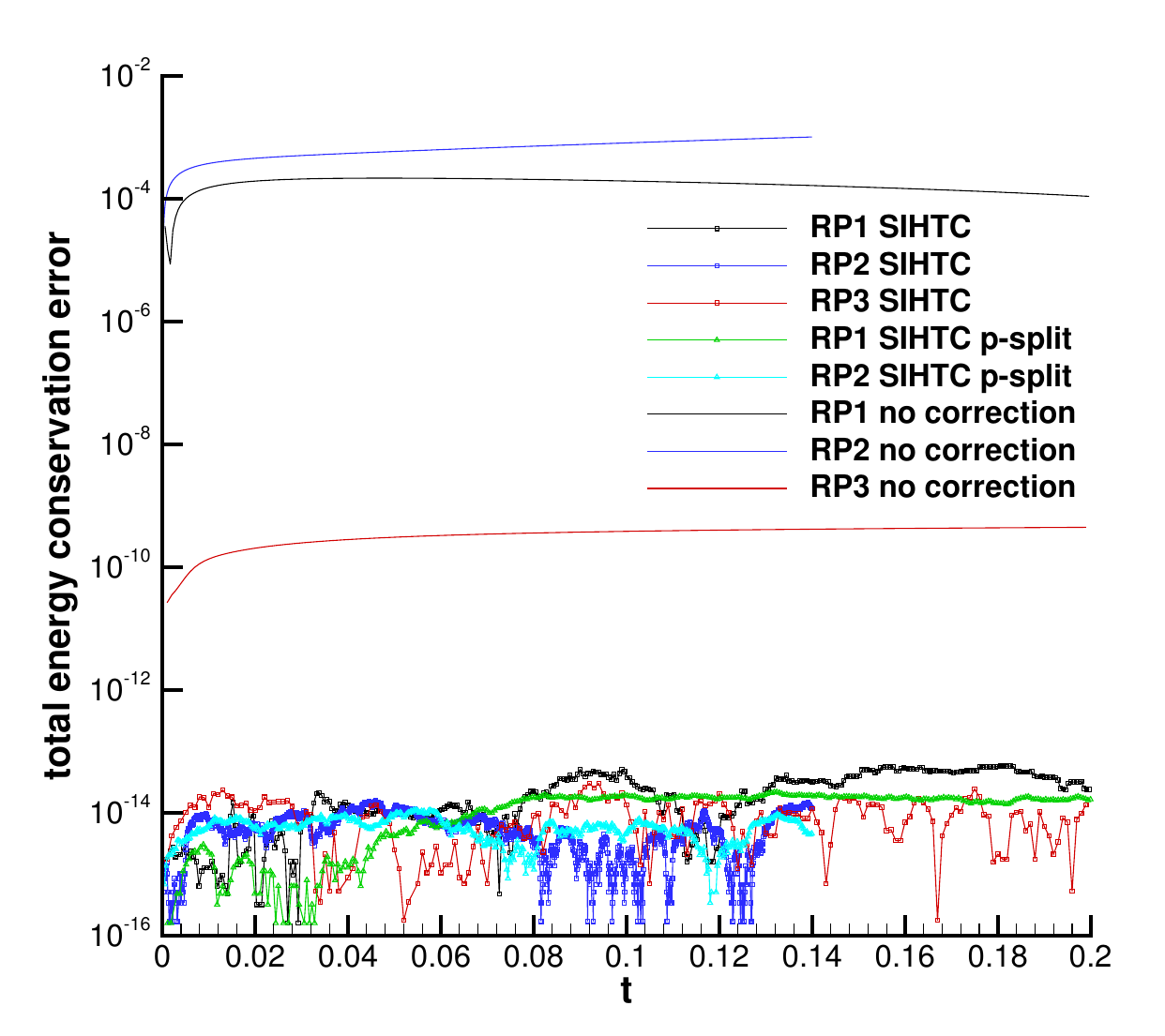}
		\caption{Energy conservation errors.}
		\label{fig:energy_error_abgrall}
	\end{subfigure}
	\hfill
	\begin{subfigure}[b]{0.48\textwidth}
		\includegraphics[width=\textwidth]{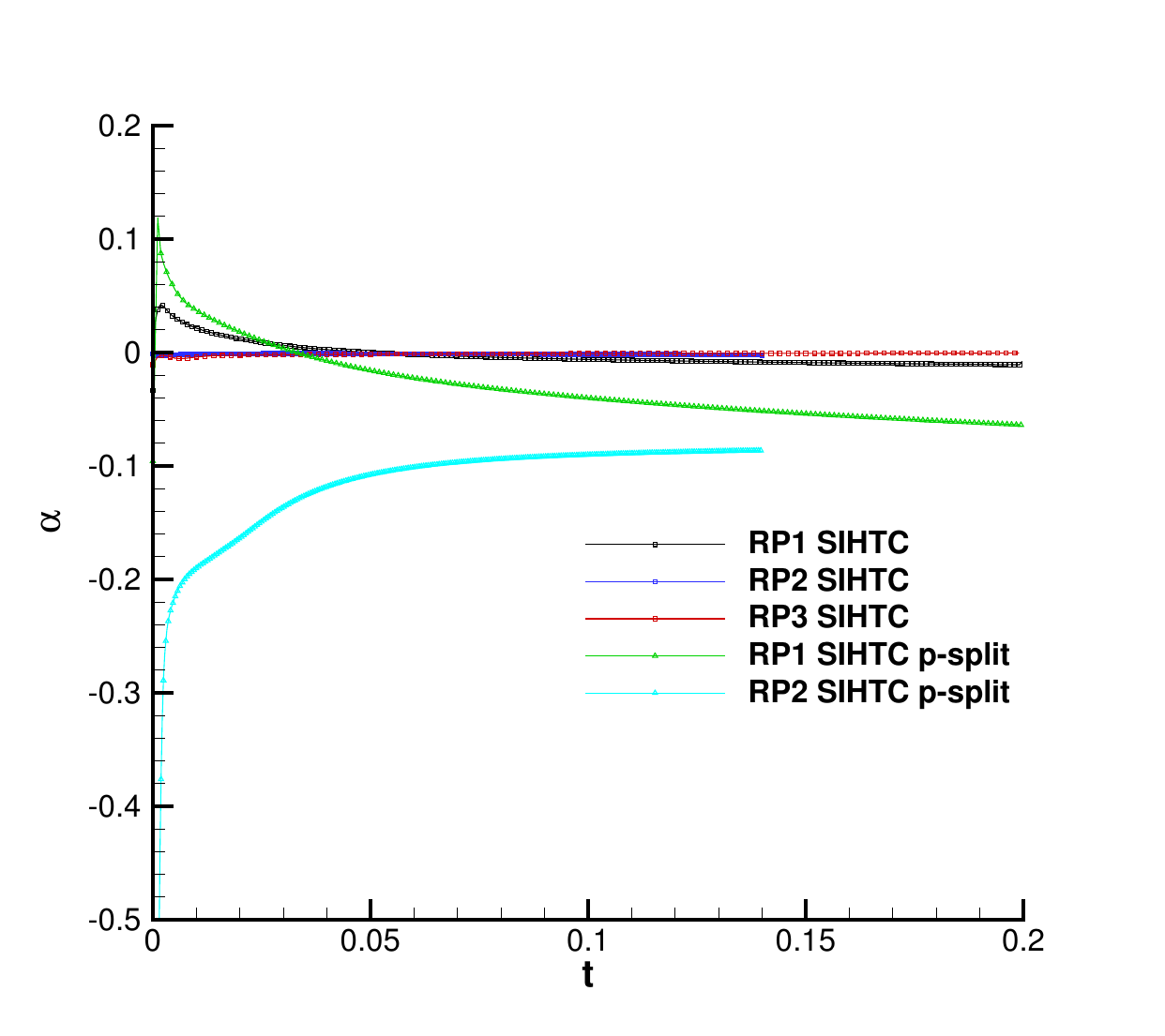}
		\caption{Evolution of the global energy correction factor $\alpha$.}
		\label{fig: a_abgrall}
	\end{subfigure}
	\caption{Effect of the global correction parameter $\alpha$ for the three Riemann problems at $N=400$.}
	\label{fig:alpha_energy}
\end{figure}

Finally, we would like to assess the behavior of the global correction parameter $\alpha$. 
Denoting by $\mathcal{E}^n = \Delta x \sum_i \mathcal{E}_i^n$ the total energy contained in the computational domain at time $t^n$, and since the boundary states $\mathbf{q}_L$ and $\mathbf{q}_R$ are constant in time, the energy conservation error reported in the following is the absolute deviation
\begin{equation}
	\label{eq: energy_error}
	\delta \mathcal{E}^n
	=
	\left|
	\mathcal{E}^n - \mathcal{E}^0
	+
	t^n
	\left(
	F(\mathbf{q}_R)
	-
	F(\mathbf{q}_L)
	\right)
	\right|,
\end{equation}
in which the energy transported through the two domain boundaries is taken into account.
In Figure~\ref{fig:energy_error_abgrall}, the effect of the global correction parameter $\alpha$ on total energy conservation for the proposed numerical semi-implicit schemes is shown. Omitting the correction leads to an energy error of order $\mathcal{O}(10^{-4})$ for RP1 and of order $\mathcal{O}(10^{-3})$ for RP2, whereas applying the global correction reduces the error to machine precision in both cases, as expected. Similarly, for RP3, without any correction, the global energy error is of order $\mathcal{O}(10^{-9})$ and again a global correction factor leads to machine precision errors. 
The logarithmic scale shows that the energy conservation error increases by several orders of magnitude when the correction procedure is omitted, which confirms that the global correction is essential to achieve energy conservation throughout the whole simulation.
 Figure~\ref{fig: a_abgrall} shows the magnitude of the correction parameter $\alpha$ at each time step for the different Riemann problems considered here. 
For the SIHTC scheme, the mean over time of the absolute value of the parameter $\alpha$ is of order $10^{-2}$ for RP1 and of order $10^{-3}$ for RP2 and RP3. For the SIHTC p-split scheme, where the correction acts on the conservative variables instead of the dual ones, it is of order $10^{-1}$ for RP1 and of order $10^{-2}$ for RP2.

\subsection{Two-dimensional test cases}
\label{sec:2D}

\subsubsection{Radial Sod shock tube}
\label{sec:radial_sod}

\begin{figure}[htpb]
	\centering
	\begin{subfigure}[b]{0.48\textwidth}
		\includegraphics[width=\textwidth]{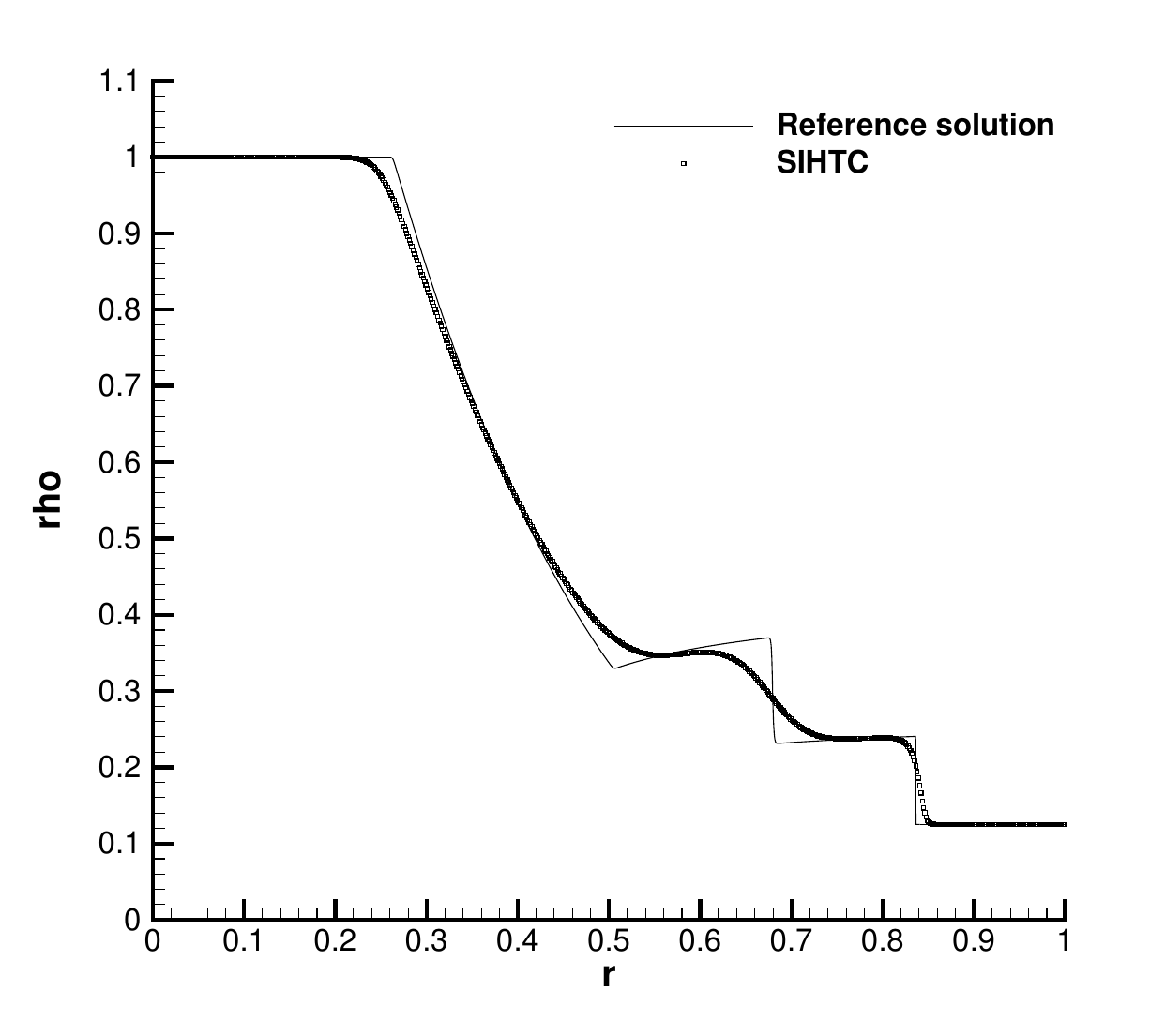}
		\caption{Density $\rho$}
	\end{subfigure}
	\hfill
	\begin{subfigure}[b]{0.48\textwidth}
		\includegraphics[width=\textwidth]{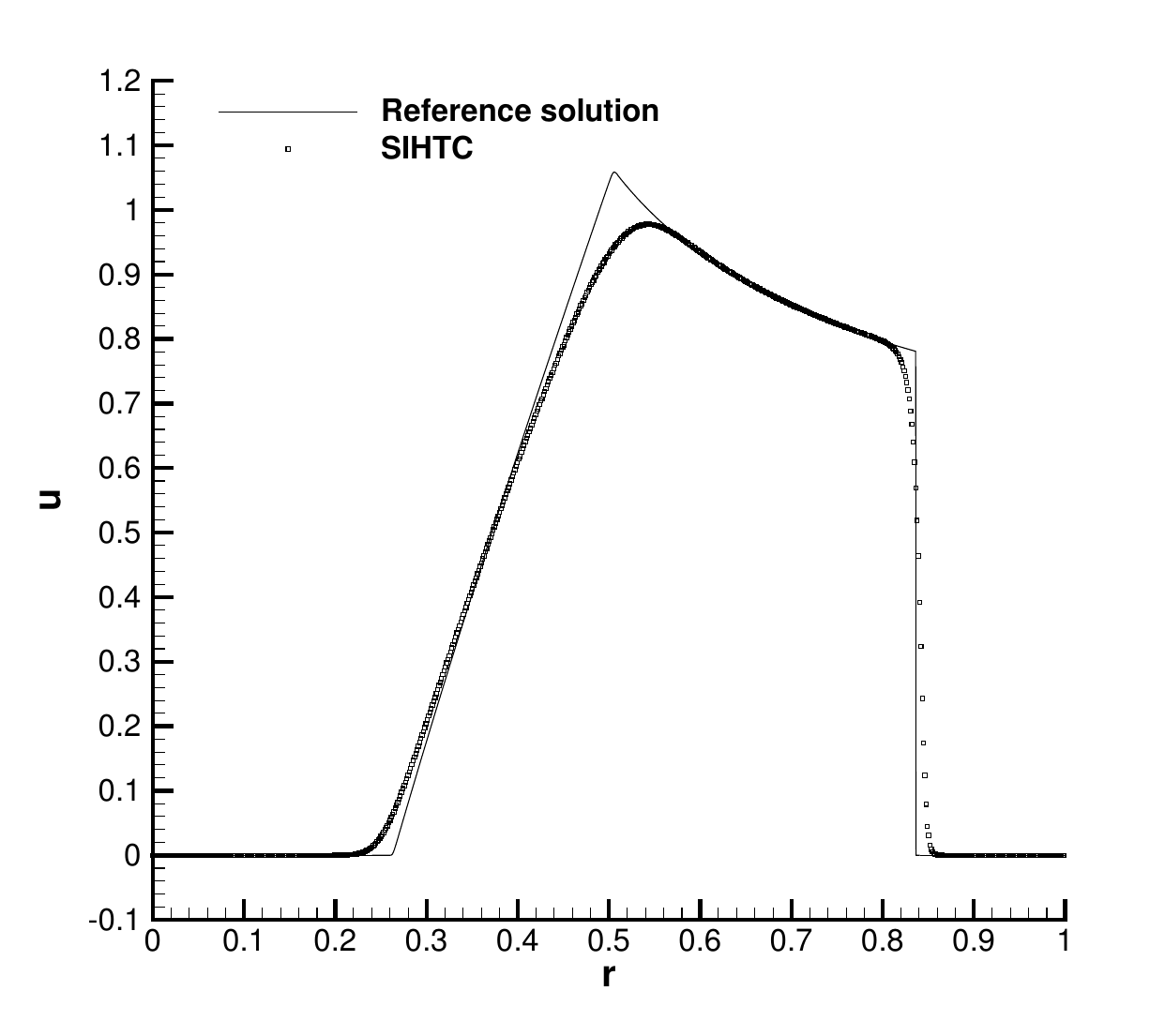}
		\caption{Radial velocity $u_r$}
	\end{subfigure}
	
	\vspace{1ex}
	
	\begin{subfigure}[b]{0.48\textwidth}
		\includegraphics[width=\textwidth]{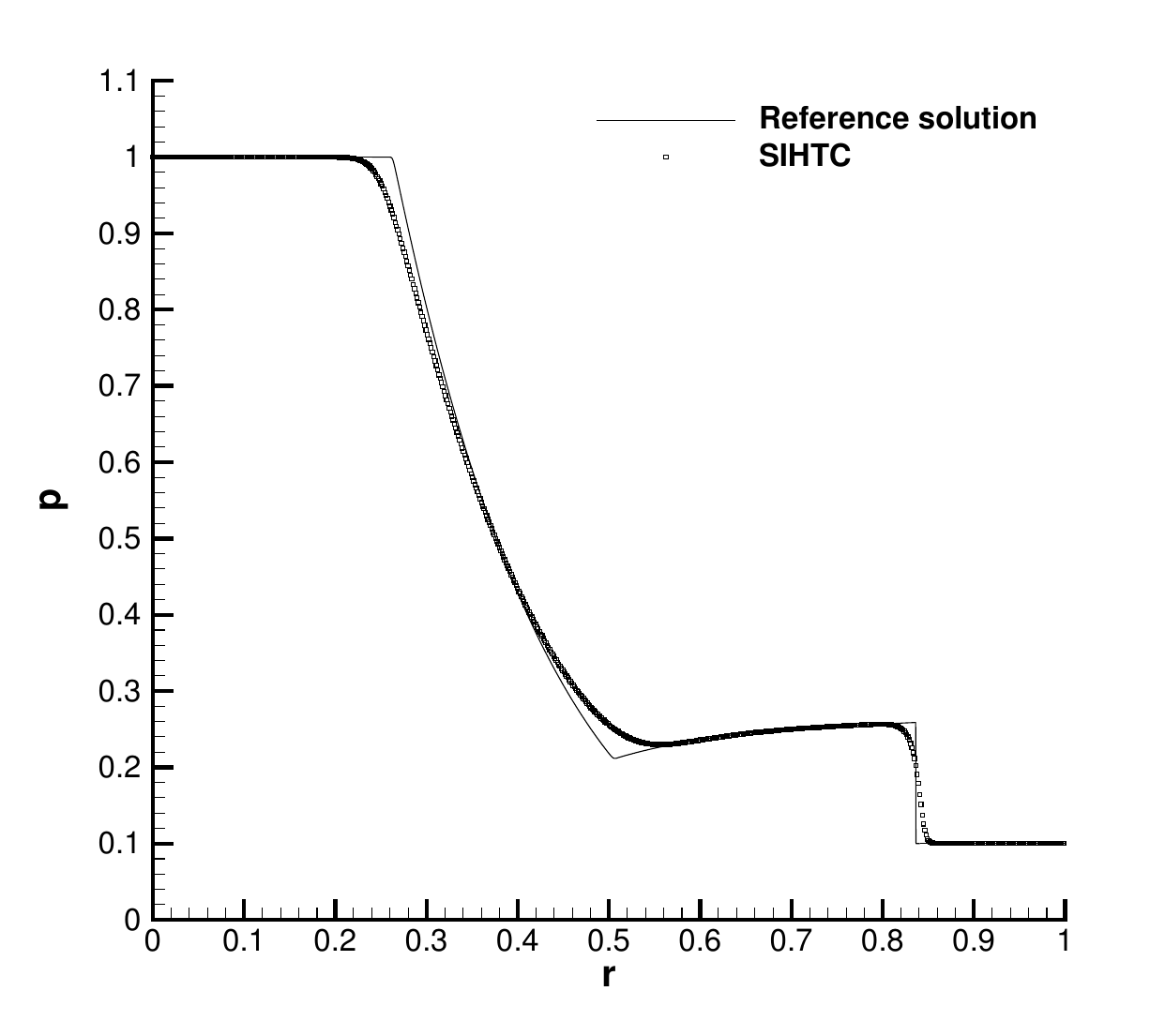}
		\caption{Pressure $p$}
	\end{subfigure}
	\hfill
	\begin{subfigure}[b]{0.48\textwidth}
		\includegraphics[width=\textwidth]{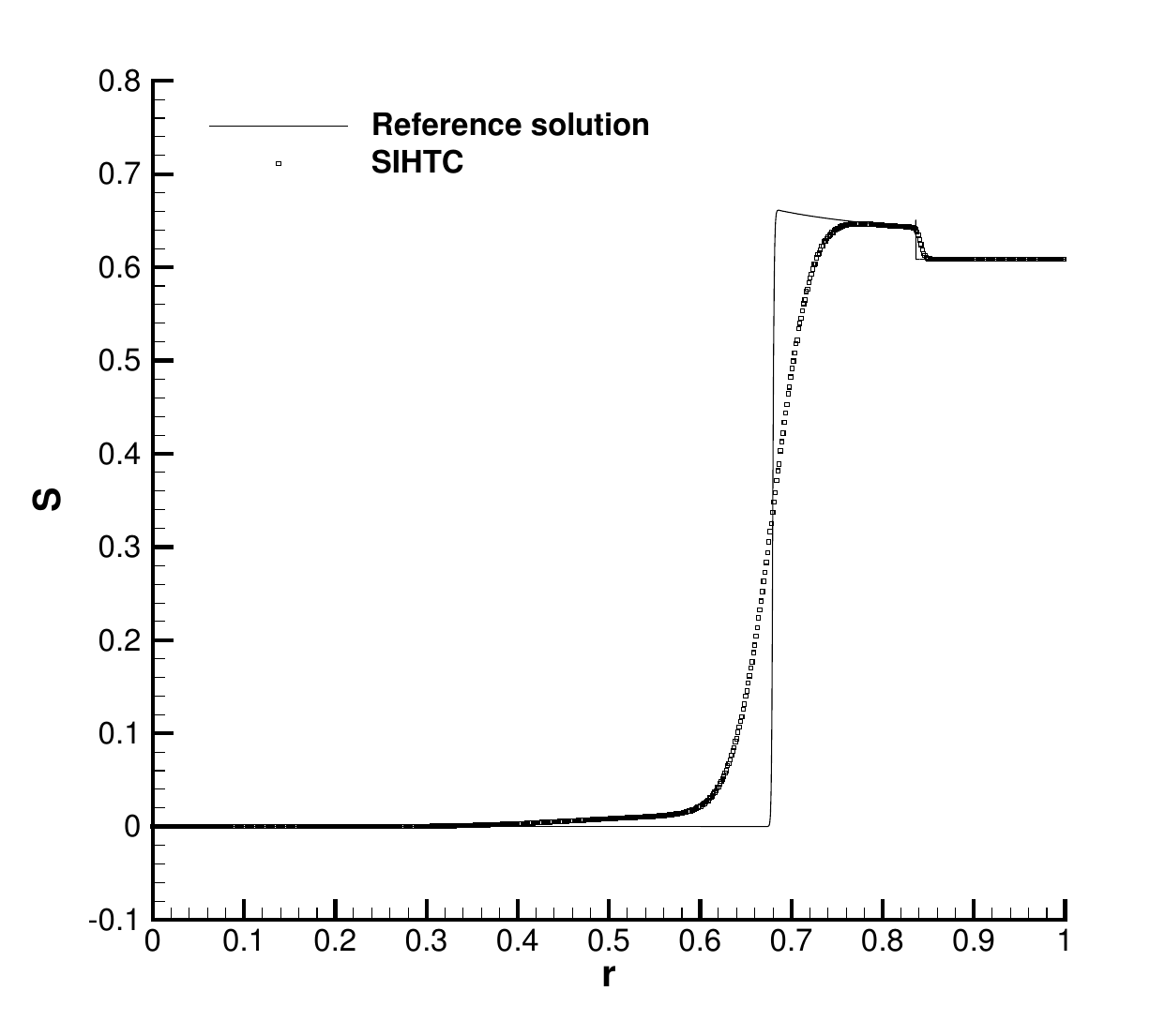}
		\caption{Entropy $S$}
	\end{subfigure}
		\caption{Radial Sod: scatter plot of the numerical solution along the radial 
		direction compared with the reference solution 
		(solid line) at $t = 0.2$ on $[-1,1]\times[-1,1]$ with $N = 500$ in each direction.}
	\label{fig:radial_sod_1D}
\end{figure}

\begin{figure}[htpb]
	\centering
	\begin{subfigure}[b]{0.48\textwidth}
		\includegraphics[width=\textwidth]{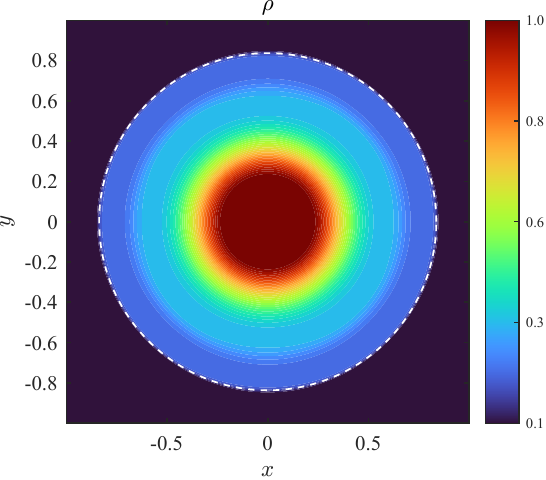}
		\caption{Density $\rho$}
	\end{subfigure}
	\hfill
	\begin{subfigure}[b]{0.48\textwidth}
		\includegraphics[width=\textwidth]{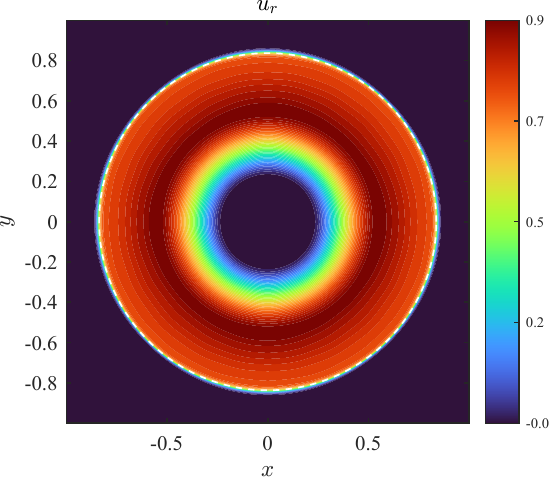}
		\caption{Radial velocity $u_r$}
	\end{subfigure}
	
	\vspace{1ex}
	
	\begin{subfigure}[b]{0.48\textwidth}
		\includegraphics[width=\textwidth]{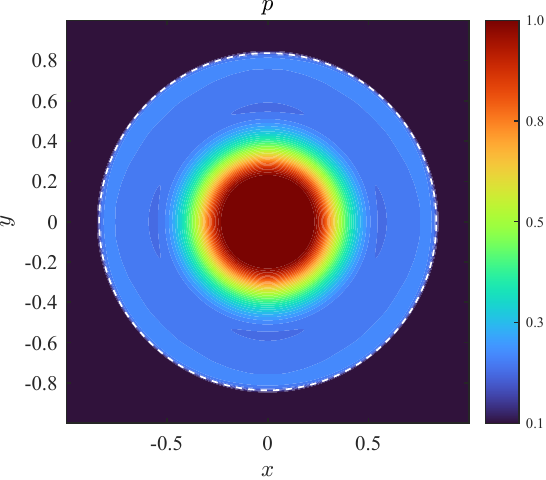}
		\caption{Pressure $p$}
	\end{subfigure}
	\hfill
	\begin{subfigure}[b]{0.48\textwidth}
		\includegraphics[width=\textwidth]{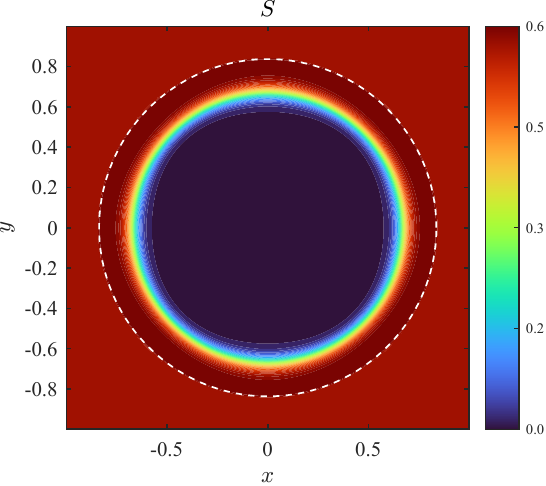}
		\caption{Entropy $S$}
	\end{subfigure}
	\caption{Radial Sod: 2D color maps at $t = 0.2$. The dotted line indicates the shock position of the reference solution.}
	\label{fig:radial_sod_2D}
\end{figure}

\begin{figure}[htpb]
	\centering
	\begin{subfigure}[b]{0.48\textwidth}
		\includegraphics[width=\textwidth]{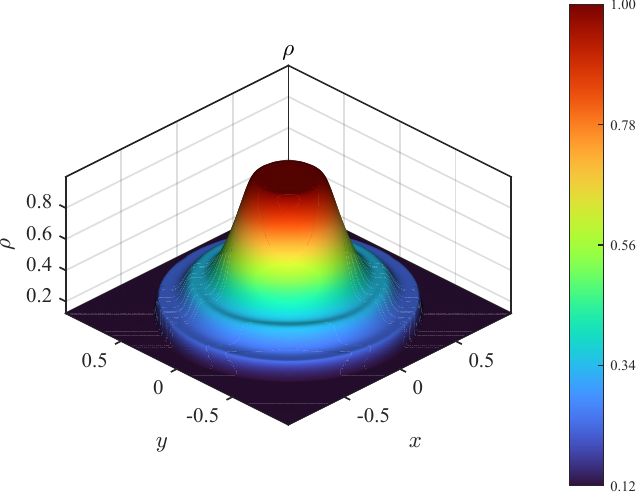}
		\caption{Density $\rho$}
	\end{subfigure}
	\hfill
	\begin{subfigure}[b]{0.48\textwidth}
		\includegraphics[width=\textwidth]{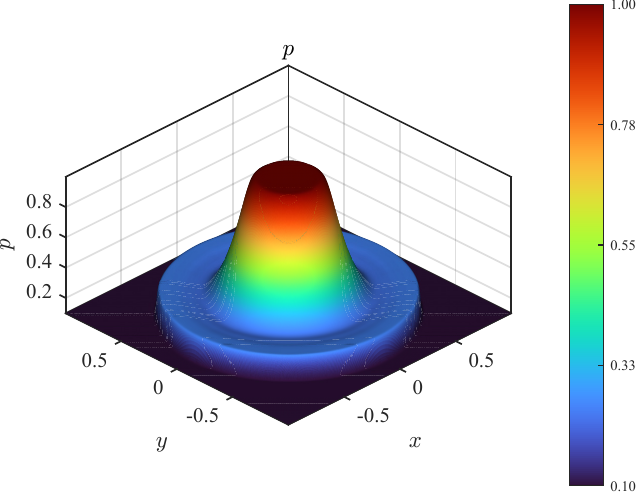}
		\caption{Pressure $p$}
	\end{subfigure}
	\caption{Radial Sod: 3D surface plots of density and pressure at $t = 0.2$.}
	\label{fig:radial_sod_3D}
\end{figure}

We extend the Sod problem to two space dimensions by imposing a cylindrically symmetric initial condition.
The computational domain is $\Omega = [-1,1]^2$, and the initial data depend only on the
radial coordinate $r = \sqrt{x^2+y^2}$:
\begin{equation}
	(\rho, u, v, p)(\mathbf{x},0) =
	\begin{cases}
		(1, 0, 0, 1)       & \text{if } r < 0.5,\\
		(0.125, 0, 0, 0.1) & \text{if } r \geq 0.5.
	\end{cases}
\end{equation}
The solution is computed with the SIHTC scheme up to $t_{\mathrm{end}} = 0.2$ with $\mathrm{CFL} = 0.5$
on a uniform Cartesian mesh of $N\times N$ cells with $N=500$.
The outer state is prescribed as a fixed Dirichlet boundary condition on all boundaries.
In Figure~\ref{fig:radial_sod_1D} we show the solution along the radial direction $r = \sqrt{x^2+y^2}$ for $y=0$, compared with a numerical reference solution. The reference solution has been computed by solving the radial Euler equations with geometric source terms on a sufficiently fine mesh, using a classical second-order MUSCL-Hancock scheme, see \cite{ToroBook} for details. Figure~\ref{fig:radial_sod_2D} shows the color maps of the density, velocity, pressure and entropy together with the shock position of the reference solution shown as a dotted line. The surface plots in three dimensions for the density and pressure are shown in Figure~\ref{fig:radial_sod_3D}.
All these results support the ability of the new SIHTC scheme to capture the position and amplitude of all discontinuities correctly even when solving for the entropy instead of the total energy. 

\subsubsection{Low-Mach-number limit: Taylor--Green vortex}
\label{sec:tgv}
\begin{figure}[htpb]
	\centering
	\begin{subfigure}[b]{0.48\textwidth}
		\includegraphics[width=\textwidth]{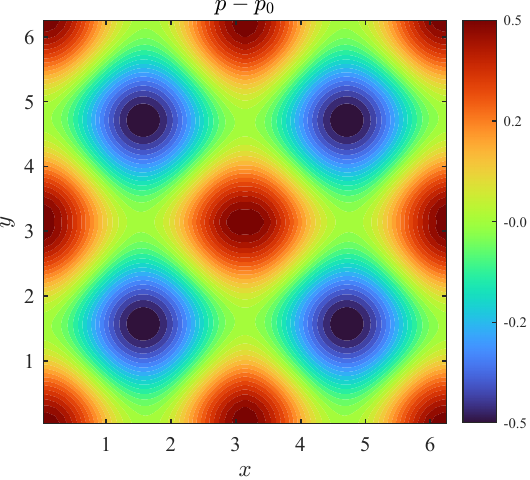}
		\caption{Pressure for Test 1}
	\end{subfigure}
	\hfill
	\begin{subfigure}[b]{0.48\textwidth}
		\includegraphics[width=\textwidth]{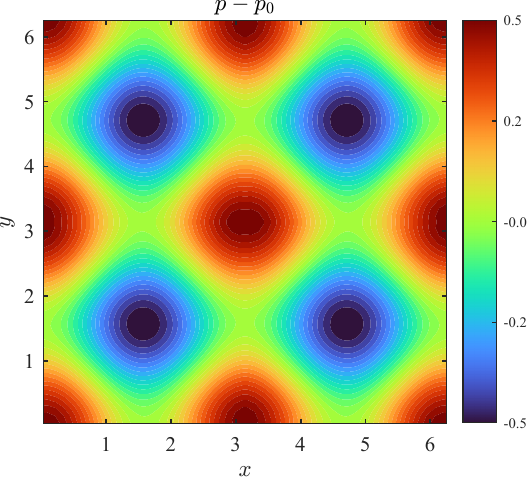}
		\caption{Pressure for Test 2}
	\end{subfigure}
	
	\vspace{1ex}
	
	\begin{subfigure}[b]{0.48\textwidth}
		\includegraphics[width=\textwidth]{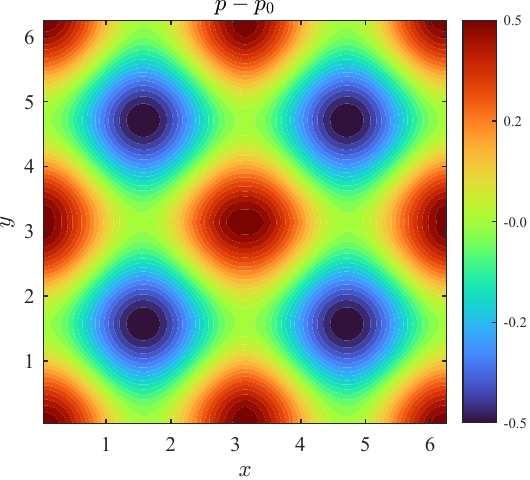}
		\caption{Pressure for Test 3}
	\end{subfigure}
	\hfill
	\begin{subfigure}[b]{0.48\textwidth}
		\includegraphics[width=\textwidth]{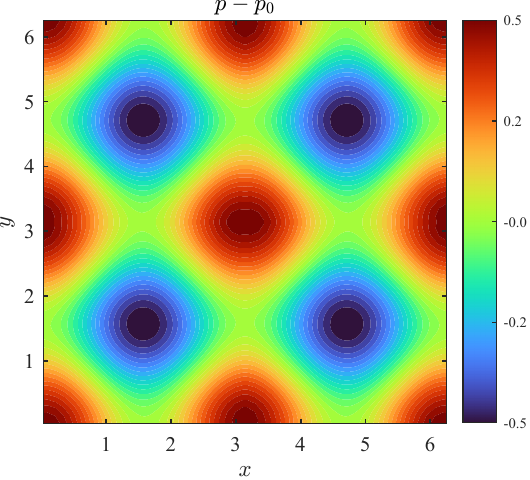}
		\caption{Pressure for Test 4}
	\end{subfigure}
	
	\caption{Taylor--Green vortex: 2D color maps of the pressure at $t = 0.2$.}
	\label{fig:tgv_p}
\end{figure}

\begin{figure}[htpb]
	\centering
	\begin{subfigure}[b]{0.48\textwidth}
		\includegraphics[width=\textwidth]{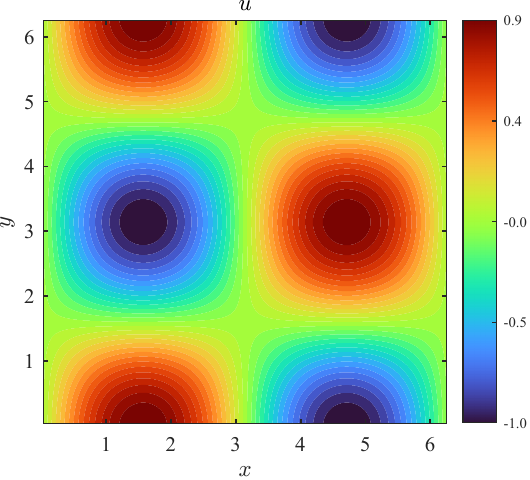}
		\caption{Velocity $u$ for Test 1}
	\end{subfigure}
	\hfill
	\begin{subfigure}[b]{0.48\textwidth}
		\includegraphics[width=\textwidth]{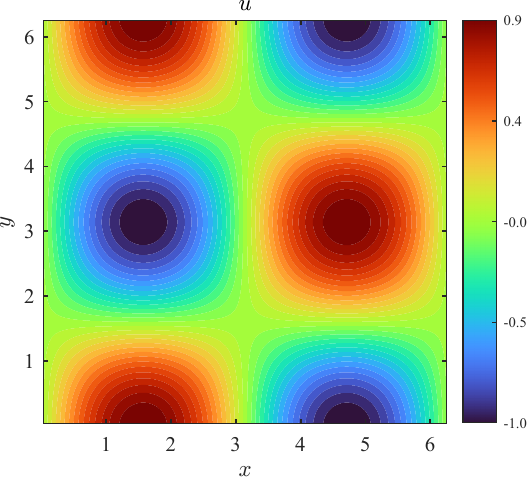}
		\caption{Velocity $u$ for Test 2}
	\end{subfigure}
	
	\vspace{1ex}
	
	\begin{subfigure}[b]{0.48\textwidth}
		\includegraphics[width=\textwidth]{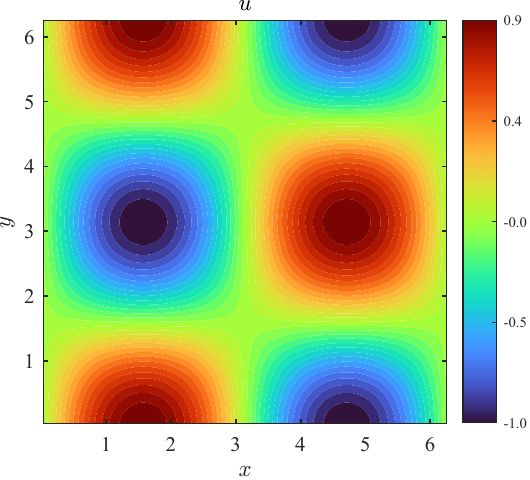}
		\caption{Velocity $u$ for Test 3 }
	\end{subfigure}
	\hfill
	\begin{subfigure}[b]{0.48\textwidth}
		\includegraphics[width=\textwidth]{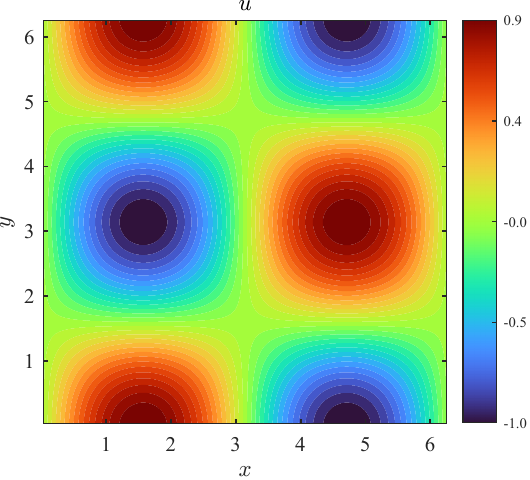}
		\caption{Velocity $u$ for Test 4}
	\end{subfigure}
	
	\caption{Taylor--Green vortex: 2D color maps of the velocity component $u$ at $t = 0.2$.}
	\label{fig:tgv_u}
\end{figure}

\begin{figure}[htpb]
	\centering
	\includegraphics[width=0.6\textwidth]{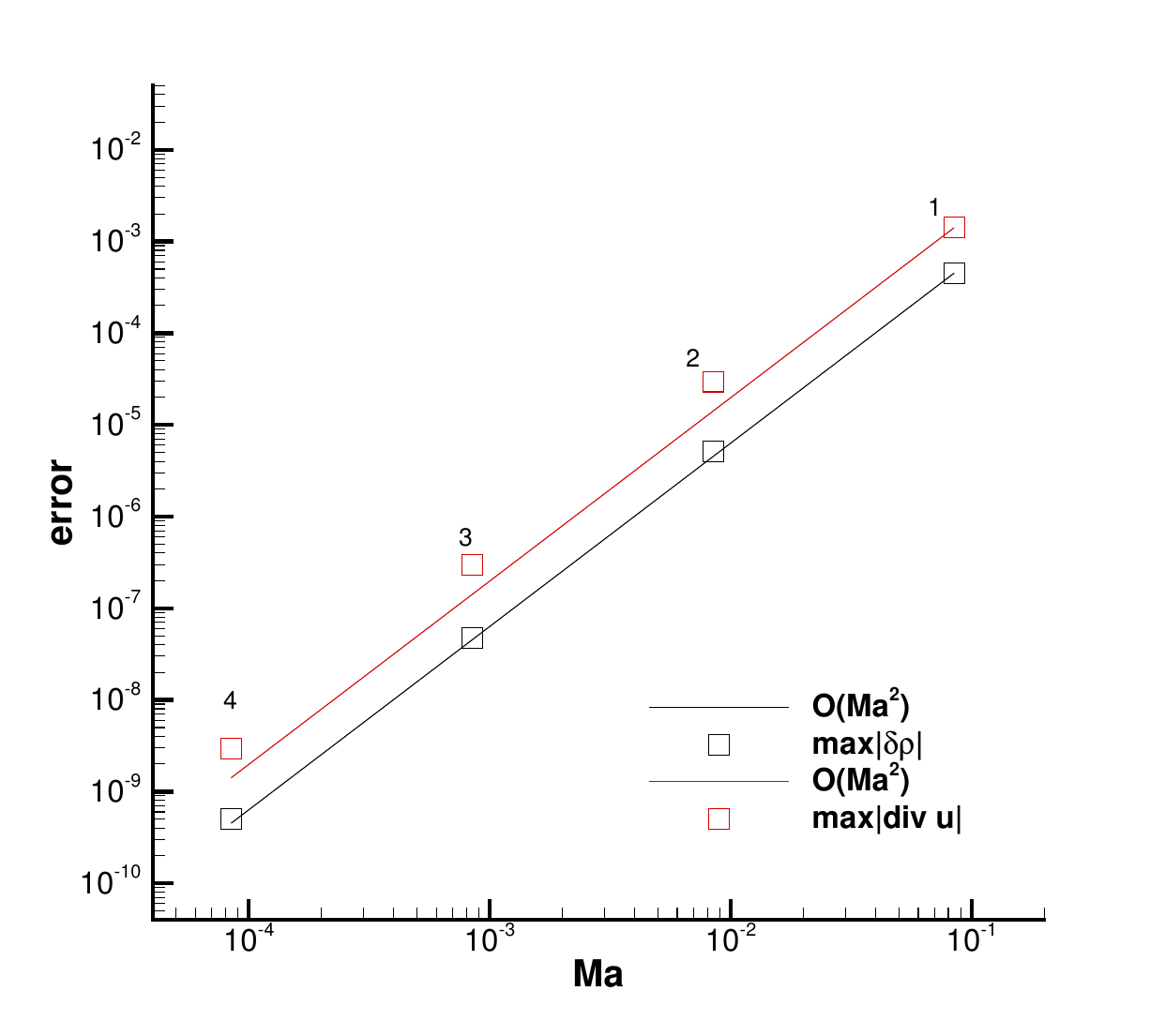}
	\caption{Taylor--Green vortex: maximum density fluctuation $\max|\delta\rho|$ and maximum velocity divergence $\max|\nabla \cdot \mathbf{u}|$ as functions of the Mach number.}
	\label{fig:tgv_drho_mach}
\end{figure}

To assess the behavior of the scheme in low-Mach-number regimes and to numerically verify the asymptotic-preserving property towards the incompressible Euler equations with constant density, which clearly emerges from the structure of the discrete pressure system \eqref{eq:pressure_system}, \eqref{eqn.J.def} and \eqref{eqn.D.def},   
we consider the Taylor--Green vortex~\cite{TaylorGreen1937,SIGPR,dumbser_2026_a}.
The computational domain is $\Omega = [0,2\pi]^2$ with periodic boundary conditions.
The initial condition is given by
\begin{equation}
	\begin{aligned}
	&\rho(x,y,0) = \rho_0, \qquad 
	p(x,y,0)    = p_0 + \frac{1}{4}\left[\cos(2 x)+\cos(2 y)\right], \\
	&u(x,y,0) =  \sin(x)\cos(y), \qquad 
	v(x,y,0) = -\cos(x)\sin(y),
	\end{aligned}
\end{equation}
where $\rho_0 = 1$ and $p_0 \in \{10^{2}, 10^{4}, 10^{6}, 10^{8}\}$. By varying the initial pressure magnitude, we obtain different Mach number regimes as reported in Table~\ref{tab:mach_p0}. We denote by Test $r$, with $r=1,\ldots,4$, the numerical simulation with initial data $p_0 = 10^{2r}$ as in Table~\ref{tab:mach_p0}.
All the numerical results have been obtained on a fixed uniform Cartesian mesh $N \times N$ with $N=100$ and $\mathrm{CFL} = 0.5$. For all the Mach numbers considered in this benchmark, the blending coefficient remains equal to one throughout all the simulations. Consequently, the p-split version reduces to the original SIHTC formulation and no distinction between the two schemes is necessary.
\begin{table}[htpb]
	\centering
	\caption{Correspondence between the reference pressure $p_0$ and the 
		Mach number $\mathrm{Ma} = 1/\sqrt{\gamma p_0 / \rho_0}$
		for $\gamma = 1.4$, $\rho_0 = 1$.}
	\label{tab:mach_p0}
	\renewcommand{\arraystretch}{1.25}
\begin{tabular}{ccc}
	\hline
	Test & $p_0$ & $\mathrm{Ma}$ \\
	\hline
	1 & $10^{2}$ & $8.452 \cdot 10^{-2}$ \\
	2 & $10^{4}$ & $8.452 \cdot 10^{-3}$ \\
	3 & $10^{6}$ & $8.452 \cdot 10^{-4}$ \\
	4 & $10^{8}$ & $8.452 \cdot 10^{-5}$ \\
	\hline
\end{tabular}
\end{table}
In Table~\ref{tab:tgv_drho} and Figure~\ref{fig:tgv_drho_mach}, we report both the maximum density fluctuation $\max|\delta\rho|$ and the maximum velocity divergence $\max|\nabla\cdot\mathbf{u}|$ as functions of the Mach number. Both quantities exhibit the expected $\mathcal{O}(\mathrm{Ma}^2)$ scaling, providing numerical evidence of the asymptotic-preserving behavior of the proposed scheme in the incompressible limit.
Moreover, the velocity fields and pressure retain the correct vortex structure as can be seen from Figures~\ref{fig:tgv_p} and \ref{fig:tgv_u}.
For each simulation, we monitor the behavior of the global Abgrall corrector $\alpha$ as the Mach number is progressively decreased, in terms of the mean and the maximum absolute value of $\alpha$ over the entire time history of the simulation. 
For the highest Mach number considered, $\mathrm{Ma} = 8.452 \cdot 10^{-2}$, the correction factor attains a small but non-zero value, with $\operatorname{mean}|\alpha| = 4.829 \cdot 10^{-7}$ and $\max|\alpha| = 8.853 \cdot 10^{-7}$. For all other cases, as the Mach number is further reduced by successive orders of magnitude, the mean and the maximum absolute value of $\alpha$ are at machine precision.
This is due to the fact that the predictor scheme is designed for low-Mach-number flows and thus the errors produced are already lower than the imposed criterion on energy conservation.

The vanishing of $\alpha$ at low Mach numbers therefore shows that the corrector is inactive when compressibility effects are negligible, ensuring a smooth and consistent transition between the compressible and incompressible formulations.

\begin{table}[htpb]
	\centering
	\caption{Taylor--Green vortex: maximum density fluctuation $\max|\delta\rho|$ and maximum velocity divergence $\max|\nabla \cdot \mathbf{u}|$ as functions of the Mach number at $t = 0.2$.}
	\label{tab:tgv_drho}
	\renewcommand{\arraystretch}{1.25}
\begin{tabular}{cc|cc|cc}
	\hline
	Test & $\mathrm{Ma}$ & $\max|\delta\rho|$ & rate & $\max|\nabla \cdot \mathbf{u}|$ & rate \\
	\hline
	1 & $8.452 \cdot 10^{-2}$ & $4.520 \cdot 10^{-4}$ & --   & $1.414 \cdot 10^{-3}$ & --       \\
	2 & $8.452 \cdot 10^{-3}$ & $5.170 \cdot 10^{-6}$ & 1.94 & $2.959 \cdot 10^{-5}$ & 1.68 \\
	3 & $8.452 \cdot 10^{-4}$ & $4.740 \cdot 10^{-8}$ & 2.04 & $2.958 \cdot 10^{-7}$ & 2.00 \\
	4 & $8.452 \cdot 10^{-5}$ & $5.018 \cdot 10^{-10}$ & 1.97 & $2.956 \cdot 10^{-9}$ & 2.00 \\
	\hline
\end{tabular}
\end{table}


\section{Conclusions}
\label{sec: conclusions}
We have presented a new semi-implicit scheme for the compressible Euler equations which is based on the entropy rather than on the total energy as evolution variable. 
Since the Euler equations can be written in hyperbolic and thermodynamically compatible (HTC) form, a semi-implicit HTC scheme was developed, which is stable under a material time step restriction. 
It is therefore especially suited for applications in low-Mach-number regimes. 
This has been achieved by splitting the Euler flux into a convective and a pressure subsystem, where the former is treated explicitly and the latter implicitly. Instead of solving the entire nonlinear implicit pressure-momentum subsystem directly, it was rewritten with the aid of the Schur complement into a single mildly nonlinear system for the pressure, making use of a staggering of the momentum. 
To ensure the correct shock strength and position, an entropy production term was introduced, taking into account the spatial and the temporal viscosity of the scheme, which yields a compatible numerical viscosity for the total energy equation that has to hold as a consequence. 
To also ensure compatibility with the total energy fluxes and finally global energy conservation, a global correction factor, the so-called global Abgrall-type correction factor, has been introduced into the numerical flux, thus ensuring global total energy conservation up to machine precision.
The properties of the numerical scheme have been numerically assessed in a series of test cases. In particular, we have assessed the consistency with Riemann problems in the compressible and in the nearly incompressible regime, as well as the asymptotic consistency and accuracy as the Mach number goes to zero. 
By means of the Taylor--Green vortex it has been verified that the numerical solution is indeed consistent with the solution of the incompressible Euler equations. 
Moreover, the density fluctuations and the divergence of the velocity converge numerically with the Mach number squared, thus achieving the expected asymptotic accuracy.

In addition, a p-split variant based on an IMEX pressure splitting has been proposed. The numerical experiments show that it retains the accuracy of the original semi-implicit HTC formulation in the compressible regime while allowing significantly larger admissible CFL numbers, thereby providing a robust transition between low and moderate Mach number flows.
This work demonstrates that it is indeed possible to construct semi-implicit numerical schemes based on the entropy formulation. 
In the case of the Euler equations it is a matter of choice whether one solves for the total energy or for the entropy. 
However, there are systems such as two-phase flows where two entropy inequalities, one for each phase, stand opposite a single conservation law for the total energy, and where the HTC formalism turns out to be particularly useful, see e.g.~\cite{thomann_2023_thermodynamically}. 
In this case, it is indispensable to solve for the entropies and conserve the energy as a consequence. The global correction does not guarantee the positivity of density and pressure by construction; a provably positivity-preserving variant is left for future work.
This work can thus be seen as a first step towards a semi-implicit HTC scheme for compressible two-phase flows and for the MHD equations, which will be the subject of future work. We also plan to extend this approach to unstructured triangular meshes using the compatible semi-implicit discretization recently forwarded in \cite{FourSplitUS,Bernardelli1}.

\vspace{-2mm}

\section*{Acknowledgments}

M.D. was funded by the Fondazione Caritro via the project SOPHOS and by the European Research Council (ERC) under the European Union's Horizon Europe research and innovation programme via the project SOPHOS, grant agreement no. ERC-ADG-2025-101265878-SOPHOS. Views and opinions expressed are however those of the authors only and do not necessarily reflect those of the European Union or the European Research Council Executive Agency. Neither the European Union nor the granting authority can be held responsible for them.
A.T. acknowledges the financial support of the Agence Nationale de la Recherche (ANR) via the project DELFIN, project no. ANR-25-CE46-7729.
M.D. and A.T. were also financially supported by the Italian Ministry of University
and Research (MUR) via the Departments of Excellence Initiative 2018--2027 attributed to DICAM of the University of Trento (grant L. 232/2016).
G.P. acknowledges financial support of the Sapienza Research fund Ateneo 2025.  

All authors are members of the Gruppo Nazionale per il Calcolo Scientifico dell'Istituto Nazionale di Alta Matematica (GNCS-INdAM).

\vspace{-2mm}

\section*{Conflict of interest}

The authors declare that they have no conflict of interest. 

\vspace{-2mm}

\section*{Data availability}

The data can be obtained from the authors on reasonable request. 

\bibliographystyle{spmpsci}      
\bibliography{biblio}

\end{document}